\documentclass[11pt]{article}
\usepackage{comment}
\usepackage[a4paper,twoside]{geometry}
\usepackage[utf8]{inputenc}

\usepackage{amssymb}
\usepackage{amsmath}
\usepackage{amsthm}

\usepackage{xcolor}
\usepackage{color}
\usepackage{tabularx}
\usepackage{ragged2e}
\usepackage{multicol}

\usepackage{thmtools}
\usepackage{thm-restate}

\includecomment{show}

\newcommand{\anonymousvariant}[2]{#2}
\begin{show}
\renewcommand{\anonymousvariant}[2]{{#1}}
\end{show}

\def\cal{\mathcal}

\def\Seq{\textup{\textsc{Seq}}}  
\def\Cyc{\textup{\textsc{Cyc}}}  
\def\Set{\textup{\textsc{Set}}}

\def\Exp{\textup{\textsc{Exp}}}

\def\bs#1{{\boldsymbol #1}}   
\def\bc#1{\boldsymbol{{\mathcal #1}}}

\def\bch#1{\boldsymbol{{\mathcal{\hat{#1}}}}}
\def\bct#1{\boldsymbol{{\mathcal{\tilde{#1}}}}}

\def\cal{\mathcal}
\def\cZ{{\cal Z}}
\def\cY{{\cal Y}}
\def\cG{{\cal G}}
\def\ini{{\bs U}}
\def\bpartial{\bs\partial}

\def\N{{\mathbb N}}
\def\Id{{\mathbf{Id}}} 

\def\val{\operatorname{val}}
\def\LT{\operatorname{LT}}
\def\al{\operatorname{alg-log}}

\def\blueroot{\cal{T}_{b}}
\def\redroot{\cal{T}_{r}}
\def\greenroot{\cal{T}_{g}}
\def\cblue{\cal{B}}
\def\cred{\cal{R}}
\def\cgreen{\cal{G}}
\def\cforest{\cal{F}}

\def\ds{\displaystyle } 

\newcounter{example}
\newenvironment{example}{\refstepcounter{example} 
\smallskip
\noindent {\bf Example~\theexample.}\xspace}   
{\smallskip}

\usepackage[ruled,titlenotnumbered,vlined,linesnumbered]{algorithm2e}
\definecolor{darkteal}{RGB}{0,70,90}

\SetCommentSty{mycommfont}
\SetKwComment{tcp}{$\lhd$~}{}
\SetKwComment{tcc}{$\triangledown$~}{}
\SetKwComment{CaseComment}{\color{darkteal}$\lhd$~}{} 

\SetNlSkip{0.5em}
\IncMargin{0.5em}

\usepackage{float}

\SetKwInput{Input}{Input}
\SetKwInput{Output}{Output}
\SetKw{KwOf}{of}
\DontPrintSemicolon

\newenvironment{method}[1][]
{ 
  \begin{algorithm}[#1]}
{ \end{algorithm}
  }
  
\newenvironment{numalgo}[1][]
{ 
  \begin{algorithm}[#1]}
{ \end{algorithm}
  }

\usepackage[colorlinks=true,citecolor=blue,linkcolor=violet]{hyperref}

\usepackage[capitalize]{cleveref} 

\SetAlCapSty{xAlCapSty}

\newtheorem{theorem}{Theorem}[section]
\newtheorem{lemma}[theorem]{Lemma}
\newtheorem{proposition}[theorem]{Proposition}
\newtheorem{corollary}[theorem]{Corollary}

\newtheorem{conjecture}[theorem]{Conjecture}
\theoremstyle{definition}
\newtheorem{definition}[theorem]{Definition}
\newtheorem*{notation}{Notation}
\theoremstyle{remark}
\newtheorem{remark}[theorem]{Remark}

\usepackage {tikz}
\usetikzlibrary {positioning}
\definecolor {processblue}{cmyk}{0.96,0.5,0,0}
\definecolor {mygreen}{cmyk}{0.96,0,0.5,0}

\title{Effective Asymptotics of Combinatorial Systems} 

\author{\anonymousvariant{Carine Pivoteau}{}\and \anonymousvariant{Bruno Salvy}{}}
\begin{document}

\maketitle


\begin{abstract}
Analytic combinatorics studies asymptotic properties of families of combinatorial objects using complex analysis on their generating functions. In their reference book on the subject, Flajolet and Sedgewick describe a general approach that allows one to derive precise asymptotic expansions starting from systems of combinatorial equations. In the situation where the combinatorial system involves only cartesian products and disjoint unions, the generating functions satisfy polynomial systems with positivity constraints for which many results and algorithms are known. We extend these results to the general situation. This produces an almost complete algorithmic chain going from combinatorial systems to asymptotic expansions. Thus, it is possible to compute asymptotic expansions of all generating functions produced by the symbolic method of Flajolet and Sedgewick when they have algebraic-logarithmic singularities (which can be decided), under the assumption that Schanuel's conjecture from number theory holds. That conjecture is not needed for systems that do not involve the constructions of sets and cycles.
\end{abstract}

\medskip

\anonymousvariant{
\begin{quote}\hfill \begin{minipage}{10.5cm}\em 
    This article is dedicated to the memory of Philippe Flajolet, who shared with us his vision of the interaction between analytic combinatorics and computer algebra.\end{minipage}
\end{quote}
}{}

\begin{small} 
\tableofcontents
\end{small}

\section{Introduction}

In their book \emph{Analytic Combinatorics}~\cite{FlajoletSedgewick2009}, Flajolet and Sedgewick describe a general approach that starts from a combinatorial description, translates this description into equations satisfied by generating functions, views these generating functions as analytic functions and exploits their singular behavior to deduce asymptotic properties of the combinatorial objects when their size becomes large. The aim of this article is to develop computational tools that automate large parts of this approach. 
We first outline the main steps and the computational issues that they raise.

\paragraph{Symbolic Method and Well-Founded Systems}
The first step of the approach is called the `Symbolic Method'.\footnote{{\sl ``The symbolic method is the `combinatorics' in analytic combinatorics''}~\cite[p.~92]{FlajoletSedgewick2009}.} It consists of a precise definition of the type of combinatorial descriptions that can be considered, in the form of a system of recursive equations relating combinatorial families of objects. The general form of such a system is 
\[\begin{cases}
\cal{Y}_1 =\cal{H}_1(\cal Z,\cal{Y}_1,\cal{Y}_2,\cdots,\cal{Y}_m), \\
\cal{Y}_2 =\cal{H}_2(\cal Z,\cal{Y}_1,\cal{Y}_2,\cdots,\cal{Y}_m), \\      
\qquad\qquad\vdots                    \\
\cal{Y}_m =\cal{H}_m(\cal Z, \cal{Y}_1,\cal{Y}_2,\cdots,\cal{Y}_m),
\end{cases}\]
that we denote concisely by $\bc Y=\bc H(\cal Z,\bc Y)$. 
Precise statements on what the $\cal H_i$ can be are given in \cref{sec:well-founded}.
For instance, permutations could be described by $\cal P=\Set(\Cyc(\mathcal Z))$ corresponding to their cycle decomposition, binary trees by the recursive equation $\cal B=\cal Z + \cal B\times \cal B$ and general (Cayley) trees by $\cal C=\cal Z\times\Set(\cal C))$. In order to automate this approach, a very first step is to determine when such a system is meaningful, excluding systems like $\cal T=\Set(\cal T)$ that would not have a useful combinatorial meaning. We define a system to be well founded when the iteration~$\bc Y^{[n+1]}=\bc H(\cal Z,\bc Y^{[n]})$ started with $\bc Y^{[0]}=\bs 0$ converges (in a combinatorial sense, see \cref{def:wf} and the following paragraphs). A simple algorithm (Algorithm \textsf{\ref{algo:wellfoundedandleadingterm}}) tests this condition. This algorithm is more efficient than \anonymousvariant{our}{the} earlier one in~\cite{PivoteauSalvySoria2012}; it can be seen as an evolution of Zimmermann's algorithm~\cite[Algorithm B, p.~33]{Zimmermann1991}, where the computation of valuations and the detection of cycles in the graph of dependence of the system are combined in a single iteration.  

This combinatorial setting is presented in \cref{sec:well-founded}, using a small amount of vocabulary from species theory that makes it easier to discuss the well-founded character of combinatorial systems. The combinatorial objects that can be defined with the constructions considered by Flajolet and Sedgewick are called constructible species, to distinguish them from the general case and the corresponding systems are called constructive. 
The algorithms presented in this article are generally stated for constructible species and their defining systems, while the mathematical results that prove their correctness are often given for arbitrary systems.

\paragraph{Generating Function Systems}
The next step of the approach is to introduce generating functions. 
In this article, we focus on the case of labeled structures and exponential generating functions. This means that the object of interest is the power series
\[\bs Y(z)=\sum_{n\ge0}\bs Y_{\!n}\frac{z^n}{n!},\]
where $\bs Y_{\!n}$ is a vector whose $k$th coordinate is the number of objects of type~$\mathcal Y_k$ with $n$ distinct atoms $\mathcal Z$ labeled from~1 to~$n$. 
For constructible species, the translation from $\bc Y=\bc H(\mathcal Z,\bc Y)$ to a system $\bs Y=\bs H(z,\bs Y)$  of generating function equations is straightforward thanks to an explicit dictionary going from combinatorial constructions to elementary functions~\cite[Chap.~II]{FlajoletSedgewick2009}, see \cref{tab:sum_esp_sg}.  

One can also define well-founded systems directly at the level of generating function equations.
This allows us to compare our definition with those studied \anonymousvariant{ by
Banderier and Drmota~\cite{BanderierDrmota2015, Drmota1997,Drmota2009} and Bell, Burris and Yeats~\cite{BellBurrisYeats2010,BellBurrisYeats2011}}{in \cite{BanderierDrmota2015, Drmota1997,Drmota2009,BellBurrisYeats2010,BellBurrisYeats2011,PivoteauSalvySoria2012}}, that have extra constraints meaningful in the contexts of their analyses. This is presented \cref{subsec:WFsys}, which can be used as a starting point by readers who want to focus on the analytic aspects.

\paragraph{Radius of Convergence}
Next comes the analytic part of analytic combinatorics. The first step is to determine the radius of convergence of the generating functions, viewed as analytic functions, as this governs the exponential growth of their coefficients. For a given~$a>0$ inside the domain of convergence of~$\bs Y$, the system~$\bs Y=\bs H(a,\bs Y)$ has $\bs Y(a)$ as a solution, i.e., the value of the generating function at~$a$. It may also have other solutions with nonnegative coordinates inside the domain of convergence of~$\bs H$. However, for~$a$ larger than the radius of convergence of~$\bs Y$, we prove that no such solution exists (\cref{thm:alerho}). This is a basis for the computability of the radius of convergence by dichotomy, provided one can detect the existence of nonnegative solutions to systems (\cref{th:radius-is-computable}). By the Tarski-Seidenberg theorem, this existence is decidable when the system~$\bs Y=\bs H(z,\bs Y)$ is a polynomial system. Combinatorially, this corresponds to a context-free grammar, or a combinatorial system that does not involve $\Cyc$ or $\Set$. In the more general situation where $\Cyc$ and $\Set$ are allowed, a result of Macintyre and Wilkie~\cite{MacintyreWilkie1996} shows that the decision can also be performed, and therefore the radius of convergence computed, provided a conjecture in transcendental number theory called Schanuel's conjecture holds. As this may well be stronger than needed, we specify more precisely an oracle whose decisions are sufficient for the computation of the radius of convergence, and to which all the decisions in our algorithms reduce.

\paragraph{Newton's Iteration}
Looking more closely at the structure of the systems of generating function equations~$\bs Y=\bs H(z,\bs Y)$ allows for more precise results, presented in \cref{sec:Newton}.
Notably, these systems exhibit strong positivity properties, very similar to those considered by Etessami, Stewart, Yannakakis~\cite{EtessamiYannakakis2009,EtessamiStewartYannakakis2013,StewartEtessamiYannakakis2015} for polynomial systems. This has an important consequence for Newton's iteration: for $a\in[0,\rho)$ where $\rho$ is the radius of convergence of~$\bs Y$, Newton's iteration started at~$\bs Y=\bs 0$ converges to~$\bs Y(a)$~\cite{PivoteauSalvySoria2012}. \Cref{thm:Newton} is a converse of this result: if Newton's iteration started at~$\bs Y=\bs 0$ converges to a point~$\bs B$ with nonnegative coordinates inside the domain of convergence of~$\bs H$, then $a\le \rho$  and $\bs B=\bs Y(a)$. In practice, one still needs to determine when to stop the iteration: it can well be that $a\le\rho$ and that the current iterate is close to~$\bs Y(a)$ but it could also happen that $a>\rho$ and the ultimate divergence of the iteration has not been detected yet. For this purpose, a Newton-Kantorovich-type theorem bounds the distance of the last iterate to the solution~(\cref{thm:Kantorovich}); positivity of the system allows this bound to be obtained from a single evaluation, instead of requiring an upper bound on a disk. \Cref{thm:speed} bounds the speed of convergence in the spirit of a similar result by Esparza, Kiefer, Luttenberger~\cite{EsparzaKieferLuttenberger2010}.

Newton's iteration can also be used to find the radius of convergence. However, the corresponding system does not have a strong positivity property anymore. Unconditional convergence starting from the origin does not hold in that case. Still,  \cref{thm:newton-radius,thm:newton-u-eq-1} show quadratic convergence in a neighborhood of the solution. A difficulty is that the characteristic system having the radius of convergence as a solution (\cref{def:characteristic-system}) may also have other solutions with positive coordinates. Bell, Burris and Yeats~\cite{BellBurrisYeats2010} have shown that the solution corresponding to the radius of convergence is the unique one where the largest eigenvalue of the Jacobian matrix is equal to~1. This result extends to our more general class of irreducible systems (\cref{thm:Lambda}).

\paragraph{Nonreal Dominant Singularities}
While the exponential growth of the coefficients of a generating function is governed by its radius of convergence~$\rho$, subexponential terms in their asymptotic behavior are governed by the local behavior of the generating function at all the singularities on the circle $|z|=\rho$, called dominant singularities. By a classical result of Pringsheim~\cite[Thm. IV.6]{FlajoletSedgewick2009}, as generating functions have nonnegative coefficients, the point $\rho$ itself is a singularity. There might be other ones on the circle~$|z|=\rho$. A simple example  (see \cref{ex:badsing}) is
\[f(z)=\frac{1-3z+6z^2}{(1-3z)(1-2z+9z^2)}=1+2z+z^2+2z^3+49z^4+\dotsb,\]
a power series with positive integers for coefficients, whose singularities lie at $e^{i\theta}/3$, with $\theta\in\{0,\pm \arccos(1/3)\}$. Generating functions coming from combinatorial systems behave more nicely than this example: the arguments of their dominant singularities are rational multiples of~$\pi$ (\cref{thm:BBY72i,thm:period-linear}) and moreover, these rational numbers are induced by periodicity properties that can be computed from the combinatorial equations. Results of that type are already known for $\mathbb N$-rational power series~\cite{Berstel1971}; for a class of nonrecursive combinatorial generating functions~\cite[Prop.~4.6]{FlajoletSalvyZimmermann1991}; for polynomial combinatorial systems~\cite[Lemma~26]{BellBurrisYeats2006}. This is the object of \cref{sec:dom-irred} and the resulting algorithm for constructible generating functions is given in \cref{sec:dominant}.
 
\paragraph{Singular Behavior} Once the dominant singularities have been located, the next step (\cref{sec:singular-behavior}) is to compute the expansions of the generating functions in the neighborhood of these singularities. In the case of positive systems of polynomial or entire functions, Banderier and Drmota have classified the possible singular behaviors at finite singularities~\cite{BanderierDrmota2015}. In particular, they proved that the exponents that can occur in the expansions are rational numbers whose denominator is a power of~2. 

For combinatorial systems that may involve~$\Cyc$ and~$\Set$, the possibilities are much more diverse. For instance, the generating function of the species $\Set(\mathcal Z\times\Seq(\mathcal Z))$ of fragmented permutations~\cite[Ex. VIII.7]{FlajoletSedgewick2009} has dominant singularity at~1 and behaves like $\exp(1/(1-z))$ as~$z\rightarrow1$. The species~$\Set(\Cyc(\mathcal Z)\times\Cyc(\mathcal Z))$ (that could be called `bicycle park') also has a generating function with dominant singularity at~1 and behaves like~$\exp(\ln^2(1-z))$ as~$z\rightarrow1$. \Cref{thm:alg-log-real} extends the classification of Banderier and Drmota to the general situation. There is a gap property  that was already known in the nonrecursive case~\cite[Prop.~4.8]{FlajoletSalvyZimmermann1991}: roughly speaking, as $z\rightarrow\rho$, either the generating function grows at least as fast as in this last example, or it has an algebraic-logarithmic behavior, i.e., either it has an infinite limit and is equivalent to
\begin{equation}\label{eq:alg-log-type}
C(1-z/\rho)^\alpha\ln^k\!\left(\frac1{1-z/\rho}\right),\qquad C>0,\alpha\le0,k\in\mathbb N,\alpha k\neq0
\end{equation}
or it behaves like
\[c_0+c_1(1-z/\rho)^{1/r}+c_2(1-z/\rho)^{2/r}+\dotsb\]
otherwise. In this last situation, as in the classification of Banderier and Drmota, the exponent~$1/r$ is a power of~2. In the first one, less regular exponents occur. \Cref{ex:transcendental-exponent} shows that 
the generating function of the species\footnote{We write $\cal A^k$ for the species $\underbrace{\cal A\times\dotsb\times \cal A}_{k \mbox{ \tiny times}}$.}
\[\mathcal Y=\Set\!\left(\mathcal Z\times\Cyc\!\left(\cal Z^2+\cal Z^2\right)\times\Set\!\left(\Seq(\cal Z^3)\right)\right)\]
has a dominant singularity at~$\rho=1/\sqrt2$, and as $z\rightarrow\rho-$, its behavior is
\[Y(z)\sim{2^\alpha}(1-z\sqrt2)^{-\alpha}\quad\text{with}\quad
\alpha=\frac{\sqrt2}2\exp\!\left(\frac{8+2\sqrt2}{7}\right).\]
Another difference with the case of systems of polynomial or entire functions is that the exponents at nonreal dominant singularities can be nonreal as well (\cref{thm:alg-log-dom}).

Despite these seemingly intricate exponents, the proofs of \cref{thm:alg-log-real,thm:alg-log-dom} give rise to algorithms that compute the singular expansions of constructible generating functions, provided they are of algebraic-logarithmic type. No new equality or inequality involving these constants need to be decided; all the delicate decisions that might possibly require an oracle have been made when computing the radius of convergence.

\paragraph{Asymptotic Expansion of the Coefficients}
Expansions of algebraic-logarithmic type~\eqref{eq:alg-log-type} are precisely those to which the transfer theorems of Flajolet and Odlyzko~\cite{FlajoletOdlyzko1990a} apply. The end result of this article (\cref{thm:everything}) is that asymptotic expansions at arbitrary precision can be computed for the coefficients of all constructible generating functions that have algebraic-logarithmic dominant singularities. The only hypothesis is that Schanuel's conjecture holds, or that the system contains neither~$\Set$ nor~$\Cyc$. This article is concluded by
\Cref{sec:conclusion}, where further effectivity questions in this area are discussed.

\section{Well-Founded Combinatorial Systems and their Characterizations}\label{sec:well-founded}
The symbolic method of Flajolet and Sedgewick~\cite{FlajoletSedgewick2009} relies on systems of combinatorial equations, called \emph{specifications}, between what they call \emph{classes of combinatorial structures}. 
Algorithms to determine that such a specification does define a class of structures are known~\cite{Zimmermann1991,PivoteauSalvySoria2012}. This section presents two  simpler characterizations (\cref{thm:gl-implicit-leading} with \cref{algo:iswf}, and \cref{thm:gl-implicit}). 
The language of the symbolic method is first recalled in \cref{sec:WF-combinatorial-systems,sec:constructive-systems-FS}.
It is sufficient to give a first characterization in~\cref{sec:characterization-leading}.
Then, the combinatorial derivative of species theory is recalled in \cref{sec:derivative} in order to present a second characterization in \cref{sec:characterization-nilpotent} with its proof that extends Joyal's implicit species theorem~\cite{Joyal81}. This second characterization is the one that translates more easily at the level of analytic properties; it provides us with a useful tool for many proofs. Readers familiar with analytic combinatorics but not with species theory should be able to read this section with an intuitive understanding of species. For more information, the book of Bergeron, Labelle and Leroux~\cite{BergeronLabelleLeroux1998} is our reference. A short introduction very close to our needs is given in \anonymousvariant{our previous work}{} \cite[\S1.1]{PivoteauSalvySoria2012}.

\subsection{Well-founded Combinatorial Systems}\label{sec:WF-combinatorial-systems}
The definition of specification used by Flajolet and Sedgewick~\cite[Def. I.7]{FlajoletSedgewick2009} uses the words `class', `term', `construction'. In the language of species theory~\cite{Joyal81,BergeronLabelleLeroux1998}, the situation becomes simpler as these words all correspond to (possibly multisort) species.
The combinatorial class described by the specification is then the species solution of a system of equations on species. 
\begin{definition}
A \emph{combinatorial system} for an $m$–tuple $\bc Y = (\cal{Y}_1, \cdots , \cal{Y}_m)$ of species is a
vector of $m$ equations,
\begin{equation}\label{eq:spec}
\begin{cases}
\cal{Y}_1 =\cal{H}_1(\cal Z,\cal{Y}_1,\cal{Y}_2,\cdots,\cal{Y}_m), \\
\cal{Y}_2 =\cal{H}_2(\cal Z,\cal{Y}_1,\cal{Y}_2,\cdots,\cal{Y}_m), \\      
\qquad\qquad\vdots                    \\
\cal{Y}_m =\cal{H}_m(\cal Z, \cal{Y}_1,\cal{Y}_2,\cdots,\cal{Y}_m),
\end{cases}
\end{equation}
where each $\cal H_i$ denotes a multisort species. (See \cref{tab:sum_esp_sg} for simple examples of possible $\cal H_i$.)
\end{definition}

\begin{notation} 
Boldfaced characters are used for vectors, matrices, or tuples of species; for example, a multisort species ${\mathcal H}({\mathcal Y}_1,{\mathcal Y}_2, \dots, {\mathcal Y}_k)$ is written ${\mathcal H}(\bc Y)$, where $\bc Y$ stands for the tuple $({\mathcal Y}_1,{\mathcal Y}_2, \dots, {\mathcal Y}_k)$; a vector of multisort species $({\mathcal H}_1(\bc Y), {\mathcal H}_2(\bc Y),\dots, {\mathcal H}_m(\bc Y))$  is consistently written $\bc H(\bc Y)$.

When several sorts of atoms are used, they will usually be denoted $(\cal Z_1,\cal Z_2,\dots)$. Thus, the most general form of combinatorial system with several sorts of atoms considerered here is written~$\bc Y=\bc H(\bc Z, \bc Y)$.

Gantmacher's notation $a_{1:k}$ is used to denote the $k$-tuple $(a_1,\dots, a_k)$. For instance, the species~$\bc H(\bc Y)$ can also be written~$\bc H_{1:m}(\bc Y_{1:k})$ if its dimensions need to be made explicit.
\end{notation}

The following general definition of well-founded combinatorial systems captures the families of systems considered in earlier literature. For a version of this definition directly at the level of generating functions, see \cref{subsec:WFsys}.

\begin{definition}\label{def:wf}
The combinatorial system~$\bc Y=\bc H(\bc Z,\bc Y)$ is \emph{well founded}  when the recurrence
\begin{equation}\label{eq:ite-def-wf}
                \bc Y^{[0]}=\bs{0}\quad\text{and}\quad\bc Y^{[n+1]}=\bc H(\bc Z,\bc Y^{[n]}),\quad n\ge0
\end{equation}
is well defined for all~$n$ and defines a sequence~$(\bc Y^{[n]})_{n\ge0}$ that is convergent. The limit~$\bc S$ of this sequence is called \emph{the combinatorial solution} of the system.
\end{definition}
In terms of combinatorial structures, `\emph{well defined}' means that the composition is possible, {i.e.}, for each~$n$ and each size~$k$, there are finitely many combinatorial structures of size~$k$ in $\bc H(\bc Z,\bc Y^{[n]})$. Convergence means that there is a limit class~$\bc S$ such that for all~$p\ge0$, $\bc S$ is isomorphic to $\bc Y^{[n]}$ up to the structures of size~$p$, for large enough~$n$. Alternatively, the corresponding sequence of generating functions converges as a sequence of power series for total degree.

Equivalently, in the context of species, `\emph{well defined}' means that for each coordinate $\mathcal Y_i^{[n]}$, $i=1,\dots,m$, either the species~$\bc H$ is polynomial in~$\mathcal Y_i^{[n]}$ or $\mathcal Y_i^{[n]}(0)=0$. The class $\bc S$ being isomorphic to $\bc Y^{[n]}$ up to the structures of size~$p$ translates into the species $\bc S$ having contact of order~$p$ with $\bc Y^{[n]}$, see~\cite{Joyal81,BergeronLabelleLeroux1998,PivoteauSalvySoria2012}.

\paragraph{Zero-free Systems and Comparison with~\cite{PivoteauSalvySoria2012}}
It is natural to focus on combinatorial systems such that the limit of the iteration in \cref{def:wf} has no zero coordinate.
\begin{definition}\label{def:zero-free}
A combinatorial system is called \emph{zero-free} if there exists $N>0$ such that for each coordinate $\mathcal Y_i$, $i=1,\dots,m$, $\mathcal Y_i^{[N]}\neq 0$, with $\bc Y^{[n]}$ defined by \cref{eq:ite-def-wf}.
\end{definition}

\anonymousvariant{Our}{The} definitions of well-founded systems in~\cite{PivoteauSalvySoria2012} included that condition. Here, we remove this constraint from our definitions, which leads to a slightly larger class of systems that satisfy the conditions of \cref{thm:gl-implicit-leading} and  \cref{thm:gl-implicit} (e.g., the equation $\mathcal Y=\mathcal Z\mathcal Y$ does). 

All the results of \anonymousvariant{our}{the} earlier work \anonymousvariant{}{\cite{PivoteauSalvySoria2012}} that are used in the present article, for this new definition, have been checked to hold with these extended conditions.
Moreover, these coordinates can be detected algorithmically (for instance, using Algorithm \textsf{\ref{algo:iswf}} below) and removed from the system. 
This will be used as hypothesis when needed or when it lightens the presentation, knowing that any well-founded system can be split into two sub-systems, one being zero-free and the other one defining only zero-coordinates (see \cref{coro:remove-zero}).

\subsection{Constructive Systems in the  Symbolic Method of Flajolet \& Sedgewick}\label{sec:constructive-systems-FS}

In order to distinguish the combinatorial systems used in the symbolic method from the general situation, we call \emph{specifications} those that are considered by Flajolet and Sedgewick~\cite[Def.~I.8]{FlajoletSedgewick2009}.
\begin{definition}\label{def:specification} 
A \emph{specification} is a combinatorial system of the form $\bc Y=\bc H(\cal Z,\bc Y)$ where the multisort species~$\mathcal H_i$ are obtained by combinations of~$1,\cal Z, +,\times,\Seq,\Cyc,\Set$ (see \cref{tab:sum_esp_sg}). \end{definition}

In species language, $\cal Z$ is the species of singletons; $1$ is the characteristic species of the empty set; $+$ denotes the sum (or disjoint union); the product $\times$ correspond to the cartesian product; sets, sequences and  cycles are easily defined~\cite{BergeronLabelleLeroux1998,PivoteauSalvySoria2012}.
This is summarized in \cref{tab:sum_esp_sg}, that also recalls how to compose the `construction' species with the species~0 (that corresponds to the empty class) and the species~1 (that corresponds to the class containing only~$\cal E$) and how to obtain the derivatives discussed in the following section. For simplicity, we omit here a discussion on cardinality constraints that can be added to the sequences, sets and cycles. They can also be handled, see the discussion in the conclusion.

\begin{table}[ht]
\begin{small} 
\begin{center}
\renewcommand{\arraystretch}{1.3}
$ \begin{array}{|l|l|c|c|l|l|} 
\hline
\mbox{species} & \mbox{construction} & ~\cG =0~ &  ~\cG = 1~ &  \mbox{derivative at } \cal Z&  \mbox{exponential g.s.}  \\[2pt] 
\hline  
\hline
\mbox{Disjoint union} & \cal F+ \cal G & \cal F & \cal F+ 1 & \cal F'+ \cal G' &F(z)+G(z)\\[2pt]
\mbox{Cartesian product} & \cal F\times \cal G & 0 & \cal F & \cal F'\cal G+ \cal F\cal G'&F(z)\times G(z)\\[2pt]   
\hline  
&&&& \mbox{derivative at }\cal G = \cal Z &   \\
\hline  
\mbox{Sequence (linear order)} &  \Seq(\cal G) & 1 & - &\Seq(\cal Z)\times\Seq(\cal Z) &(1-G(z))^{-1} \\[2pt]      
\mbox{Cycle} & \Cyc(\cG) & 0 & - &\Seq(\cal Z)&\ds\ln\frac1{1-G(z)}\\[2pt]                     
\mbox{Set}  & \Set(\cG) & 1 & -  & \Set(\cal Z) &\exp(G(z))\\[2pt]
\hline
\end{array}
$ 
\end{center}
\end{small}
\caption{Constructions involved in combinatorial systems in normal form, with special cases for composition, derivatives and associated generating functions {\footnotesize (``$-$'' stands for ``undefined'')}.
}
\label{tab:sum_esp_sg}
\end{table}

\begin{definition}\label{def:constructive} 
A well-founded zero-free specification is called a \emph{constructive} system. The species solution are called \emph{constructible}.
\end{definition}
The name of these systems comes from the fact that they are exactly the systems that define constructible combinatorial classes in the work of Flajolet and Sedgewick \cite[Def. I.7]{FlajoletSedgewick2009}.

In our algorithms, it is convenient to make use of systems of a special form similar to the Chomsky normal form for context-free grammars, see also~\cite[p.~21]{Zimmermann1991}, \cite{FlajoletZimmermanVanCutsem1994}. 
\begin{definition}\label{def:normal_form} 
    A specification~$\bc Y=\bc H_{1:m}(\cZ,\bc Y)$ is \emph{in normal form} when its equations are of one of the forms:
    \[
    \cY_i = 1; \quad \cY_i = \cZ; \quad 
    \cY_i = \cal G_1 + \,\cal G_2 + \dotsb; \quad 
    \cY_i = \cal G_1 \times \,\cal G_2 \times \dots; \quad 
    \cY_i = \cal F(\cal G_i); 
    \]   
    where each $\cal G_j \in \{ 1, \cZ, \cY_1,\dots,\cY_m \}$ and $\cal F$ is one of the constructions in $ \{ \Seq, \Cyc, \Set \}$ 
    (see \cref{tab:sum_esp_sg}).
\end{definition}

Any specification can easily be transformed into a system in normal form. Thus all the algorithms presented in this article can take as input a specification in normal form. Most of them  also require the zero-free constraint that is equally easy to satisfy.

\begin{example}\label{ex:colored_forest} The following specification will be a running example.\footnote{A Maple session dedicated to this example \anonymousvariant{can be found in the arXiv along with this article.}{ exists, but we have not found how to transmit it anonymously together with this submission.}} It defines unordered forests of trees with red, green and blue nodes and various other constraints. Its normal form is given on the right. For this example, it is possible to obtain closed forms for the generating functions. While these are not used by our algorithms, they allow to check the results easily.

\noindent\begin{minipage}{.3\textwidth}
\begin{small}
\begin{equation}\label{eq:colored_forest}
\begin{cases}
\cforest = \Set(\redroot),\\
\redroot = \cZ + \cred \times \Seq(\blueroot),\\
\blueroot = \cblue \times \Seq(\greenroot),\\
\greenroot = \cgreen \times \Seq(\redroot),\\
\cred = \cZ^3 + \cZ \times\Seq(\cblue),\\
\cblue = \cZ \times \Seq(\cred),\\
\cgreen = \cZ + \cgreen \times \cgreen.
\end{cases}
\end{equation}
\end{small}
\end{minipage}\hfill
\begin{minipage}{.65\textwidth}
\begin{small}
\begin{equation}\label{eq:normal_colored_forest}
\begin{cases}
\cforest = \Set(\redroot),\\
\redroot = \cZ + \cal{N}_1,\quad \cal{N}_1 = \cred \times \cal{N}_2,\quad \cal{N}_2 = \Seq(\blueroot),\\
\blueroot = \cblue \times \cal{N}_3,\quad \cal{N}_3 = \Seq(\greenroot),\\
\greenroot = \cgreen \times \cal{N}_4,\quad \cal{N}_4 = \Seq(\redroot),\\
\cred = \cal{N}_5 + \cal{N}_6, \!\!\!\quad \cal{N}_5 = \cZ^3,\quad \cal{N}_6 = \cZ \times \cal{N}_7, \quad \cal{N}_7 = \Seq(\cblue),\\
\cblue = \cZ \times \cal{N}_8,\quad \cal{N}_8 = \Seq(\cred),\\
\cgreen = \cZ + \cal{N}_9,\quad \cal{N}_9 = \cgreen \times \cgreen.
\end{cases}
\end{equation}  
\end{small}
\end{minipage}  
\end{example}

\subsection{Characterization of Well-founded Systems by Leading Terms}\label{sec:characterization-leading}

A first characterization of the systems~$\bc{Y}=\bc{H}(\cal{Z},\bc{Y})$ that are well founded can be given in terms of leading terms of the iterates.  In \cref{sec:characterization-nilpotent}, a second characterization, based on the Jacobian matrix of~$\bc{H}$, is given. It will be used in the proofs of~\cref{sec:generating-functions}.

\begin{definition}\label{def:leading}
Let $\cal F$ be a species. The valuation $\val({\cal F})$ of $\cal F$ is the size of smallest $\cal F$-structure (the valuation of the species~0 is infinite). The \emph{leading term} of $\cal F$ is the species $\cal F$ restricted to $\cal F$-structures of size $\operatorname{val}({\cal F})$ (and the leading term of the species~0 is itself).  The \emph{leading terms} of a vector of species $\bc F$ is the vector of the leading terms of $\bc F$. It is denoted $\LT(\bc F)$.
\end{definition}
\SetKwFor{RepTimes}{repeat}{times}{end}

\SetNlSkip{1em}
\IncMargin{1em}
\begin{algorithm}
\SetAlgoRefName{WellFoundedAndLeadingTerm}
\SetKwFunction{LeadingTerm}{LeadingTerm}
\caption{Leading terms of the generating functions\label{algo:wellfoundedandleadingterm}\label{algo:iswf}}
\Input{$\bc{Y}=\bc{H}(\cal{Z},\bc{Y})$ a specification in normal form with $m$ equations}
\Output{Vector of pairs (count, valuation) corresponding to the leading terms of $\bc{Y}$; FAIL if the system is not well founded\vspace{10pt}}
    $w_{1:m}:=((0,\infty),\dots,(0,\infty))$
    
    \RepTimes{m+1}{
        $v:=w$
        
        \For{$i=1$ to $m$}{
            \lIf{$\mathcal{H}_i=1$}{$w_i:=(1,0)$}
            \lIf{$\mathcal{H}_i=\mathcal Z$}{$w_i:=(1,1)$}
            \lIf{$\mathcal{H}_i=\Set(\mathcal Y_j)$ or $\mathcal{H}_i=\Seq(\mathcal Y_j)$}{$w_i:=(1,0)$ if $v_{j,2}\neq0$, otherwise FAIL}\label{line:fail_set}
            \lIf{$\mathcal{H}_i=\Cyc(\mathcal Y_j)$}{$w_i:=v_j$ if $v_{j,2}\neq 0$, otherwise FAIL}
            \lIf{$\mathcal H_i=\mathcal Y_{j_1}\times\cdots\times\mathcal Y_{j_k}$}{
                $w_i:=(v_{j_1,1}\dotsm v_{j_k,1},v_{j_1,2}+\dots+v_{j_k,2})$}
            \If{$\mathcal H_i=\mathcal Y_{j_1}+\dots + \mathcal Y_{j_k}$}{
                $w_{i,2}:= \operatorname{min}(v_{j_1,2},\dots,v_{j_k,2})$\\
                $\displaystyle w_{i,1}:= \sum_{\ell\in[1,k],v_{j_\ell,2} = w_{i,2}} v_{j_\ell,1} + \dots + v_{j_k,1}$
            }
        }
         \lIf{$v=w$}{\Return{$v$}}
    }
    \lIf{$v\neq w$}{FAIL}
    \Return{$v$}
\end{algorithm}
\begin{restatable}
[Characterization by Leading Terms]{theorem}{glimplicitleading}\label{thm:gl-implicit-leading}
Let $\bc H_{1:m}(\bc Z,\bc Y)$ be a vector of species such that 
the iteration \eqref{eq:ite-def-wf} is well defined for $n=0,\dots,m$. Then the system $\bc Y=\bc H(\bc Z,\bc Y)$ is well founded if and only if 
the leading terms of $\bc Y^{[m]}$ and $\bc Y^{[m+1]}$ are equal. In this situation, the combinatorial solution~$\bc S$ of the system satisfies $\LT(\bc S)=\LT(\bc Y^{[m]})$.
\end{restatable}
\noindent The proof is given in \cref{appendix:prooftheorem1.7}.

\medskip

Thus in at most $m+1$ iterations, one can detect whether a system is well founded or not.
This characterization at the level of species yields Algorithm \textsf{\ref{algo:iswf}} for specifications (where there is only one sort of atom~$\cal Z$). This algorithm is close in spirit with early work on these questions~\cite{Zimmermann1991,FlajoletSalvyZimmermann1991} for zero-free specifications and results in a much simpler algorithm than \anonymousvariant{our previous}{the}  one\anonymousvariant{~}{in~}\cite{PivoteauSalvySoria2012}. The algorithm given by Zimmermann~\cite{Zimmermann1991} (Algorithm B, p.~33) and summarized in~\cite{FlajoletSalvyZimmermann1991} first computes the valuations (this corresponds to computing the second coordinate in Algorithm \textsf{\ref{algo:iswf}}). Next, a cycle detection is performed on the dependency graph of the system (for a restricted notion of dependency). As explained in the next section, this corresponds to checking the nilpotence of the Jacobian matrix of~$\bc{H}$ at different iterates, which will appear in our second characterization.
Our new characterization by leading terms allows us to bound the number of iterations in the computation of valuations and results in a more straightforward algorithm. Also, it does not require removing the zero coordinates beforehand and can be used to detect them (\cref{coro:remove-zero} below).

\begin{theorem}\label{thm:algo-leading} Given as input a specification in normal form, Algorithm \textsf{\ref{algo:iswf}} returns FAIL if and only if the system is not well founded. Otherwise, it returns a vector of pairs that, for each coordinate, indicate the minimal size of the structures defined by that coordinate (its valuation) and the number of structures of that size.
\end{theorem}
\begin{proof}
At any step $k\in\{1,\dots,m+1\}$, either Algorithm \textsf{\ref{algo:iswf}} computes the vector of pairs (count, valuation) corresponding to the leading terms of $\bc{Y}^{[k]}$ or the system is not well founded. This is proved by induction using \cref{lemma:leading-leading}: if $\bc{Y}^{[k+1]}$ is well defined, $\LT(\bc{Y}^{[k+1]})=\LT(\bc{H}(\cal Z,\bc{Y}^{[k]}))=\LT(\bc{H}(\cal Z,\LT(\bc{Y}^{[k]})))$; the rest follows from basic properties of the constructions of \cref{def:specification}. It is straightforward for most of the constructions:
\[
\LT(1)=1,\quad
\LT(\cal Z)=\cal Z,\quad
\LT(\Set(\cal Y))=\LT(\Seq(\cal Y))=1,\quad
\LT(\Cyc(\cal Y))=\LT(\cal Y),
\]
an slightly more involved for products:
\[
\LT(\cal Y_1\times\dots\times\cal Y_n)=\sum_{\gamma_1\in\LT(\cal Y_1),\dots,\gamma_n\in\LT(\cal Y_n)} (\gamma_1,\dots,\gamma_n),
\]
and sums:
\[
\LT(\cal Y_1+\dots+\cal Y_n)=\sum_{\gamma\in\LT(\cal Y_1) + \dots + \LT(\cal Y_n),\;|\gamma|=v_+} \gamma
\]
where $v_+:=\val(\cal Y_1+\dots+\cal Y_n)=\min\left(\val(\cal Y_1),\dots,
\val(\cal Y_n)\right)$.
If the algorithm does not fail, a vector corresponding to the leading terms of $\LT(\bc{Y}^{[k]}) = \LT(\bc{Y}^{[k+1]}) = \dots=\LT(\bc{Y}^{[m]})$, for $k\leq m$, is returned. By \cref{thm:gl-implicit-leading}, these are the leading terms of the species solution of the system.
\end{proof}

\begin{corollary}\label{coro:remove-zero}
Given a well-founded specification $\bc Y=\bc H(\mathcal Z,\bc Y)$, Algorithm \textsf{\ref{algo:iswf}} detects its zero-coordinates as those for which it returns $(0,\infty)$. If $k$ is such that $\bc Y_{k+1:m}$ are the zero-coordinates, then the system $\bc Y_{1:k}=\bc H_{1:k}(\mathcal Z,\bc Y_{1:k},\bs 0)$ defines the same solutions as the nonzero ones in the initial system.
\end{corollary}
\begin{proof}
The zero coordinates are those that have an infinite valuation.
\end{proof}

\begin{example}On the system~\eqref{eq:normal_colored_forest} of \cref{ex:colored_forest}, the first iterations of the main loop in Algorithm \textsf{\ref{algo:iswf}} produce the following vectors, with each pair $(c,e):=$(coefficient, exponent) written as $c \cZ^{e}$:
\[\arraycolsep=1.4pt\def\arraystretch{1.2}
\begin{array}{@{[\,}*{16}{r}@{\,]}}
1,& \cZ,& 0,& 1,& 0,& 1,& 0,& 1,& 0,& \phantom{\,\,}\cZ^{3},& 0,& 1,& 0,& 1,& \cZ,& 0 \\
1,& \cZ,& 0,& 1,& 0,& 1,& \cZ,& 1,& \cZ^{3},& \cZ^{3},& \cZ,& 1,& \cZ,& 1,& \cZ,& \cZ^{2}\\
1,& \cZ,& \cZ^{3},& 1,& \cZ,& 1,& \cZ,& 1,& \cZ,& \cZ^{3},& \cZ,& 1,& \cZ,& 1,& \cZ,& \cZ^{2}\\
1,& \cZ,& \cZ,& 1,& \cZ,& 1,& \cZ,& 1,& \cZ,& \cZ^{3},& \cZ,& 1,& \cZ,& 1,& \cZ,& \cZ^{2} \\
1,& 2\cZ,& \cZ,& 1,& \cZ,& 1,& \cZ,& 1,& \cZ,& \cZ^{3},& \cZ,& 1,& \cZ,& 1,& \cZ,& \cZ^{2} \\
1,& 2\cZ,& \cZ,& 1,& \cZ,& 1,& \cZ,& 1,& \cZ,& \cZ^{3},& \cZ,& 1,& \cZ,& 1,& \cZ,& \cZ^{2}
\end{array}
\]
The last two vectors are equal and therefore, the system of 16~equations is found to be well founded in 6~iterations.
\end{example}

\begin{example}If the last two equations of the system defined by \cref{eq:normal_colored_forest} are replaced with another species of trees, $\cgreen = 1 + \cal N_9$, $\cal N_9 = \cZ \times \cgreen \times \cgreen$, the system is not well founded anymore. Indeed, the second iteration produces
\[\arraycolsep=1.4pt\def\arraystretch{1.2}
\begin{array}{@{[\,}*{16}{r}@{\,]}}
1,&{\cZ},&0,&1,&0,&1,&{\color{blue}1},&1,&{\cZ}^{3},&{\cZ}^{3},&\cZ,&1,&\cZ,&1,&1,&\cZ
\end{array}
\]
which means that $\LT(\cal T_g) = 1$ (in blue), which will break the condition in line~\ref{line:fail_set} for equation $\cal{N}_3 = \Seq(\greenroot)$ at the next iteration.
However, if this equation is also rewritten as $\cal{N}_3 = 1 + \cal M$, $\cal M = \greenroot\times\cal N_3$, the condition in the loop now holds. But in this case, the last two vectors after $m+1=18$ iterations are not equal:
\[\arraycolsep=1.4pt\def\arraystretch{1.2}
\begin{array}{@{[\,}*{18}{r}@{\,]}}
&1,&2\cZ,&\cZ,&1,&8\cZ,&{\color{blue}8},&8,&1,&1,&\cZ,&{\cZ}^{3},&\cZ,&1,& \cZ, & 1 ,&1,& \cZ\\
&1,&2\cZ,&\cZ,&1,&8\cZ,&{\color{blue}9},&8,&1,&1,&\cZ,&{\cZ}^{3},& \cZ,&1,&\cZ, & 1,&1,& \cZ
\end{array}
\]
showing that the new system is not well founded.
\end{example}
\begin{remark}
    The number of iterations in Algorithm \textsf{\ref{algo:iswf}} is bounded by the number of equations in the input system. Using a system in normal form allows for a simple presentation of the algorithm but it can increase significantly the number of iterations, in comparison with the original system (as between the ones of \cref{eq:colored_forest} and \cref{eq:normal_colored_forest}). It is always possible to work with a system that is not in normal form (as in~\cite{Zimmermann1991}), by making explicit how to compute the leading terms of any $\cal H_i$ defined as in \cref{def:specification}.
\end{remark}

\subsection{Combinatorial Derivative}\label{sec:derivative}

The derivative of species theory is crucial in this work. In terms of classes of structures, 
the derivative of a combinatorial class $\cal A$ (with respect to $\cal Z$) is another combinatorial class containing all the possible structures obtained by erasing a $\cal Z$ in a structure of $\cal A$ and replacing it by a distinguished $\cal E$. 
Then, the product of a derivative with a class $\cal B$ can be interpreted as replacing a distinguished $\cal E$ in the derivative by a structure in $\cal B$. While Flajolet and Sedgewick do not define a derivative, their closely related pointing construction can be stated in terms of derivatives: $\Theta \cal A(\cal Z) = A'(\cal Z)\times \cal Z$.

Partial derivatives are obtained when considering constructions or terms with several sorts of atoms. For example, erasing an atom of the sort $\cal{Y}_1$ in all possible ways in the class $\cal{H}_i(\cal Z,\cal{Y}_1,\dots,\cal{Y}_m)$ taken from \cref{eq:spec} produces the class $\partial\cal H_i/\partial\cal Y_1$. Structures of $\partial\cal H_i/\partial\cal Y_1\times\partial\cal H_1/\partial\cal Y_j$ are structures of $\cal{H}_i(\cal Z,\cal{Y}_1,\dots,\cal{Y}_m)$ where a $\cal{Y}_1$ has been replaced by a structure of the derivative of $\cal H_1$ with respect to the sort $\cal Y_j$.

\begin{example}Let $\cal H = \Seq(\cal Y_1+\cal Y_2)$. The class $\partial\cal H/\partial\cal Y_1\times \cal E$ can be interpreted 
as the set of sequences of atoms of sort $\cal Y_1$ or $\cal Y_2$ where one of them as been erased and its position remembered. 
Similarly, the class $\partial\cal H/\partial\cal Y_1\times \cal Y_1$ is viewed as the set of sequences of $\cal Y_1$ or $\cal Y_2$ where one~$\cal Y_ 1$ is distinguished, and 
the class $\partial\cal H/\partial\cal Y_1\times \cal Y_3$ as the set of sequences of $\cal Y_1$ or $\cal Y_2$ and only one~$\cal Y_3$
\end{example}

Matrices and vectors of combinatorial classes are defined as usual 
and the classical notion of Jacobian matrix extends to combinatorial systems.
The product of a matrix by a matrix or a vector is also given by the usual rules (with disjoint unions and cartesian products as sums and products), which allows to define \emph{nilpotent} combinatorial matrices as usual.
This notion from species theory is the key element to decide whether a system is well founded.

\begin{example}The Jacobian matrix for the combinatorial system of \cref{eq:colored_forest} is
\[\resizebox{\textwidth}{!}{$
\begin{pmatrix} 
0\hspace*{4mm} & \Set(\redroot)& 0\hspace*{4mm} & 0 & 0 & 0 & 0 \\
0\hspace*{4mm} & 0\hspace*{4mm} & \hspace*{-4mm}\cred \times\Seq(\blueroot)^2 & 0 & \Seq(\blueroot) & 0 & 0 \\
0\hspace*{4mm} & 0\hspace*{4mm} & 0\hspace*{4mm} & \cblue\times\Seq(\greenroot)^2 & 0 & \Seq(\greenroot) & 0 \\
0\hspace*{4mm} & \hspace*{-4mm}\cgreen\times\Seq(\redroot)^2 & 0\hspace*{4mm} & 0 & 0 & 0 & \Seq(\redroot) \\
0\hspace*{4mm} & 0\hspace*{4mm} & 0\hspace*{4mm} & 0 & 0 & \cZ\times\Seq(\cblue)^2 & 0 \\
0\hspace*{4mm} & 0\hspace*{4mm} & 0\hspace*{4mm} & 0 & \cZ\times\Seq(\cred)^2 & 0 & 0 \\
0\hspace*{4mm} & 0\hspace*{4mm} & 0\hspace*{4mm} & 0 & 0 & 0 & 2\cgreen 
\end{pmatrix}.$} 
\]  
\end{example}

\begin{definition}\label{def:dependency-graph}
The {\em dependency graph}~\cite[p. 33]{FlajoletSedgewick2009} of a combinatorial system $\bc Y=\bc H(\cZ,\bc Y)$ is the directed graph whose vertices are the $\cal Y_i$s and such that there is  an edge from $\cal Y_i$ to $\cal Y_j$ when $\cal Y_j$ 
appears explicitly on the right-hand side of the equation for  $\cal Y_i$. Equivalently, the dependency matrix of this graph is given by the location of the nonzero entries of the Jacobian matrix ${\bpartial \bc H}/{\bpartial\bc Y}(\cZ,\bc Y)$.
\end{definition}

\begin{example}\label{ex:dependency_graph}For the system of \cref{eq:colored_forest}, the dependency graph of $(\cforest,\redroot,\blueroot,\greenroot,\cred,\cblue,\cgreen)$ is
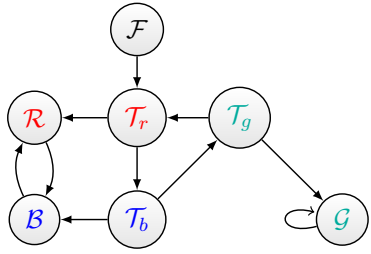
\begin{center}
\scalebox{0.9}{
\begin{tikzpicture}[,-latex ,auto ,node distance =1.5 cm and 1.5cm ,on grid , semithick ,
state/.style ={ circle ,top color =white , bottom color = gray!20 ,
draw , minimum width =0.7cm}]

\node[state] (F) {$\cforest$};
\node[state] (Tr) [text=red,  below =of F,yshift=2mm] {$\redroot$};
\node[state] (Tb) [text=blue, below =of Tr] {$\blueroot$};
\node[state] (Tg) [text=mygreen,right =of Tr]{$\greenroot$};

\node[state] (G)[text=mygreen,below right =of Tg] {$\cgreen$};
\node[state] (R)[text=red,left =of Tr] {$\cred$};
\node[state] (B)[text=blue,left =of Tb] {$\cblue$};

\path (F) edge  node[left] {} (Tr);

\path (Tr) edge  node[left] {} (Tb);
\path (Tb) edge  node[left] {} (Tg);
\path (Tg) edge  node[left] {} (Tr);

\path (Tr) edge  node[left] {} (R);
\path (Tb) edge  node[left] {} (B);
\path (Tg) edge  node[left] {} (G);

\path (R) edge  [bend left =25] node[left] {} (B);
\path (B) edge  [bend left =25] node[left] {} (R);
\path (G) edge  [loop left] node[left] {} (G);
\end{tikzpicture}}
\end{center}
\end{example}

\subsection{Characterization of Well-founded Systems by Jacobian Matrices}\label{sec:characterization-nilpotent}

The following characterization and its proof use the language of species theory, drawing in part 
from \anonymousvariant{our previous}{the} work\anonymousvariant{~}{ in~}\cite{PivoteauSalvySoria2012}. 
The starting point is the following theorem by Joyal~\cite{Joyal81}.

\begin{theorem}[Joyal's implicit species theorem]\label{th:IST}
If $\bc H(\bs{0},\bs{0})=\bs{0}$ and ${\bpartial \bc H}/{\bpartial\bc Y}(\bs{0},\bs{0})$ is nilpotent, then
the system $\bc Y=\bc H(\bc Z,\bc Y)$ is well founded and has a solution~$\bc S$ with $\bc S(\bs 0)=\bs 0$, that is unique up to isomorphism.
\end{theorem}
A partial converse of this result (for zero-free systems) holds, showing that its condition is essentially tight.
\begin{proposition}\cite[Thm.~3.6]{PivoteauSalvySoria2012}\label{prop:well-founded-at0-implies-nilp}
If $\bc Y=\bc H(\bc Z,\bc Y)$ is well founded with $\bc H(\bs 0,\bs 0)=\bs 0$ and the limit~$\bc S$ of Iteration~\eqref{eq:ite-def-wf} has no zero coordinate, then ${\bpartial \bc H}/{\bpartial\bc Y}(\bs{0},\bs{0})$ is nilpotent.
\end{proposition}

These results generalize to the case when $\bc H(\bs 0,\bs 0)$ is not necessarily~$\bs 0$. The analogue of Joyal's theorem is the following.
\begin{theorem}[Implicit Species Theorem]\label{thm:gl-implicit}
Let $\bc H_{1:m}(\bc Z,\bc Y)$ be a vector of species. If 
\begin{enumerate}
    \item letting $\bc U^{[0]}:=\bs 0$, for each $k=0,\dots,m$, $\bc U^{[k+1]}:=\bc{H}(\bs 0,\bc U^{[k]})$ is well defined;
    \item the matrix $\bpartial\bc H/\bpartial\bc Y(\bs 0,\bc U^{[m]})$ is nilpotent,
\end{enumerate}
then the system $\bc Y=\bc H(\bc Z,\bc Y)$ is well founded and has a solution $\bc S$ with $\bc S(\bs 0)=\bc U^{[m]}$, which is unique up to isomorphism.
\end{theorem}
\begin{proof}

The sequence $(\bc U^{[k]})_{k\geq0}$ 
satisfies $\bc U^{[k]}\subset\bc U^{[k+1]}$ by induction: this is true for~$k=0$ and for larger~$k$, the inclusion is preserved by composition with~$\bc H$.
Next, the identity
\[
\bc U^{[k+1]}-\bc U^{[k]}\subset{\bpartial\bc H}/{\bpartial\bc Y}(\bs 0,\bc U^{[k]}) \times(\bc U^{[k]}-\bc U^{[k-1]}),\qquad k>0
\] 
shows that any $\bc U^{[k+1]}$-structure that is not a $\bc U^{[k]}$-structure is an $\bc H$-assembly of~$\bc U^{[k]}$-structures, where at least one lies in the difference~$\bc U^{[k]}-\bc U^{[k-1]}$. Iterating and using the inclusion~$\bc U^{[k]}\subset\bc U^{[k+1]}$ for all~$k\ge0$ proves
\[
\bc U^{[m+1]}-\bc U^{[m]}\subset\left({\bpartial \bc H}/{\bpartial\bc Y}(\bs 0,\bc{U}^{[m]})\right)^m\times(\bc U^{[1]}-\bc U^{[0]}).
\]
Since ${\bpartial \bc H}/{\bpartial\bc Y}(\bs 0,\bc{U}^{[m]})$ is nilpotent of order~$p\leq m$, the right-hand side is~$\bs{0}$ and~$\bc U^{[m+1]}=\bc U^{[m]}$.

Let us now consider the shifted system $\bc Y=\bc K(\bc Z,\bc Y)$ with $\bc K(\bc Z,\bc Y)=\bc H(\bc Z,\bc U^{[m]}+\bc Y)-\bc U^{[m]}$; the subtraction is well defined since 
\[\bc U^{[m]}=\bc U^{[m+1]}= \bc H(\bs{0},\bc U^{[m]})\subset \bc H(\bc Z,\bc U^{[m]})\subset \bc H(\bc Z,\bc U^{[m]}+\bc Y).\]

This shifted system satisfies the conditions of Joyal's Implicit Species Theorem~\ref{th:IST}, that are $\bc K(\bs{0},\bs{0})=\bc H(\bs 0,\bc U^{[m]})-\bc{U}^{[m]}=\bs{0}$ and ${\bpartial \bc K}/{\bpartial\bc Y}(\bs{0},\bs{0})=\bpartial\bc H/\bpartial\bc Y(\bs 0,\bc U^{[m]})$ nilpotent. 
Therefore, it admits a solution~$\bc S$ such that $\bc S(\bs 0)=\bs 0$, unique up to  isomorphism. Now, $\bc T:=\bc U^{[m]}+\bc S$ is a solution of~$\bc Y=\bc H(\bc Z,\bc Y)$ such that~$\,\bc T(\bs{0})=\bc U^{[m]}$. Any other solution~$\bc T'$ is such that~$\bc T'-\bc U^{[m]}$ is isomorphic to~$\bc S$, which concludes the proof.
\end{proof}

\begin{example}The system $\cal Y_1=\cal Z+\cal Z\cal Y_1,~\cal Y_2=\cal Y_1 + \cal Y_2^2$ defines binary trees with leaves that are nonempty sequences. \Cref{thm:gl-implicit} shows that this is well founded since the first condition naturally holds and defines $\bc U^{[0]}=\bc U^{[1]}=\bc U^{[2]}=(0,0)$ and the Jacobian matrix
$\left(\begin{smallmatrix}
\cal Z & 0\\ 
1 & 2\cal Y_2
\end{smallmatrix}\right)$ is nilpotent at $\cal Z = 0$ and $\bc Y = (0,0)$.

On the contrary, the system $\cal Y_1=1+\cal Z\cal Y_1,~\cal Y_2=\cal Y_1 + \cal Y_2^2$ uses possibly empty sequences for the leaves, which can create infinitely many trees of size 0. The system is not well founded, since $\bc U^{[2]}$ is now $(1,1)$ and the Jacobian matrix is not nilpotent at $\cal Z = 0$ and $\bc Y = (1,1)$.
\end{example}

\begin{remark}
It is not sufficient that $\bc U^{[m]}=\bc U^{[m+1]}$ to ensure that the system is well founded. 
For instance, the system defined by $\bc H = (\cal Z + \cal Y_1\cal Y_2, \cal Y_3, 1)$ has $\bc U^{[1]}=(0,0,1)$, $\bc U^{[2]}=(0,1,1)=\bc U^{[3]}=\bc U^{[4]}$. However, this system produces infinitely many $\cal Y_1$-structures of size~$1$. The second condition of the theorem does not hold: the Jacobian matrix has a block decomposition with~1 in the top-left corner. (This is also detected by Algorithm \textsf{\ref{algo:iswf}} since $\LT(\bc Y^{[3]})=(2\mathcal Z,1,1)$ and $\LT(\bc Y^{[4]}=(3\mathcal Z,1,1)$).
This example shows also that in the second point of \cref{thm:gl-implicit}, one cannot replace $\bc U^{[m]}$ by $\bc H(\bs 0, \bs 0)$ in the Jacobian matrix.  
\end{remark}

\noindent Actually, the conditions of \cref{thm:gl-implicit} are not far from optimal, as shown by the following.
\begin{proposition}\label{prop:convergent-implies well-founded}
If the system $\bc Y=\bc H(\bc Z,\bc Y)$ is zero-free and well founded, then the conditions of Theorem~\ref{thm:gl-implicit} hold.
\end{proposition}
\begin{proof}By induction, for all $k=0,\dots,m+1$, $\bc U^{[k]}\subset\bc Y^{[k]}$. Part 1. of \cref{thm:gl-implicit} follows from $\bc H(\bs 0,\bc U^{[k]})\subset \bc H(\bc Z,\bc Y^{[k]})$, which is well defined.

Next, Theorem~5.5 of~\cite{PivoteauSalvySoria2012} shows that, under the assumptions of this proposition,  the companion system~$\bc Y=\bc K(\mathcal{Z}_1,\bc Z,\bc Y)=\bc H(\bc Z,\bc Y)-\bc H(\bs 0,\bs0)+\mathcal{Z}_1\bc H(\bs 0,\bs 0)$ is well founded; it has a solution $\bc S(\mathcal Z_1,\bc Z)$ such that the following properties hold: $\bc S_1(0,\bs 0)=\bs 0$; the limit~$\bc S(\bc Z)$ is $\bc S(1,\bc Z)$; $\bc S(\mathcal Z_1,\bc Z)$ is polynomial in~$\mathcal Z_1$. 
By~\cite[Prop. 4.5]{PivoteauSalvySoria2012}, this implies that $\bpartial\bc K/\bpartial\bc Y(\mathcal{Z}_1,\bs 0,\bc S(\mathcal Z_1,\bs 0))=\bpartial\bc H/\bpartial\bc Y(\bs 0,\bc S(\mathcal Z_1,\bs 0))$ is nilpotent. This implies that $\bpartial\bc H/\bpartial\bc Y(\bs 0,\bc S(1,\bs 0))=\bpartial\bc H/\bpartial\bc Y(\bs 0,\bc S(\bs 0))$ is nilpotent too. Finally, by~\cite[Thm.~5.7]{PivoteauSalvySoria2012}, $\bc S(\bs 0)=\bc H^m(\bs 0,\bs 0)=\bc U^{[m]}$.
\end{proof}

\noindent\begin{minipage}{.59\textwidth}
\begin{example}Consider again the combinatorial system of \cref{eq:colored_forest}. $\bc U^{[7]}$ is the vector whose coordinates indexed by~$(\cal F,\cal T_r,\cal T_b,\cal T_g,\cal R,\cal B,\cal G)$ are $(1, 0, 0, 0, 0, 0,  0)$. Here is the Jacobian matrix $\bpartial\bc H/\bpartial\bc Y(\bs 0,\bc U^{[7]})$. It is nilpotent, showing that the system is well founded.
\end{example}
\end{minipage}\hfill
\begin{minipage}{.4\textwidth}
\begin{scriptsize}
    
\[
\begin{pmatrix} 
0 & 1 & 0 & 0 & 0 & 0 & 0 \\
0 & 0 & 0 & 0 & 1 & 0 & 0 \\
0 & 0 & 0 & 0 & 0 & 1 & 0 \\
0 & 0 & 0 & 0 & 0 & 0 & 1 \\
0 & 0 & 0 & 0 & 0 & 0 & 0 \\
0 & 0 & 0 & 0 & 0 & 0 & 0 \\
0 & 0 & 0 & 0 & 0 & 0 & 0
\end{pmatrix} 
\]  
\end{scriptsize}
\end{minipage}

\begin{example}If the last two equations of the system of \cref{eq:normal_colored_forest} are replaced with another species of trees, $\cgreen = 1 + \cal N_9$, $\cal N_9 = \cZ \times \cgreen \times \cgreen$, the system is no longer well founded.  Indeed,  $\bc U^{[3]}$ breaks the first condition of \cref{th:IST} (in the equation $\cal N_3 =  \Seq(\greenroot)$).
However, if the equation $\cal{N}_3 = \Seq(\greenroot)$ is also rewritten as $\cal{N}_3 = 1 + \cal M$, $\cal M = \greenroot\times\cal N_3$, the condition in the loop now holds. But in this case, the resulting Jacobian matrix is not nilpotent.
\end{example}

\begin{example}As regards 0-coordinates, it is straightforward to check that the conditions of \cref{thm:gl-implicit} hold for the system $\cal Y = \cal Z \cal Y$ which is well founded and defines 0. Still, the condition `zero-free' in \cref{prop:convergent-implies well-founded} is necessary: the system $\cal Y = \cal Y + \cal Z\cal Y$ is also well founded as the sequence defined by \cref{eq:ite-def-wf} converges to $(0)$, even though the second condition of \cref{thm:gl-implicit} does not hold. 
\end{example} 

\subsection{Decomposition into Strongly Connected Components}\label{sec:decomp}

\noindent Some characteristics of a well-founded system  can be found on its dependency graph:
\begin{definition}\label{def:combi_recursive}
A well-founded system $\bc Y=\bc H(\cZ,\bc Y)$ is \emph{recursive} when its dependency graph has a cycle, or equivalently when the matrix $\bpartial\bc H/\bpartial\bc Y$ is not nilpotent.
\end{definition}
\begin{definition}\label{def:combi_irreducible}
A combinatorial system $\bc Y=\bc H(\cZ,\bc Y)$ is \emph{irreducible} when it is well founded, zero-free, recursive and its dependency graph is strongly connected (every vertex can be reached from any of the other ones).
\end{definition}

Note that these definitions would also apply to combinatorial systems that are not well founded but  this situation is not considered in the rest of this article.

\begin{example}\label{ex:scc}The combinatorial system of \cref{eq:colored_forest} is recursive : the corresponding Jacobian matrix is not nilpotent and, indeed, the corresponding dependency graph has many cycles.
It is not irreducible since four strongly connected components can be identified : $(\redroot,\blueroot,\greenroot)$, $(\cred,\cblue)$ and $(\cgreen)$ are irreducible while $(\cforest)$ is nonrecursive.\end{example}

The analysis of well-founded systems is made simpler by focusing on \emph{irreducible} systems first, and reducing the general case to the irreducible and the reducible cases separately, through the classical decomposition of a directed graph into a {directed \emph{acyclic} graph} (DAG) of \emph{strongly connected components}. 
Since the components form a DAG, without loss of generality, the species can be relabeled in such a way that the $k$~successive parts compose a block-triangular system of the form
\begin{equation}\label{eq:triangular}
\bc{Y}_{1:i_1}=\bc{H}_{1:i_1}(\cal{Z},\bc{Y}_{1:i_1}),
\bc{Y}_{i_1+1:i_2}=\bc{H}_{i_1+1:i_2}(\cal{Z},\bc{Y}_{1:i_2}),\dots,
\bc{Y}_{i_{k-1}+1:m}=\bc{H}_{i_{k-1}+1:m}(\cal{Z},\bc{Y}_{1:m}).\quad
\end{equation}
These blocks will be analyzed from the first one to the last, so that any species is analyzed before it occurs in the definitions of other species.
\begin{definition}\label{def:irreducible-component-com} By analogy with graphs, the blocks in the system \eqref{eq:triangular} are called \emph{components}. The ones consisting of a single equation with $\partial \mathcal H_i/\partial \mathcal Y_i=0$ are the \emph{nonrecursive} components. 
The other ones are called the \emph{irreducible} components of the system. 
Components of constructive systems are called \emph{constructive components}.
{The definition of {\em normal form} also extends to any component 
$\bc{Y}_{i+1:j}=\bc{H}_{i+1:j}(\cal{Z},\bc{Y}_{1:i},\bc{Y}_{i+1:j})$
by including 
$\bc{Y}_{1:i}$ in the set from where $\cal G_j$ can be taken in \cref{def:normal_form}.}
\end{definition}

\begin{example}The dependency graph of the system given in~\cref{ex:dependency_graph} decomposes into four strongly connected components identified in the left of the following picture.
The corresponding DAG is given on the right. 
Both the subsystems $(\cgreen)$ and $(\cred, \cblue)$ are irreducible (and thus, irreducible components).
Note that the system $\{
\redroot =\cZ + \cred \times \Seq(\blueroot),
\blueroot = \cblue \times \Seq(\greenroot),
\greenroot = \cgreen \times \Seq(\redroot)
\}$ is an irre\-ducible component defining $(\redroot, \blueroot, \greenroot)$ while the subsystem  $(\redroot, \blueroot, \greenroot, \cgreen,\cred, \cblue)$ is not irreducible. The system~$\{\cforest = \Set(\redroot)\}$ defining $\cforest$ is a strongly connected component in the graph, but being nonrecursive, it does not form an irreducible component.

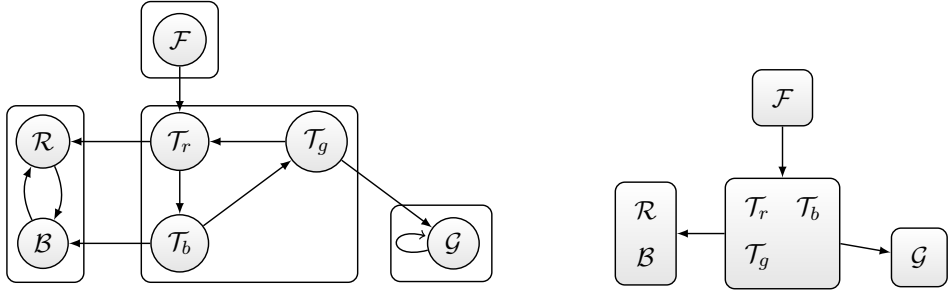
\scalebox{0.9}{
\begin{minipage}{.62\textwidth}
\qquad\quad\begin{tikzpicture}[,-latex ,auto ,node distance =1.5 cm and 2cm ,on grid , semithick ,
state/.style ={ circle ,top color =white , bottom color = gray!20 ,
draw , minimum width =0.7cm}]

\node[state] (F) {$\cforest$};
\node[state] (Tr) [below =of F] {$\redroot$};
\node[state] (Tb) [below =of Tr] {$\blueroot$};
\node[state] (Tg) [right =of Tr]{$\greenroot$};

\node[state] (G)[below right =of Tg] {$\cgreen$};
\node[state] (R)[left =of Tr] {$\cred$};
\node[state] (B)[left =of Tb] {$\cblue$};

\path (F) edge  node[left] {} (Tr);

\path (Tr) edge  node[left] {} (Tb);
\path (Tb) edge  node[left] {} (Tg);
\path (Tg) edge  node[left] {} (Tr);

\path (Tr) edge  node[left] {} (R);
\path (Tb) edge  node[left] {} (B);
\path (Tg) edge  node[left] {} (G);

\path (R) edge  [bend left =25] node[left] {} (B);
\path (B) edge  [bend left =25] node[left] {} (R);
\path (G) edge  [loop left] node[left] {} (G);

\node[draw, rectangle, rounded corners, inner sep=16pt] at (F) {};
\node[draw, rectangle, rounded corners, inner sep=16pt, minimum width=42pt, xshift=-5] at (G) {};
\node[draw, rectangle, rounded corners, inner sep=16pt, minimum height=74pt, yshift=-22, xshift=1] at (R) {};
\node[draw, rectangle, rounded corners, minimum width=90pt, minimum height=74pt, yshift=-22, xshift=29] at (Tr) {};
\end{tikzpicture}
\end{minipage}
\begin{minipage}{.34\textwidth}

\vspace*{1.1cm}
\begin{tikzpicture}[,-latex ,auto ,node distance =2 cm and 2cm ,on grid , semithick ,
state/.style ={ rounded corners ,top color =white , bottom color = gray!20 ,
draw , minimum width =0.7cm, inner sep = 8pt}]

\node[state] (F) {$\cforest$};
\node[state] (T) [below =of F,align=left] {$\redroot\quad\blueroot$\\[6pt]$\greenroot$};

\node[state] (G)[right =of T,yshift=-10pt] {$\cgreen$};
\node[state] (RB)[left =of T,align=center] {$\cred$\\[6pt]$\cblue$};

\path (F) edge  node[left] {} (T);

\path (T) edge  node[left] {} (RB);
\path (T) edge  node[left] {} (G);

\end{tikzpicture}
\end{minipage}
}
\end{example}


\section{Well-Founded Systems of Generating Function Equations}\label{sec:generating-functions}

A combinatorial system translates into a system $\bs Y=\bs H(z,\bs Y)$ over generating functions. In the case of constructive systems, this is direct using the dictionary presented in \cref{tab:sum_esp_sg}. A few examples help understand the type of systems considered in this section.

\begin{example}The system of combinatorial equations from \cref{ex:colored_forest} translates into the system
\begin{multline}\label{eq:gf-colored-forests}
\left\{G(z)=z+G(z)^2,\;B(z)=\frac z{1-R(z)},\;R(z)=z^3+\frac z{1-B(z)},\right.\\
\left.T_r(z)=z + \frac{R(z)}{1 - T_b(z)},\;T_b(z)=\frac{B(z)}{1 - T_g(z)},\; T_g(z)=\frac{G(z)}{1 - \textbf{}T_r(z)},\;F(z)=\exp(T_r(z))\right\},
\end{multline}
where the equations here are sorted by connected components of the system, following the graph in \cref{ex:dependency_graph}.
\end{example}

\begin{example}Functional graphs~\cite[II.5.2]{FlajoletSedgewick2009} decompose as sets of connected components, each consisting of a cycle of trees. This leads to the system of generating function equations
\[\left\{F(z)=\exp(K(z)),\;K(z)=\ln\frac1{1-T(z)},\;T(z)=z\exp(T(z))\right\}.\]
\end{example}

More generally, we give a definition of well-founded analytic systems that reflects the combinatorial conditions of \cref{thm:gl-implicit} and extends similar definitions that were considered in the literature~\cite{Drmota1997,Drmota2009,BellBurrisYeats2006,BellBurrisYeats2010,PivoteauSalvySoria2012}. 
This section then develops the mathematical tools that are used in the development and proof of the algorithms operating on these well-founded systems in the rest of this article.

A difficulty is that for any given~$z\in\mathbb C$, such a system of generating function equations may have many solutions $\bs Y\in\mathbb C^m$ and the only solution of interest is the one that corresponds to the value of the generating function, which is the one we denote by $\bs{Y}(z)$. Even restricting to solutions with nonnegative coordinates at a real~$a>0$ may still result in several possible choices. Still, \cref{thm:alerho} shows that when such solutions exist, then $a$ is at most the radius of convergence of the generating function. This is the basis for a dichotomy that reduces the computability of the radius of convergence to the existence of nonnegative real solutions of systems of equations in \cref{subsec:radius}. 

The computation of the radius of convergence~---~indeed, much of the analysis of these systems~---~reduces to the special case of irreducible systems. When these are polynomial systems, their analysis is summarized by Flajolet and Sedgewick~\cite[Ch.~VII]{FlajoletSedgewick2009}; a more general situation is handled by Drmota~\cite{Drmota1997,Drmota2009}. A key insight of Bell, Burris and Yeats~\cite{BellBurrisYeats2010} is that one can distinguish which of the nonnegative real solutions of an irreducible system corresponds to the value of the generating function by considering the dominant eigenvalue of the Jacobian matrix of~$\bs H$ with respect to~$\bs Y$. Their hypotheses are slightly more restrictive than those used here, but most of their proofs extend to our more general situation (\cref{thm:Lambda}). 

\begin{notation}[for power series] This section focuses on generating functions in one variable~$z$, which explains why the boldfaced~$\bs z$ is not used. In order to insist on the values taken by the generating functions at real points, sometimes $\bs Y(x)$ is used instead of $\bs Y(z)$. 
The ring of formal power series with real coefficients is denoted~$\mathbb R[[z]]$ and similarly $\mathbb C[[z]]$ when the coefficients are complex. Their field of fractions are denoted $\mathbb R((z))$ and $\mathbb C((z))$ and their elements are called formal Laurent series. 
The subrings $\mathbb R\{z
\}$ and $\mathbb C\{z\}$ denote the convergent power series, while the subfields $\mathbb R\{\{z\}\}$ and $\mathbb C\{\{z\}\}$ denote power series that are convergent in a punctured neighborhood of~0. When $r$ is a positive integer, changing $z$ into~$z^{1/r}$ gives the corresponding rings and fields of Puiseux series. The elements of $\mathbb R\{\{z^{1/r}\}\}$ and $\mathbb C\{\{z^{1/r}\}\}$ converge in a slit neighborhood of~0. Finally, for $\sigma\in\mathbb C\setminus\{0\}$, changing~$z$ into~$1-z/\sigma$ gives the rings and fields used to express expansions in the neighborhood of~$\sigma$.
\end{notation}

\subsection{Well-founded Systems}\label{subsec:WFsys}

\subsubsection{Definition}
Systems of generating function equations that are well founded are defined in analogy with~\cref{thm:gl-implicit}.

\begin{definition}\label{def:wf-analytic}
Let $\bs H_{1:m}:\mathbb{C}^{m+1}\rightarrow\mathbb{C}^m$ be analytic in a neighborhood of~$\bs 0$. 
The system $\bs Y=\bs H(z,\bs Y)$ is called \emph{well founded} when
\begin{enumerate}
    \item each coordinate of $\bs H$ has nonnegative real Taylor coefficients at $\bs 0$; 
    \item letting $\bs U^{[0]}:=\bs 0$, for $k=0,\dots,m$, $(0,\bs U^{[k]})$ belongs to the domain of convergence of the Taylor series of~$\bs H$ at $\bs 0$, allowing one to define $\bs U^{[k+1]}:=\bs H(0,\bs U^{[k]})$;
     \item\label{def:wf-an-part3} the matrix $\bpartial\bs H/\bpartial\bs Y$ is nilpotent at $(0,\bs U^{[m]})$.
\end{enumerate}
\end{definition}
\noindent In the common case when $\bs H(0,\bs 0)=\bs 0$, condition~2 is automatically satisfied and $\bs U^{[m]}=0$.

\begin{definition}\label{def-irred}
The well-founded system~$\bs Y=\bs H(z,\bs Y)$ is called \emph{recursive} when $\bpartial\bs H/\bpartial\bs Y$ is not nilpotent and \emph{irreducible} when moreover it is well-founded and zero-free and its dependency graph (defined as in \cref{def:dependency-graph}) is strongly connected. This is equivalent to the irreducibility of the matrix $\bpartial\bs H/\bpartial\bs Y$ in the sense used in the Perron-Frobenius theory. 
\end{definition}
\Cref{def:wf-analytic} is analogous to \cref{def:combi_irreducible} in the combinatorial situation.
Similarly, nonrecursive and irreducible \emph{components} 
are defined in the same way as in \cref{def:irreducible-component-com}.
Also, the system is called \emph{constructive} when it is the translation of a constructive system in the sense of \cref{def:constructive} by the rules of \cref{tab:sum_esp_sg}.

\subsubsection{Solutions}

This definition of well-founded systems is motivated by the following analytic analogue of \cref{thm:gl-implicit}. 
\begin{proposition}\label{prop:well-founded-analytic}
If $\bs Y=\bs H(z,\bs Y)$ is a well-founded system, then it admits a unique solution~$\bs S$ analytic in a neighborhood of~$0$ with $\bs S(0)=\bs U^{[m]}$ given in \cref{def:wf-analytic}. This solution has nonnegative Taylor coefficients at~$0$.
\end{proposition}
\begin{definition}
The solution $\bs S$ is called the \emph{generating function solution} of the system. When the system is constructive, it is called a \emph{constructible generating function}.
\end{definition}
\begin{proof}
If $\bs H(0,\bs 0)=\bs 0$, then $\bs U^{[m]}=0$ and the existence of an analytic solution with $\bs S(0)=\bs 0$ is granted by the implicit function theorem at~$(0,\bs 0)$ (see, e.g., \cite[Chap.~IV]{Cartan1995}). Nonnegativity of the Taylor coefficients of $\bs S$ comes from the fixed-point iteration $\bs Y\mapsto\bs H(z,\bs Y)$ over power series, starting 
with~$\bs 0$.

In general, the implicit function theorem at~$(0,\bs 0)$ shows that the system $\bs Y=\bs H(z,\bs U^{[m]}+\bs Y)-\bs H(0,\bs U^{[m]})$ admits an analytic solution~$\tilde{\bs S}$ with~$\tilde{\bs S}(0)=\bs 0$ and again nonnegative Taylor coefficients. If $\bs U^{[m+1]}=\bs U_m$, the function~$\bs S=\bs U^{[m]}+\tilde{\bs S}$ is then a solution of~$\bs Y=\bs H(z,\bs Y)$ with nonnegative Taylor coefficients and~$\bs S(0)=\bs U^{[m]}$. 

It remains to show that~$\bs U^{[m+1]}\textbf{}=\bs U^{[m]}$. 
The proof follows the steps of that of \cref{thm:gl-implicit}. All inequalities between vectors mean that they hold for all coordinates.

The sequence~$(\bs U^{[k]})_{k\ge0}$ satisfies $\bs U^{[k]}\le\bs U^{[k+1]}$. For $k=0$, this follows from the nonnegativity of the Taylor coefficients of~$\bs H$ and for higher~$k$ by the fact that this nonnegativity makes each coordinate increasing. Next, by the lemma below,
\[
\bs U^{[k+1]}-\bs U^{[k]}\le \bpartial\bs H/\bpartial\bs Y(0,\bs U^{[k]})\times (\bs U^{[k]}-\bs U^{[k-1]}),\quad k>0.\]
Iterating and using $\bs U^{[k]}\le\bs U^{[k+1]}$ gives
\[\bs U^{[m+1]}-\bs U^{[m]}\le \left(\bpartial\bs H/\bpartial\bs Y(0,\bs U^{[m]})\right)^m\times(\bs U^{[1]}-\bs U^{[0]})=\bs 0.\qedhere\]
\end{proof}

\begin{lemma}\label{lemma:convex}
Let $\bs H(z,\bs Y)$ be as in \cref{def:wf-analytic} and let $(z,\bs A)$ and $(z,\bs B)$ belong to the domain of convergence of~$\bs H$ with $\bs A\ge \bs B$. Then
\[\frac{\bs\partial\bs H}{\bs\partial\bs Y}(z,\bs B)(\bs A-\bs B)\le \bs H(z,\bs A)-\bs H(z,\bs B)\le \frac{\bs\partial\bs H}{\bs\partial\bs Y}(z,\bs A)(\bs A-\bs B).\]
\end{lemma}
\begin{proof}
Since $\bs H$ has nonnegative Taylor coefficients at $\bs 0$, then so do all its derivatives, which implies in particular that $\bs H$ is convex on the intersection of its domain of convergence with~$\mathbb R^{m+1}$. Both inequalities are  consequences of that convexity.
\end{proof}

\subsubsection{Comparison with Previous Definitions}\label{sec:comparison}
Drmota~\cite{Drmota1997,Drmota2009} requires $\bs H(0,\bs Y)=\bs 0$ and thus $\bpartial\bs H/\bpartial\bs Y(0,\bs 0)=\bs 0$, but also that $\bs H(z,\bs 0)\neq\bs 0$ and the solution has no zero-coordinates. 
Bell, Burris, Yeats~\cite{BellBurrisYeats2010,BellBurrisYeats2011} have two definitions: \emph{elementary} systems of power series equations whose dependency graph is not necessarily strongly connected and where it is only required that $\bs H(0,\bs 0)=\bs 0$ and $\bpartial\bs H/\bpartial\bs Y(0,\bs 0)=\bs 0$; and \emph{well-conditioned} systems with Drmota's conditions and where moreover
$\bs H(z,\bs 0)$ should not have zero-coordinates, but with no condition on the solution. This is summarized in~\cref{table:system_types}.
Despite these variations, some of their results can be used in our context, by the following reductions.

\begin{table}
\renewcommand{\arraystretch}{1.2}
\setlength{\tabcolsep}{5pt}
\begin{tabularx}{\linewidth}{ |>{\raggedright}p{10em}|>{\raggedright}p{5.7em}|>{\raggedright}p{12.5em} |>{\RaggedRight}X| } 
\hline
Bell, Burris, Yeats, \qquad {\em well conditioned}~\cite{BellBurrisYeats2010}& 
$\bs H(0,\bs Y)=\bs 0$ & 
$\bpartial\bs H/\bpartial\bs Y(0,\bs 0)=\bs 0$ \\ Strongly connected graph  & 
$\bs H(z,\bs 0)$ has no zero-coordinates\\ 
\hline
Drmota~\cite{Drmota1997,Drmota2009} & 
$\bs H(0,\bs Y)=\bs 0$ & 
$\bpartial\bs H/\bpartial\bs Y(0,\bs 0)=\bs 0$ \\ Strongly connected graph  & 
Zero-free system and $\bs H(z,\bs 0)\neq\bs 0$ \\
\hline
Bell, Burris, Yeats, \qquad {\em elementary}~\cite{BellBurrisYeats2011} & 
$\bs H(0,\bs 0)=\bs 0$ & 
$\bpartial\bs H/\bpartial\bs Y(0,\bs 0)=\bs 0$ \\  & 

\\ 
\hline
Pivoteau, Salvy, Soria, \qquad \qquad {\em well founded at 0}~\cite{PivoteauSalvySoria2012} & 
$\bs H(0,\bs 0)=\bs 0$ & 
$\bpartial\bs H/\bpartial\bs Y(0,\bs 0)$ nilpotent\\  & 
Zero-free system\\ 
\hline
Pivoteau, Salvy, Soria, \qquad \qquad {\em well founded}~\cite{PivoteauSalvySoria2012} & 
\multicolumn{3}{p{11cm}|}{Let $\bc K= \bc H(\bc Z,\bc Y)-\bc H(\bs{0},\bs{0})+{\mathcal Z}_1\bc H(\bs{0},\bs{0})$; \,$\bc Y=\bc K({\mathcal Z}_1, \bc Z,\bc Y)$  well founded at 0, plus a finiteness condition on its solution}
\\ 
\hline
Banderier, Drmota,\quad\quad{\em analytically well defined}~\cite{BanderierDrmota2015}&
$\bs H(0,\bs 0)=\bs 0$
&$\bpartial\bs H/\bpartial\bs Y(0,\bs 0)$ has spectral radius~$<1$&$\bs H$ entire, no polynomial coordinate in~$\bs Y$
\\ 
\hline
This article,\qquad \qquad  {\em well founded} & 
$\bs H^{[m]}(0,0)$ well defined & 
$\bpartial\bs H/\bpartial\bs Y(0,\bs Y(0))$ nilpotent \\ & 
\\ 
\hline
\end{tabularx}
\caption{Conditions on the system~$\bs Y=\bs H_{1:m}(z,\bs Y)$ in previous works. In all cases, each coordinate of $\bs H$ has nonnegative Taylor coefficients at $\bs 0$. 
}
\label{table:system_types}
\end{table}

\paragraph{Zero coordinates at $\bs Y=\bs 0$ and zero-free systems.}\label{def:0free} 
As stated earlier, these coordinates can be detected algorithmically (using Algorithm \textsf{\ref{algo:iswf}}) and removed from the system.

\paragraph{Nonzero nilpotent Jacobian matrix.}
The case when $\bpartial\bs H/\bpartial\bs Y(0,\bs 0)$ is nilpotent but not~$\bs 0$ reduces to the case when it is~$\bs 0$~\cite[Section~4.4]{BellBurrisYeats2011}:  writing $\bs H=\bpartial\bs H/\bpartial\bs Y(0,\bs 0)\bs Y+\hat{\bs H}$ leads to a new well-founded system for the same~$\bs Y$,
\[\bs Y=\left(\Id-\frac{\bpartial\bs H}{\bpartial\bs Y}(0,\bs 0)\right)^{-1}\hat{\bs H}(z,\bs Y),\]
where that property holds. Note that if $\bs H(0,\bs Y)$ is the constant $\bs 0$, then $\bpartial\bs H/\bpartial\bs Y(0,\bs 0)=\bs 0$ and then $\hat{\bs H}(0,\bs Y)=\bs0$, which shows that this property is preserved by the reduction.

\paragraph{Nonzero value at $(0,0)$.}
Similarly, if $\bs H(0,\bs 0)\neq\bs 0$, a new system can be formed for~$\hat{\bs Y}(z):=\bs Y(z)-\bs Y(0)$:
\[\hat{\bs Y}=\hat{\bs H}(z,\hat{\bs Y}):=\bs H(z,\bs Y(0)+\hat{\bs Y})-\bs Y(0).\]
Recall that $\bs Y(0)$ is obtained as~$\bs U^{[m]}$.

This system satisfies $\hat{\bs H}(0,\bs 0)=\bs 0$, which follows from $\bs H(0,\bs Y(0))=\bs Y(0)$. It is well founded: 
the Taylor coefficients of $\hat{\bs H}$ are nonnegative, which can be seen by expansion of~$\bs H$; $\bpartial\hat{\bs H}/\bpartial\hat{\bs Y}(0,\bs 0)=\bpartial\bs H/\bpartial\bs Y(0,\bs Y(0))$ is nilpotent.

\paragraph{Nonzero value at $z=0$.} This case is a generalization of the previous one. It contains in particular combinatorial systems that put weights on leaves rather than internal nodes, such as binary trees defined by~$y=z+y^2$. 

A new system can be formed for~$\hat{\bs Y}$ defined by
\[\bs Y(z)=\bs Y(0)+z\bs Y'(0)+z\hat{\bs Y}(z),\quad\text{with $\hat{\bs Y}(0)=0$}.\]
As before, $Y(0)$ is obtained as~$\bs U^{[m]}$ and from there $\bs Y'(0)$ as the solution of the linear system
\[\bs Y'(0)=\frac{\partial\bs H}{\partial z}(0,\bs Y(0))+\frac{\bs\partial\bs H}{\bs\partial\bs Y}(0,\bs Y(0))\bs Y'(0).\]
Then, $\hat{\bs Y}$ satisfies
\[\hat{\bs Y}=\hat{\bs H}\!(z,\hat{\bs Y}):=\frac1z\left(\bs H\left(z,\bs Y(0)+z\bs Y'(0)+z\hat{\bs Y}\right)
-\bs Y(0)-z\bs Y'(0)\right).\]
By Taylor expansion, $\hat{\bs H}(0,\hat{\bs Y})=0$. Moreover, this new system is well founded: as before, the Taylor coefficients of~$\hat{\bs H}$ are nonnegative; the Jacobian matrix is
\[\frac{\bpartial\hat{\bs H}}{\bpartial \hat{\bs Y}}(z,\hat{\bs Y})=\frac{\bpartial\bs H}{\bpartial \bs Y}\!\left(z,\bs Y(0)+z\bs Y'(0)+z\hat{\bs Y}\right),\]
whose value at $(0,\bs 0)$ is $\bpartial\bs H/\bpartial\bs Y(0,\bs Y(0))$, which is nilpotent.

\subsubsection{Irreducible Systems}
The following lemma ensures that results on  irreducible systems apply to irreducible components of well-founded systems.
\begin{lemma}\label{lemma:wf-irreducible-components}
Let $\bs Y=\bs H(z,\bs Y)$ be a well-founded system that is not irreducible and splits into two sub-systems
\[\bs Y_{1:i}=\bs H_{1:i}(z,\bs Y_{1:i}),\quad\bs Y_{i+1:m}=\bs H_{i+1:m}(z,\bs Y_{1:i},\bs Y_{i+1:m})\]
where the second one is an irreducible component.
Let $\bs S_{1:i}$ be the solution of the first component given by \cref{prop:well-founded-analytic}. Then the system
\[\bs Y_{i+1:m}=\bs H_{i+1:m}(z,\bs S_{1:i},\bs Y_{i+1:m})=:\hat{\bs H}_{i+1:m}(z,\bs Y_{i+1:m})\]
is an irreducible system. 
\end{lemma}
\begin{proof}
\newcommand\w{{}_{\phantom{x}}}
By definition, $\bs Y_{i+1:m}=\hat{\bs H}_{i+1:m}(z,\bs Y_{i+1:m})$ is an irreducible system.
It is also well founded: (1) $\hat{\bs H}_{i+1:m}$ has nonnegative Taylor coefficients at $\bs 0$ by composition of $\bs H_{i+1:m}(z,\bs Y)$ with the generating functions $\bs Y_{1:i}$;
(2) recall from~\cref{def:wf-analytic},
that $\bs U\w^{[0]} =\bs 0$ and for $k=0,\dotsc,m$, $\bs U\w^{[k+1]}=\bs H(0,\bs U\w^{[k]})$. 
The corresponding iteration for $\hat{\bs H}$ at $z=0$ is defined by $\hat{\bs U}\w^{[0]}=\bs 0$  and for $k\leq m-i$,  $\hat{\bs U}\w^{[k+1]}={\bs H}_{i+1:m}(0,\bs U_{1:i}^{[m]},\hat{\bs U}\w^{[k]})$. By positivity, for $k\leq m-i$, $\hat{\bs U}\w^{[k]} \leq \bs U_{i+1:m}^{[m]}$. The properties required by~\cref{def:wf-analytic} can now be checked.
Indeed, as $(0,\bs U\w^{[m]})$ belongs to the domain of convergence of the Taylor series of~$\bs H$ at $\bs 0$, so does $(0,\hat{\bs U}\w^{[k]})$ for~$\hat{\bs H}$ at~$\bs 0$, for $k\leq m-i$;
(3) the bottom-right block of size  $(m-i)\times(m-i)$ in $\bs\partial\bs H/\bs\partial\bs Y(0,\bs U\w^{[m]})$ is precisely $\bs\partial\hat{\bs H}/\bs\partial\bs Y_{i+1:m}(0,{\bs U}_{i+1:m}^{[m]})$,  while the bottom-left block is zero. Since $\bs\partial\bs H/\bs\partial\bs Y(0,\bs U\w^{[m]})$ is nilpotent, so is its bottom-right block. Thus 
\[
0 \le (\bs\partial{\hat{\bs H}}/\bs\partial\bs Y_{i+1:m})^{m-i}(0,\hat{\bs U}{}^{[m-i]}) \le (\bs\partial\hat{\bs H}/\bs\partial\bs Y_{i+1:m})^{m-i}(0,\bs U^{[m]})=0.\qedhere
\]
\end{proof}

\subsection{Characteristic System and Dominant Eigenvalues}
\label{sec:ift}
We consider a well-founded system $\bs Y=\bs H(z,\bs Y)$ and denote by $\ini=\bs{Y}(0)$ the value of the generating function at~0, which is given by \cref{prop:well-founded-analytic}.
The implicit function theorem implies that the solution~$\bs Y$ is analytic on the positive real axis as long as $(z,\bs{Y}(z))$ belongs to the domain of convergence of $\bs H$ and is not a solution of the characteristic system in the following sense.
\begin{definition}\label{def:characteristic-system}
The \emph{characteristic system} of the system~$\bs Y=\bs H(z,\bs Y)$ is the extended system
\begin{equation}\label{eq:characteristic-system}
\bs Y=\bs H(z,\bs Y),\quad \det\!\left(\Id-\frac{\bpartial\bs H}{\bpartial\bs Y}\right)=0.
\end{equation}
\end{definition}
The characteristic system may possess several solutions with positive coordinates. An example due to Burris is
\begin{equation}\label{ex:Burris}
\{y_1=z(1+y_1^2+y_2),y_2=z(1+y_1+y_2^2)\},
\end{equation}
whose characteristic system has both $(z,y_1,y_2)=(1/3,1,1)$ and $(z,y_1,y_2)=((2\sqrt{2}-1)/7,\sqrt{2}+1,\sqrt{2}+1)$ as 
solutions~\cite[Ex.VII.28]{FlajoletSedgewick2009} and only the first one, with larger value of~$z$, is of interest. Thus, even if an algorithm converges to a solution of the characteristic system with positive coordinates, more work is required in order to prove that this solution corresponds to a dominant singularity of the generating function solutions. 

An idea due to Bell, Burris and Yeats~\cite{BellBurrisYeats2010} solves this problem for irreducible systems. We build on their study to design our algorithms.
For a well-founded system, at any point $(x,y_1,\dots,y_m)\in\mathbb{R}_{\ge0}^{m+1}$ inside the domain of convergence of $\bs H$, the Jacobian matrix $\bpartial\bs H/\bpartial\bs Y$ has nonnegative entries. By the Perron-Frobenius theorem, this implies that its spectral radius~$\Lambda_{\bs H}(x,\bs y)$ is also its largest real eigenvalue, that is called its \emph{dominant} eigenvalue. In particular, the dominant eigenvalue gives a criterion to check that a solution of the system $\bs Y=\bs H(a,\bs Y)$ is the value of the generating function at~$a\in\mathbb R_{\ge0}$ for irreducible systems. This will be a building block for our algorithm for the general case in the next section.

\begin{restatable}{theorem}{thLambda}\label{thm:Lambda}
Let $\bs Y=\bs H(z,\bs Y)$ be an \emph{irreducible} system. Let $(a,\bs B)\in{\mathbb R_{\ge0}^{m+1}}$ belong to the domain of convergence of $\bs H$ or its boundary and satisfy  $\bs B=\bs H(a,\bs B)$. Let~$\rho$ be the radius of convergence of the generating function solution $\bs Y$. Then, $\rho>0$ and 
\begin{enumerate}
    \item if $\bs H$ is linear in the variable $\bs Y$ then $a\le\rho$ and $\bs B=\bs Y(a)$;
    \item in all cases \label{part2}
    \begin{enumerate}
        \item\label{part2b} $\Lambda_{\bs H}(a,\bs B)\le1\Leftrightarrow a\le\rho$ and $\bs B=\bs Y(a)$;
        \item\label{part2a} $\Lambda_{\bs H}(a,\bs B)=1\Rightarrow (a,\bs B)=(\rho,\bs Y(\rho))$.
    \end{enumerate}
\end{enumerate}
\end{restatable}

\begin{notation}For $0\le a\le\rho$, we write $\Lambda_{\bs H}(a)$ for $\Lambda_{\bs H}(a,\bs Y(a))$.
\end{notation}
\smallskip

An application of Part \ref{part2} of this theorem is to filter out solutions of the characteristic system~\eqref{eq:characteristic-system} that do not correspond to a dominant singularity. This was the motivation for the work of Bell, Burris and Yeats, who use slightly stronger hypotheses. An adaptation of their proof to \cref{thm:Lambda} is given in \cref{appendix:prooftheorem2.3}.

\begin{corollary}\label{coro:lambda=1} Let $\bs Y=\bs H(z,\bs Y)$ be an irreducible system. Let $\rho$ be the radius of convergence of the generating function~$\bs Y$. Assume that $\bs\tau=\lim_{u\rightarrow\rho-}\bs Y(u)$ is finite. Then either $(\rho,\bs\tau)$ is on the boundary of the domain of convergence of $\bs H$ or $\Lambda_{\bs H}(\rho,\bs\tau)=1$.
\end{corollary}
\begin{proof}When $\bs \tau$ is finite, 
$\bs Y(\rho)=\bs\tau$. By continuity, $(\rho,\bs\tau)$ is either in the domain of convergence of $\bs H$ or on its boundary. \Cref{part2b} of the theorem implies that $\Lambda_{\bs H}(\rho,\bs Y(\rho))\le1$. If moreover $(\rho,\bs\tau)$ lies inside the domain of convergence of~$\bs H$, then by the implicit function theorem, the singularity occurs at a solution of the characteristic system~\eqref{eq:characteristic-system}, where $\Lambda_{\bs H}=1$.
\end{proof}

\begin{corollary}\label{coro:lambda-irr-comp}\cref{thm:Lambda} and \cref{coro:lambda=1}  hold more generally if the system is an irreducible component of a well-founded system.
\end{corollary}
\begin{proof}
    By \cref{lemma:wf-irreducible-components}.
\end{proof}

\subsection{Recognizing Points Inside the Domain of Convergence}
Using the decomposition of general systems into irreducible components, it is possible to ensure that a point is inside the domain of convergence of the generating function solution of a well-founded system that is not necessarily irreducible, by the following. 

\begin{restatable}{proposition}{alerho}\label{thm:alerho}
Let $\bs Y=\bs H(z,\bs Y)$ be a zero-free well-founded system, $\mathcal D(\bs H)\subset\mathbb C^{m+1}$  the domain of convergence of~$\bs H$ and~$\rho$ the radius of convergence of the generating function solution~$\bs Y$. 
If there exists $(a,\bs B)\in\mathbb R_{\ge0}^{m+1}\cap\mathcal D(\bs H)$ such that $\bs B=\bs H(a,\bs B)$, then $a\le\rho$, $\bs Y(a)$ exists and $\bs B\ge\bs Y(a)$.
\end{restatable}

Note however that, except in the linear case, one cannot conclude from the hypotheses that $\bs B=\bs Y(a)$. For instance, the equation for binary trees $Y=z+Y^2$ has the positive solution $B=1$ at~$z=0$. This is not equal to $Y(0)$, which is~0. Since this system is irreducible, by \cref{thm:Lambda}, detecting that $\bs B=\bs Y(a)$ is achieved by checking that $\Lambda_{\bs H}(a,\bs B)\le1$. In this example, $\Lambda(0,0)=0\le1$ and $\Lambda(0,1)=2>1$. 

In general, the inequality $\bs B\ge\bs Y(a)$ is similar to the classical Kleene fixed point theorem. The main use of this proposition however is that it ensures that the system does not have nonnegative solutions beyond the radius of convergence of the generating function solution~$\bs Y$.

\medskip

\noindent As the proof of \cref{thm:alerho} uses several of the lemmas proved in \cref{appendix:prooftheorem2.3}, it is postponed to \cref{appendix:proofprop3.11}.

\subsection{Radius of Convergence of Constructible Generating Functions}\label{subsec:radius}
The results on general well-founded systems are now applied to the constructive systems of Flajolet and Sedgewick (\cref{def:constructive}).
\subsubsection{Finiteness}\label{sec:finiteness}
It is easy to detect when the radius of convergence is infinite.
\begin{proposition}\label{prop:finiteness}
    The radius of convergence of the constructible generating function solution to~$\bs Y=\bs H_{1:m}(z,\bs Y)$ is infinite if and only if the system is not recursive and~$\bs H$ does not contain $\Seq$ or $\Cyc$.
\end{proposition} 
\begin{proof}
Let $\bs S$ be the generating function solution.
When the system is recursive and irreducible, $\bpartial\bs H/\bpartial\bs Y$ is not nilpotent while its value at $(0,\bs S(0))$ is. This implies that~$\partial\bs H/\partial z$ is nonzero at~$(z,\bs S(z))$ for~$z>0$. Since
\[\bs Y'(z)=\left(\Id-\frac{\bpartial\bs H}{\bpartial\bs Y}\right)^{-1}\frac{\partial\bs H}{\partial z},\]
it follows by evaluation at $(z,\bs S(z))$ that $\bs S'(z)$ is positive. Then $\bs S(z)$ is a nonconstant power series with nonnegative coefficients. If its radius of convergence were infinite, the series would tend to infinity as~$z\rightarrow\infty$. Then, so would the dominant eigenvalue $\Lambda(z,\bs Y(z))$, which is bounded by~1, a contradiction.

When the system is recursive and not irreducible, it has an irreducible component to which the previous reasoning applies. 

Finally, if the system is not recursive, then $\bs S$ is the limit of the iteration of power series $\bs Y^{[i+1]}=\bs H(z,\bs Y^{[i]})$ with $\bs Y^{[0]}=\bs 0$. By \cref{lemma:convex},
\[\bs Y^{[m+1]}-\bs Y^{[m]}\le\frac{\bpartial\bs H}{\bpartial\bs Y}(z,\bs Y^{[m]})(\bs Y^{[m]}-\bs Y^{[m-1]})
\le \left(\frac{\bpartial\bs H}{\bpartial\bs Y}(z,\bs Y^{[m]})\right)^{\!m}(\bs Y^{[1]}). 
\] 
Since the system is nonrecursive, the matrix $\bpartial\bs H/\bpartial\bs Y$ is nilpotent and the right-hand side is~$\bs 0$, which implies that $\bs S=\bs Y^{[m]}$, the $m$th iterate of~$\bs H$ at~$\bs 0$. Thus each coordinate of~$\bs S$ is an explicit combination of $1,z,z\mapsto\exp(z),z\mapsto1/(1-z),z\mapsto\ln1/(1-z)$ (the last two correspond to $\Seq$ and $\Cyc$) with additions, multiplications, and the constraints that $z\mapsto 1/(1-z)$ and $z\mapsto \ln1/(1-z)$ are applied to nonzero power series (since the system is zero free) of positive valuation. By induction such a function is entire if and only if neither $z\mapsto 1/(1-z)$ nor $z\mapsto \ln1/(1-z)$ occur in its expression. (See also \cite[Prop.~4.2]{FlajoletSalvyZimmermann1991}.)
\end{proof}

\subsubsection{Decidability}\label{sec:decidability}
\Cref{thm:alerho} gives a basis for an algorithmic characterization of the \textbf{}radius of convergence of generating functions. When the system~$\bs B=\bs H(a,\bs B)$ has a \emph{polynomial}~$\bs H$, then its domain of convergence is all of~$\mathbb C^{m+1}$ and the question of the existence of a solution in~$\mathbb R^{m+1}_{\ge 0}$ is decidable by the Tarski-Seidenberg theorem. For the constructive systems of \cref{def:constructive} that also involve the function~$x\mapsto1/(1-x)$, testing that a point belongs to the domain of convergence is also given by a polynomial inequality $x<1$ and thus also decidable. 

In general, constructive systems have generating function equations that also involve logarithms and exponentials (see \cref{tab:sum_esp_sg}). Provided the system is written in normal form (\cref{def:normal_form}), the hypothesis $(a,\bs B)\in\mathbb R_{\ge0}^{m+1}\cap \mathcal D(\bc H)$ is expressed by inequalities~$a\ge0$, $B_i\ge0$ for all~$i$ and~$B_i<1$ for the~$B_i$ that are arguments to the functions $x\mapsto1/(1-x)$, $x\mapsto\ln(1/(1-x))$. These inequalities, together with 
the existence of a solution of the system $\bs B=\bs H(a,\bs B)$ form a formula in the theory of the real numbers with the exponential. This theory was proved decidable by Macintyre and Wilkie~\cite{MacintyreWilkie1996}, under the assumption that the following classical conjecture in number theory holds.
\begin{conjecture}[Schanuel] If $x_1,\dots,x_n$ are real numbers that are linearly independent over~$\mathbb{Q}$, then the field $\mathbb{Q}(x_1,\dots,x_n,e^{x_1},\dots,e^{x_n})$ has transcendence degree at least~$n$ over~$\mathbb{Q}$.
\end{conjecture}
A direct consequence of these strong results from mathematical logic is the following.
\begin{theorem}\label{th:radius-is-computable}
If the radius of convergence of a constructible generating function is finite, then it is a computable number, provided either that its defining system does not contain $\Cyc{}$ or $\Set$, or that Schanuel's conjecture holds. 
\end{theorem}
\begin{proof}
Recall that by \cref{def:constructive}, constructive systems are zero-free.
A simple process sufficient for the computability is a search by dichotomy, using \cref{thm:alerho} that reduces the decision $a\le\rho$ to the existence of nonnegative solutions of the system $\bs B=\bs H(a,\bs B)$ in the domain of convergence of~$\bs H$.
\end{proof}

\subsubsection{Constant Oracle}\label{sec:constant-oracle}

The existence problems that must be decided in the use of \cref{thm:alerho} are of a specific form and might well be decidable in a more direct way than by the result of Macintyre and Wilkie. 
For this reason, we state our algorithms in terms of a more restricted oracle.

\begin{definition}\label{def:oracle}
The \textsf{Constant Oracle} is an oracle that takes as input functions~$\bs H:(z,\bs Y_{1:m})\mapsto \bs H_{1:m}$ with $H_k$s of one of the forms~$\{1,z,Y_i+Y_j,Y_i\times Y_j,1/(1-Y_i),\ln 1/(1-Y_i),\exp(Y_i)\}$ ($1\le i,j\le m$) and the inequalities~$Y_i\ge0$ for all~$i=1,\dots,m$, and $Y_i<1$ for the coordinates involved in~$1/(1-Y_i),\ln 1/(1-Y_i)$. It is able to decide:
\begin{enumerate}
    \item given~$z\ge0,\bs H$ and the inequalities, whether the system~$\bs Y=\bs H(z,\bs Y)$ has a solution satisfying the inequalities;
    \item if so, whether it has a solution with $\Lambda(z,\bs Y)=1$;
    \item given two such systems $\bs Y_{\!a}=\bs H_{\!a}(z_a,\bs Y_{\!a})$ and $\bs Y_{\!b}=\bs H_{\!b}(z_b,\bs Y_{\!b})$ with solutions that satisfy their inequalities, which of $z_a<z_b,z_a=z_b,z_a>z_b$ holds.
\end{enumerate}
\end{definition}
\noindent A weaker oracle was used in \cite{FlajoletSalvyZimmermann1991}, where only nonrecursive equations were considered.
The following is a reformulation of \cref{th:radius-is-computable}.
\begin{proposition}\label{prop:oracle-gives-radius}
If the radius of convergence of a constructible generating function is finite, then it can be computed with the help of the \textsf{Constant Oracle}.
\end{proposition}

Again, in the case of context-free systems, i.e., when the combinatorial system contains neither~$\Cyc$ nor $\Set$, all the numbers involved are real algebraic numbers and the \textsf{Constant Oracle} is computable. Efficient algorithms and implementations are available, such as Safey El Din's \texttt{RagLIB} package.\footnote{Available at~\url{https://www-polsys.lip6.fr/~safey/RAGLib/}.}

In practice, in the more general situation, one can use numerical approximations and assume that if two constants agree for their first few hundred digits, then they coincide. This will very likely be satisfactory unless the system has been designed for the purpose of countering this method.

\section{Newton Iteration for Well-Founded Systems}
\label{sec:Newton}

The classical Newton iteration for solving  equations or systems of equations numerically is particularly powerful when applied to combinatorial systems, whose positivity allows for more control over the convergence than in a general situation. For $a\ge0$ inside the disk of convergence of the generating functions, when started at~$\bs 0$, Newton's iteration always converges to the values of the generating functions at~$a$~\cite{PivoteauSalvySoria2012}. A converse holds: if Newton's iteration started at~$\bs 0$ converges, then $a\le\rho$, where $\rho$ is the radius of convergence (\cref{subsec:convergence-newton}). Positive systems with polynomial entries have been considered in detail by Etessami, Yannakakis, Esparza, Kiefer, Luttenberger~\cite{EtessamiYannakakis2009,EsparzaKieferLuttenberger2010,EtessamiStewartYannakakis2013,StewartEtessamiYannakakis2015}. Their algorithmic approach is similar to ours and faces the same difficulties: if the iteration is stopped prematurely, it could be that $a>\rho$ while the divergence of the iteration has not manifested itself yet. For this situation, we give a Newton-Kantorovich-type \emph{a posteriori} bound on the output of the iteration that can be used to certify the convergence (\cref{subsec:Kantorovich}). 
The positivity of combinatorial systems allows to derive such a bound from only one evaluation, instead of the required bound on a whole disk in the general situation. An analysis of the speed of convergence in the spirit of~\cite{EsparzaKieferLuttenberger2010} is also possible (\cref{subsec:speed}).

With these results, one can compute radii of convergence by dichotomy, but this is not as fast as a direct use of Newton's iteration.
For irreducible systems, the radius of convergence is given as a solution of a characteristic system (\cref{def:characteristic-system}). Unfortunately, this is not a positive system and the unconditional convergence of Newton's iteration does not hold for these systems. Still, there is a neighborhood of the solution where the convergence is quadratic (\cref{sec:newton-radius-irred}). The same property holds for systems expressing that one of the coordinates is~1, that occur when looking for singularities of $\Seq$ and $\Cyc$ (\cref{subsec:newton-u-eq-1}).

\subsection{Convergence and Divergence of Newton's Iteration}\label{subsec:convergence-newton}
The following strengthens \anonymousvariant{our previous result}{the result of}~\cite[Thm.~9.13]{PivoteauSalvySoria2012} when~$a$ is real. 
It shows that for the values computed by Newton's iteration, one has more precise information than is provided in general by \cref{thm:alerho}.
\begin{theorem}\label{thm:Newton}
Let $\bs Y=\bs H_{1:m}(z,\bs Y)$ be a well-founded  system and let $\rho$ be the radius of convergence of the generating function~$\bs Y$ of the solution. For $a\ge0$, let $\bs y^{[n]}$ denote the sequence defined by $\bs y^{[0]}=\bs 0$ and
\begin{equation}\label{eq:Newton-gf-val}\bs y^{[n+1]}=\bs y^{[n]}+\left(\Id-\frac{\bpartial\bs H}{\bpartial\bs Y}(a,\bs y^{[n]})\right)^{\!-1}\cdot (\bs H(a,\bs y^{[n]})-\bs y^{[n]}),\qquad n\ge0,
\end{equation}
whenever $(a,\bs y^{[n]})$ belongs to the domain of convergence of~$\bs H$, and undefined otherwise.
Then,
\begin{itemize}
\item[---] If $a\le\rho$, then $(a,\bs y^{[n]})$ belongs to the domain of convergence of~$\bs H$ for all~$n$  and $\bs y^{[n]}$ is a strictly increasing sequence converging to $\bs Y(a)$;
\item[---] conversely, if $\bs y^{[n]}$ converges to a vector~$\bs B\in\mathbb R_{\ge0}^{m}$ such that $(a,\bs B)$ belongs to the domain of convergence of~$\bs H$ or its boundary, then $a\le\rho$ and $\bs B=\bs Y(a)$;
\item[---] if $a>\rho$, then $\bs y^{[n]}$ is an increasing sequence for $n=0,\dots,N$ until a value~$N$ such that
$(a,\bs y^{[N]})$ does not belong to the domain of convergence of~$\bs H$ or such that $\Lambda_{\bs H}(a,\bs y^{[N]})\ge1$.
\end{itemize}
\end{theorem}
\begin{proof}
The first part is \cite[Thm.~9.13]{PivoteauSalvySoria2012}, without the constraint of being zero-free on the system (which can be removed using \cref{thm:gl-implicit}). The second is a consequence of \cref{thm:alerho} and of the first part. 

For the last part, we first show that as long as $(a,\bs y^{[n]})$ belongs to the domain of convergence of $\bs H$ and $\Lambda_{\bs H}(a,\bs y^{[n]})<1$,  the sequence $\bs y^{[n]}$ increases. Indeed, by positivity of~$\bs H$, it follows that $\bs H(a,\bs 0)\ge\bs 0$ and by induction $\bs H(a,\bs y^{[n]})\ge\bs H(a,\bs y^{[n-1]})=\bs y^{[n]}$ showing that the last vector in \cref{eq:Newton-gf-val} has nonnegative coordinates. Next, since $\Lambda_{\bs H}(a,\bs y^{[n]})< 1$, the inverse matrix in the right-hand side of \cref{eq:Newton-gf-val} has nonnegative entries, showing that the sequence~$\bs y^{[n]}$ is increasing. 

Finally, assume by contradiction that $a>\rho$, that the sequence $\bs y^{[n]}$ is increasing, that for all~$n$, $(a,\bs y^{[n]})$ belongs to the domain of convergence of~$\bs H$ and that $\Lambda_{\bs H}(a,\bs y^{[n]})< 1$.
If the increasing sequence $\bs y^{[n]}$ were bounded, it would converge, but this is impossible by the second part of the theorem. Thus, one of the coordinates of $\bs y^{[n]}$ tends to~$+\infty$, implying that $\bs H$ and $\bpartial \bs H/\bpartial\bs Y$ depend on this coordinate. The monotonicity of $\Lambda_{\bs H}$ makes this incompatible with the fact that $\Lambda_{\bs H}(a,\bs y^{[n]})<1$ for all~$n$.
\end{proof}

\begin{method}[ht]
\SetAlgoRefName{NewtonIterationMethod} 
\caption{Newton Iteration Method\label{method:newton-method}}

\DontPrintSemicolon
\Input{$\bs{Y}=\bs{H}(z,\bs{Y})$ a well-founded system;\\
\qquad\qquad $a$ a nonnegative real number.}
\Output{An approximation of $\bs Y(a)$ if $a\le\rho$ or `$a>\rho$' in some cases when $a>\rho$\vspace{10pt}}
    Set $\bs y^{[0]}:=\bs 0$\;
    \For{$n=0,1,2,\dotsc$}{
        Compute $\bs y^{[n+1]}$ by Newton's iteration from \cref{eq:Newton-gf-val}\;
        \lIf{$(a,\bs y^{[n+1]})$ is not in the domain of convergence of $\bs H$ or\\\quad\,\,$\bs y^{[n+1]}\not\ge \bs y^{[n]}$ or \\\quad\,$\Lambda_\bs H(a,\bs y^{[n+1]})\ge1$}{\Return{`$a>\rho$'}}
        \lIf{$\bs y^{[n+1]}$ is sufficiently close to $\bs y^{[n]}$}{\Return $\bs y^{[n+1]}$}
    }
\end{method}

This theorem gives a practical criterion for detecting that a value is larger than the radius of convergence: the vector $\bs y^{[n]}$ computed by the Newton iteration will ultimately have a decreasing coordinate or be outside the domain of convergence of~$\bs H$ or be such that ${\Lambda_{\bs H}(a,\bs y^{[n]})>1}$. An explicit variant of the method deduced from this theorem is given in \textsf{\ref{method:newton-method}}.

\begin{remark} This method is very similar to the approach taken by Etessami and Yannakakis~\cite{EtessamiYannakakis2009} for positive polynomial systems, which it extends to the case of well-founded systems. Moreover,
for systems of polynomial equations coming from recursive Markov chains, they use a `decomposed Newton's method', where the decomposition of the system into strongly connected components is used to solve one system per component by Newton iteration along the dag of the decomposition. The above theorem shows that, in the context of well-founded combinatorial systems, this is not necessary for convergence. Still, our implementation uses the decomposition, as this leads to smaller Jacobian matrices.
\end{remark}

\begin{remark}A practical implementation of Newton's iteration~\eqref{eq:Newton-gf-val}  uses floating-point numbers. This introduces further issues due to rounding. While the convergence of Newton's iteration is known to be robust,  inequalities like increasing coordinates may fail to hold close to the radius of convergence unless some care is taken.
Thus, it is important to use directed rounding towards~0 or~$-\infty$ (these rounding modes are provided by the IEEE-754 standard) during the evaluation. In addition, it is necessary that the implementation of the evaluation of~$\bs H$ preserves the inequality $\bs H(a,\bs y^{[n]})\ge \bs y^{[n]}$ and the positivity of $\left(\Id-\frac{\bpartial\bs H}{\bpartial\bs Y}(a,\bs y^{[n]})\right)^{\!-1}$. In that case, the monotonicity of the rounding operation implies that one can stop the iteration when a fixed point is reached on floating-point values. We do not consider these issues further in this article.
\end{remark}

\subsection{A Posteriori Bounds}\label{subsec:Kantorovich}

\Cref{thm:Newton} gives information on the limit of the sequence of Newton iterations. In practice, when the iteration is stopped after a finite number of steps as in \textsf{\ref{method:newton-method}}, the theorem does not help determine whether $a\le\rho$ or not, except when the algorithm returns `$a>\rho$'.\footnote{This is the same problem as in~\cite[Cor.~7.5 and Footnote~6]{EtessamiYannakakis2009}.} Still, the positivity of our systems is a strong property, that makes it easier to apply the Newton-Kantorovich theorem: in a more general situation, that theorem requires control of the system in a whole neighborhood of the iterate, while here it is only necessary to perform one evaluation at one point in order to certify that $a<\rho$.
\begin{theorem}\label{thm:Kantorovich}
With the hypotheses and notation of \cref{thm:Newton}, let $r=2\|\bs y^{[n+1]}-\bs y^{[n]}\|_\infty$ and $\bs u=\bs y^{[n]}+\bs r$, where $\bs r=(r,\dots,r)$. If $(a,\bs u)$ belongs to the domain of convergence of~$\bs H$, $\bs y^{[n]}-\bs r\ge\bs 0$ and 
\begin{equation}\label{eq:kanto}
r\left\|\left(\Id-\frac{\bpartial\bs H}{\bpartial\bs Y}(a,\bs y^{[n]})\right)^{\!-1}
\begin{pmatrix}
\sum_{i,j}\frac{\partial^2 H_1}{\partial Y_i\partial Y_j}(\bs u)\\
\vdots\\
\sum_{i,j}\frac{\partial^2 H_m}{\partial Y_i\partial Y_j}(\bs u)\\
\end{pmatrix}\right\|_\infty\le1,\end{equation}
then $a<\rho$ and for all $k\ge0$,
\[\|\bs Y(a)-\bs y^{[n+k]}\|_\infty\le \frac{\kappa}{2^{2^k-1}}r,\quad\text{where}\quad \kappa=\sum_{i=0}^\infty{2^{-2^i}}< 0.8165.\]
\end{theorem}
In other words, when the conditions hold, the current iterate $\bs y^{[n]}$ is at distance less than twice the length of the Newton step $\bs y^{[n+1]}-\bs y^{[n]}$ from~$\bs Y(a)$ and moreover,  from that point on, the convergence is bounded by that of a quadratically convergent sequence.

\begin{proof}
The proof reduces to a variant of the Newton-Kantorovich theorem~\cite[\S3.2]{Dedieu2006}. In order to simplify notations, we write
\[\bs x_k:=\bs y^{[n+k]},\quad \bs F(\bs x):=\bs x-\bs H(a,\bs x)\]
and observe that $D\bs F=\Id-D\bs H$ and $D^2\bs F=-D^2\bs H$.

By positivity of~$\bs H$ and all its derivatives for nonnegative arguments inside the domain of convergence of~$\bs H$, for any $\bs x$ such that $\|\bs x-\bs x_0\|_\infty\le r$,  for each $i,j,\ell$ in $\{1,\dots,m\}$, 
\[\left|\frac{\partial^2H_\ell}{\partial Y_i\partial Y_j}(\bs x)\right|
\le \frac{\partial^2H_\ell}{\partial Y_i\partial Y_j}(\bs u).\]
For any two vectors $\bs U,\bs V$ in $\mathbb R^m$, 
\[D\bs F(\bs x_0)^{-1}D^2\bs F(\bs x):(U,V)\mapsto -(\Id-D\bs H(\bs x_0))^{-1} 
\begin{pmatrix}
    \bs U^\top\operatorname{Hess}(H_1(\bs x))\bs V\\
    \vdots\\
    \bs U^\top\operatorname{Hess}(H_m(\bs x))\bs V    
\end{pmatrix},\]
where $\operatorname{Hess}(H_i(\bs x))$ is the Hessian matrix of the $i$th coordinate $H_i$ of $\bs H$ at~$\bs x$.
As $(\Id-D\bs H(\bs x_0))^{-1}$ is a matrix with nonnegative entries, by triangular inequality, the product on the right-hand side is a vector whose coordinates are bilinear forms in~$\bs U,\bs V$ with coefficients that are upper bounded by those 
obtained by replacing~$\bs x$ by~$\bs u$ in its expression. 
The maximal value of the coordinates of this vector over all $\bs U,\bs V$ of norm~1 is reached at~$\bs U=\bs V=(1,\dots,1)$ and thus the norm in \cref{eq:kanto} is an upper bound on the norm of 
the operator $D\bs F(\bs x_0)^{-1}D^2\bs F(\bs x)$ for $\|\bs x-\bs x_0\|_\infty\le r$.

In this situation, the hypotheses of the Newton-Kantorovich theorem~\cite[Thm. 88]{Dedieu2006} hold. It follows that:
\begin{itemize}
\item[--] there exists a unique~$\bs x$ with $\|\bs x-\bs x_0\|_\infty\le r$ such that $\bs x=\bs H(a,\bs x)$; 
\item[--] Newton's iteration converges to it at the speed given in the last inequality;
\item[--] $\Id-\frac{\bpartial\bs H}{\bpartial\bs Y}(a,\bs x)$ is invertible. 
\end{itemize}
It remains to see that~$\bs x=\bs Y(a)$. Since $\bs x_0+\bs r$ is inside the domain of convergence of~$\bs H$, so is~$\bs x$. Also, since $\bs x_0-\bs r\ge\bs 0$, it follows that $(a,\bs x)\in\mathbb R_{\ge0}^{m+1}$. Then $\bs x=\bs Y(a)$ follows from the second part of \cref{thm:Newton}. That $a\neq\rho$ is a consequence of the invertibility of $\Id-\frac{\bpartial\bs H}{\bpartial\bs Y}(a,\bs Y(a))$.
\end{proof}
\begin{example}\label{ex:sing_big_spec} Let $\bc Y = \bc H(\cal Z, \bc Y)$ be the combinatorial system defined by~\cref{eq:colored_forest} in \cref{ex:colored_forest} and $\rho$~its radius of convergence. For $a=0.1703915$, after 13 steps of the iteration, the last step has length longer than~$10^{-5}$. Still, the vector is positive, its entries are larger than the step and the inequality \eqref{eq:kanto} becomes $0.7996822\le1$, which is sufficient to conclude that $a\le\rho$. Conversely, for $b=0.1703918$, Newton's iteration \eqref{eq:Newton-gf-val} at $z=b$ produces at least one decreasing coordinate, showing that $b>\rho$ by \cref{thm:Newton}. Thus~$\rho\in[a,b]$ and this can be refined by dichotomy. A direct Newton iteration is discussed in \cref{sec:newton-radius-irred}.
\end{example}


\subsection{Speed of Convergence}\label{subsec:speed}
Even when $a<\rho$ and the convergence of Newton's iteration becomes eventually quadratic as in \cref{thm:Kantorovich}, there is no clear control \emph{a priori} of the number of steps needed before that quadratic regime takes place. Parts of the analysis of Esparza \emph{et al.}~\cite{EsparzaKieferLuttenberger2010} are now extended to our context where the equations are not necessarily polynomial. The basic result allowing for the generalization is the following.
\begin{lemma}\label{lemma:Esparza}
Let $\bs Y=\bs H(z,\bs Y)$ be a well-founded system, $0\le a$ and $\bs 0\le \bs t\le \bs y$ be such that $(a,\bs y)$ is inside the domain of convergence of $\bs H$, then
\[\bs H(a,\bs y)-\bs H(a,\bs t)\le\frac12\left(\frac{\bpartial\bs H}{\bpartial\bs Y}(a,\bs t)
+\frac{\bpartial\bs H}{\bpartial\bs Y}(a,\bs y)\right)(\bs y-\bs t).\]
\end{lemma}
\begin{proof}Since $\bs H$ has nonnegative Taylor coefficients at~0, so do all its derivatives and therefore so does its Taylor expansion at any point with nonnegative real coordinates inside its domain of convergence.

The case of dimension~1 is obtained by Taylor expansion. We write $h(y)$ for $H(a,y)$.
By Taylor expansion,
\begin{align*}
h'(y)&=h'(t)+\sum_{k\ge2}{h^{(k)}(t)\frac{(y-t)^{k-1}}{(k-1)!}}.\\
\intertext{Multiplying by $y-t$ gives}
(y-t)h'(y)&=(y-t)h'(t)+\sum_{k\ge2}{h^{(k)}(t)k\frac{(y-t)^{k}}{k!}}.\\
\intertext{By nonnegativity of the coefficients, this is}
&\ge(y-t)h'(t)+2\sum_{k\ge2}{h^{(k)}(t)\frac{(y-t)^{k}}{k!}},\\
&=(y-t)h'(t)+2\left(h(y)-h(t)-h'(t)(y-t)\right).
\end{align*}
Reorganizing terms concludes the case of dimension~1.

In higher dimension, it is sufficient to apply this reasoning to each $\partial H_i/\partial y_j$.
\end{proof}
\noindent This lemma gives control over the speed of convergence of Newton iteration as follows.
\begin{lemma}\label{lemma:induction-Esparza}
Let $\bs Y$, $\bs y^{[n]}$ and $0\le a\le\rho$ be as in \cref{thm:Newton} and let $\bs D(a)$ be an eigenvector of $\bpartial\bs H/\bpartial\bs Y(a,\bs Y(a))$ with positive entries for the eigenvalue $\Lambda_{\bs H}(a)\le1$. If $\mu>0$ is such that $\bs Y(a)-\bs y^{[n]}\le\mu\bs D(a)$ then $\bs Y(a)-\bs y^{[n+1]}\le(\mu/2)\bs D(a)$.
\end{lemma}
\begin{proof}The existence of $\bs D(a)$ is given by the Perron-Frobenius theorem. The result is obtained by combining the inequalities as follows
\begin{align*}
\bs Y(a)-\bs y^{[n+1]}
&=\bs Y(a)-\bs y^{[n]}-\left(\Id-\frac{\bpartial\bs  H}{\bpartial \bs Y}(a,\bs y^{[n]})\right)^{-1}(\bs H(a,\bs y^{[n]})-\bs y^{[n]})\\
&=\left(\Id-\frac{\bpartial\bs  H}{\bpartial\bs  Y}(a,\bs y^{[n]})\right)^{-1}\left(\bs Y(a)-\bs H(a,\bs y^{[n]})-\frac{\bpartial\bs  H}{\bpartial\bs  Y}(a,\bs y^{[n]})(\bs Y(a)-\bs y^{[n]})\right)\\
&\le\left(\Id-\frac{\bpartial\bs  H}{\bpartial\bs  Y}(a,\bs y^{[n]})\right)^{-1}
\frac12\left(\frac{\bpartial\bs  H}{\bpartial\bs  Y}(a,\bs Y(a))-\frac{\bpartial\bs  H}{\bpartial \bs Y}(a,\bs y^{[n]})\right)(\bs Y(a)-\bs y^{[n]})\\
&\le\frac\mu2\left(\Id-\frac{\bpartial\bs  H}{\bpartial\bs  Y}(a,\bs y^{[n]})\right)^{-1}
\left(\frac{\bpartial\bs  H}{\bpartial\bs  Y}(a,\bs Y(a))-\frac{\bpartial\bs  H}{\bpartial \bs Y}(a,\bs y^{[n]})\right) \bs D(a)\\
&\le\frac\mu2\left(\Id-\frac{\bpartial\bs  H}{\bpartial \bs Y}(a,\bs y^{[n]})\right)^{-1}\left(\Id-\frac{\bpartial\bs  H}{\bpartial\bs  Y}(a,\bs y^{[n]})\right)\bs D(a)\\
&=\frac\mu2\bs D(a).
\end{align*}
In the third line, the inequality follows from \cref{lemma:Esparza} and the positivity of the inverse of the matrix, which is a consequence of~$\Lambda_{\bs H}(a)\le1$. The next line comes from the induction hypothesis. The following one uses the fact that $\bs D(a)$ is an eigenvector or $\bpartial\bs H/\bpartial\bs Y(a,\bs Y(a))$ for $\Lambda_{\bs H}(a)\le1$.
\end{proof}
We can finally conclude on the speed of convergence of Newton's iteration inside the domain of convergence of the generating function.
\begin{theorem}\label{thm:speed}
Let $\bc Y=\bc H(\cal Z,\bc Y)$ be a well-founded system and let $\rho$ be the radius of convergence of the generating function~$\bs Y$ of the solution. Let $\bs y^{[n]}$ be defined by the Newton iteration from \cref{eq:Newton-gf-val}. 

Let $0\le a\le\rho$ be such that $(a,\bs Y(a))$ belongs to the domain of convergence of $\bs H$ and let $\bs D(a)$ be an eigenvector of $\bpartial\bs H/\bpartial\bs Y(a,\bs Y(a))$ with positive entries for the eigenvalue $\Lambda_{\bs H}(a)\le1$. Then there exists $\lambda>0$ such that for all $n\ge0$,
\[\bs Y(a)-\frac{\lambda}{2^n}\bs D(a)\le\bs y^{[n]}\le \bs Y(a).\]
\end{theorem}
\begin{proof}
The second inequality comes from \cref{thm:Newton}.
Since both $\bs Y(a)$ and $\bs D(a)$ have nonnegative entries, there exists $\lambda>0$ such that $\bs Y(a)\le\lambda \bs D(a)$. This value gives the result for $n=0$ and the conclusion then follows by induction from \cref{lemma:induction-Esparza}.
\end{proof}
It follows that a linear number of iterations is sufficient to obtain a linear number of digits of the values, however close $a$ is to~$\rho$  and even for $a=\rho$ if $(\rho,\bs Y(\rho))$ belongs to the domain of convergence of~$\bs H$. A more refined analysis along the lines of that of Esparza \emph{et al.}~\cite{EsparzaKieferLuttenberger2010} could be of interest. 

\subsection{Radius of Convergence of Irreducible Systems}\label{sec:newton-radius-irred}
The characteristic system of \cref{sec:ift} has a set of solutions that contains~$(\rho,\bs Y(\rho))$, where~$\rho$ is the radius of convergence of the generating function solution. The computation of $(\rho,\bs Y(\rho))$ plays an important role in random generation by Boltzmann sampling, where it is used in the \emph{singular Boltzmann samplers}~\cite[\S~7]{DuchonFlajoletLouchardSchaeffer2004}. That system is not a positive system in the sense of \cref{def:wf-analytic} and therefore \cref{thm:Newton} cannot be applied to guarantee convergence of Newton's iteration. Still, 
denoting by~$\bs D$ the matrix $\bs D=\Id-\frac{\bpartial\bs H}{\bpartial\bs Y}$, the Jacobian matrix of the characteristic system can be observed to have the block structure 
\begin{equation}\label{eq:jac-char}
\bs J_{\text{char}}(z,\bs Y):=\begin{pmatrix}
    \bs D&-\frac{\partial\bs H}{\partial z}\\
    \frac{\bpartial\det \bs D}{\bpartial\bs Y}&\frac{\partial\det\bs D}{\partial z}
\end{pmatrix}
\end{equation}
which can be used to perform a Newton iteration in the irreducible case, thanks to the following.
\begin{theorem}\label{thm:newton-radius}
Let~$\bs Y=\bs H(z,\bs Y)$ be an irreducible system, let~$\bs Y$ be its generating function solution and $\rho$ its radius of convergence. If~$(\rho,\bs Y(\rho))$ is inside the domain of convergence of~$\bs H$, then there exists a neighborhood of~$\bs\zeta:=(\bs Y(\rho),\rho)$ such that Newton's iteration for the characteristic system, i.e.,
\[\bs \zeta^{[n+1]}=\bs \zeta^{[n]}-\bs J_{\text{char}}(\zeta^{[n]})^{-1}\cdot \begin{pmatrix}\bs\zeta^{[n]}-\bs H(\bs\zeta^{[n]})\\ \det \bs D(\bs \zeta^{[n]})\end{pmatrix},\]
converges quadratically to~$\bs\zeta$ when started in that neighborhood.
\end{theorem}
Note that if Newton's iteration for the characteristic system converges to a real positive value~$(a,\bs B)$, then computing $\Lambda_{\bs H}(a,\bs B)$ can be used to certify that~$(a,\bs B)=(\rho,\bs Y(\rho))$, by \cref{thm:Lambda}. 
\begin{proof}It is classical that it is sufficient to prove that the Jacobian matrix in \cref{eq:jac-char} is invertible at~$\bs\zeta$~\cite{Dedieu2006}. By Schur's formula, its determinant is
\begin{equation}\label{eq:detJac}
\det(\bs D)\frac{\partial\det\bs D}{\partial z}+\frac{\bpartial\det\bs D}{\bpartial \bs Y}\operatorname{adj}(\bs D)\frac{\partial \bs H}{\partial z},
\end{equation}
where $\operatorname{adj}(\bs D)=\det(\bs D)\bs D^{-1}$ is the adjugate matrix of $\bs D$.

At~$(z,\bs Y(z))$ for~$z<\rho$, $\det\bs D>0$: it is the value of the characteristic polynomial $\det(\lambda \Id-\bpartial \bs H/\bpartial\bs Y)$ at $\lambda=1>\Lambda_{\bs H}(\bpartial \bs H/\bpartial\bs Y)$, where the last inequality follows from \cref{thm:Lambda}. 

Since $\bs J:=\bs\partial\bs H/\bs\partial\bs Y$ is an irreducible matrix, for each $i,j$ there is a $q\in\mathbb N$ such that the entry $(i,j)$ of $\bs J^q$ is positive~\cite[XIII.1.2]{Gantmacher1959}. For $z<\rho$, as $\bs D^{-1}$ is the sum of nonnegative integer powers of $\bs J$, it has positive entries. Thus $\operatorname{adj}\bs D$ has positive entries. This property persists for~$z=\rho$, where~1 is the dominant eigenvalue of~$\bpartial\bs H/\bpartial\bs Y$ (see~\cite[XIII.2.2]{Gantmacher1959}). 

Finally, for $z<\rho$ and $u\in\{z,Y_1,\dots,Y_m\}$, using Jacobi's formula, one has  
\[\frac{\partial\det\bs D}{\partial u}=\operatorname{tr}\left(\operatorname{adj}(\bs D)\frac{\partial\bs D}{\partial u}\right)<0,\]
by positivity of the entries of $\operatorname{adj}\bs D$
and monotonicity of the entries of $\bpartial\bs H/\bpartial \bs Y$.

Thus both terms of \cref{eq:detJac} are nonpositive. When~$z\rightarrow\rho$, the first one becomes~0. The second one remains negative: $\partial \bs H/\partial z\ge0$ is nonzero, therefore its product with $\operatorname{adj}\bs D$ has positive entries and at least one of those of  $\bpartial\det\bs D/\bpartial \bs Y$ is nonzero by irreducibility, and then it is negative.
\end{proof}

\begin{example}Consider again Burris' example from \cref{ex:Burris}. The characteristic system in that case is 
\[\{y_1 - z(y_1^2 + y_2 + 1), y_2 - z(y_2^2 + y_1 + 1), 4y_1y_2z^2 - 2y_1z - 2y_2z - z^2 + 1\}.\]
If one takes for starting point~$z=0.05375$ and the corresponding values~$(y_1,y_2)$ given by \textsf{\ref{method:newton-method}}, the first few iterates of~$(y_1,z)$ for the Newton iteration solving the characteristic system are (the value of~$y_2$ is always equal to that of~$y_1$)
\begin{multline*}
(0.05698764391,0.05375),
(1.988562422,1.766501487),
(-2.497923913,5.321620240),\\
(-1.254302784,5.312553615),
(0.555040642,11.41502810),
(0.6123418114,-0.4121335348),\\
(0.732014967,0.423558331),
(1.000037986,0.3176177955),
(0.9985108606,0.3333093468),\\
(0.9999964163,0.3333336145),(1.000000000,0.3333333333).
\end{multline*}
The fixed point is thus reached in 11~iterations. That this solution indeed corresponds to the radius of convergence is checked by verifying that $\Lambda=1$ at this point (by \cref{thm:Lambda}). Note that some of the intermediate values computed during the iteration do not correspond to values of the generating functions (for instance, for small negative~$z$, they are negative), and that the iteration uses values of~$z$ outside the domain of convergence of the generating functions. With the starting point $z=0.09$, another solution of the system is reached, namely $y_1=y_2=1-\sqrt2,z=-(2\sqrt2+1)/7$, which has negative coordinates. Empirically, the radius of convergence is found with any initial point $(z,y(z))$ such that $z\in[0.18,1/3]$.
\end{example}

\Cref{thm:newton-radius} gives an efficient iteration, but as usual with this type of result, the choice of initial point remains an issue. For constructive systems, experiments suggest to start with a point $(z,\bs Y(z))$ with $z<\rho$ and $\bs Y(z)$ computed by Newton iteration. As~$\rho$ is unknown, this seems like a circular argument, but thanks to \cref{thm:Newton} and \cref{thm:Kantorovich}, one can use dichotomy to find such a point,  
even when the domain of application of \cref{thm:newton-radius} is quite small, as in the following.

\begin{example}The grammar 
\begin{multline*}
\{\cal C_1 = \cal Z + \cal Z \times  (\cal Z + \cal Z), \quad 
\cal C_2 = \cal Z + \cal Z \times \cal C_3, \quad 
\cal C_3 = \cal Z + \cal Z \times (\cal C_2+ \Cyc(\cal C_4))
,\\ 
\cal C_4 = \cal Z + \cal Z^2 \times\Seq  (\cal C_5), \quad 
\cal C_5 = \cal Z + \cal Z^2 \times \cal C_2 \times \Cyc(\cal C_1 + \cal C_1)
\}
\end{multline*}
translates into the system of generating function equations formed of $C_1=z+2z^2$ and the irreducible system
\begin{multline*}
\left\{C_{2}\! \left(z \right) = 
z +z C_{3}\! \left(z \right),\quad  C_{3}\! \left(z \right) = 
z +z \left(C_{2}\! \left(z \right)+\ln \! \left(\frac{1}{1-C_{4}\! \left(z \right)}\right)\right)
,\right. \\ \left.
C_{4}\! \left(z \right) = z +\frac{z^{2}}{1-C_{5}\! \left(z \right)}
,\quad  C_{5}\! \left(z \right) = 
z +z^{2} C_{2}\! \left(z \right) \ln \! \left(\frac{1}{1-2(z+2z^2)}\right)
\right\}.
\end{multline*}
The difficulty here is that the singularity of the logarithm in the equation for~$C_5$ is at $(\sqrt5-1)/4\approx0.30902$, which is quite close to the radius~$\rho\approx 0.308968$.
Our implementation first uses a dichotomy to find~$\rho$ with 2~digits of precision, leading to the starting point~$(z=0.30,\bs C(z))$, but there the iteration fails by leaving the domain of convergence of the logarithm in one iteration. Next, a new dichotomy with 4~digits of precision gives a starting point~0.3089 that fails again. Finally with 8~digits, the point~$0.30896838$ is satisfactory and 20~digits of~$\rho$ are found in 5~iterations. Correctness of the result is obtained by checking the value of $\Lambda$ at that point.
\end{example}

\begin{remark}
In practice, the computation of the derivatives of~$\det\bs D$ can be performed numerically with the adjugate~$\operatorname{adj}(\bs D)$, even close to the point where $\det \bs D=0$, through a method based on the singular value decomposition of~$\bs D$~\cite{Stewart1998}.
\end{remark}
By \cref{lemma:wf-irreducible-components}, \cref{thm:newton-radius} can also be used for irreducible systems constructed from components of larger systems. With the notation of the lemma, the iteration is performed on the system~$\bs Y_{i+1:m}=\hat{\bs H}_{i+1:m}(z,\bs Y_{i+1:m})$. At each step of the iteration, the values of the variables $\bs Y_{1:i}$ are computed at the current value of~$z$ and then used to evaluate $\hat{\bs H}$ and its derivatives.
In practice however, it is easier and empirically more efficient to perform the iteration on the whole system. The convergence is still quadratic, by the following result.
\begin{theorem}\label{thm:newton-rho-irred-extended}
Let $\bs Y=\bs H(z,\bs Y)$ be a well-founded system that splits into two sub-systems satisfying the hypotheses of \cref{lemma:wf-irreducible-components} and assume that the radius of convergence~$\rho$ of its generating function solution~$\bs Y$ is smaller than that of~$\bs Y_{1:i}$ and is such that $(\rho,\bs Y(\rho))$ is inside the domain of convergence of~$\bs H$. Then the conclusion of \cref{thm:newton-radius} holds.
\end{theorem}
\begin{proof}
The proof follows the steps of that of \cref{thm:newton-radius}. The matrix~$\bs D$ has a block structure
\[\bs D=\begin{pmatrix}\bs a&\bs 0\\ -\bs c & \bs d\end{pmatrix},\]
where $\Id-\bs d$ is an irreducible matrix corresponding to the second sub-system and $\bs a$ is invertible at~$(z,\bs Y(z))$ for $z\le\rho$. The block structure of $\bs D$ induces a block decomposition of $\operatorname{adj}\bs D$:
\[\operatorname{adj}\bs D=\begin{pmatrix}\det( \bs d)\operatorname{adj}(\bs a)&\bs 0\\
\operatorname{adj}(\bs d)\,\bs c\operatorname{adj}(\bs a)&\det(\bs a)\operatorname{adj}(\bs d)
\end{pmatrix}.\]

In the factorization $\det\bs D=\det\bs a\det\bs d$, both factors are positive: the first one by invertibility for $z\le\rho$ and the fact that $\det\bs a$ is~1 at~$z=0$, since $\bs\partial\bs H_{1:i}/\bs\partial\bs Y_{1:i}$ is nilpotent at~0; the second one for $z<\rho$ by irreducibility of $\Id-\bs d$. This shows that $\det\bs D>0$ for $z<\rho$. Moreover, as in the proof of \cref{thm:newton-radius}, irreducibility implies that $\operatorname{adj}(\bs d)$ has positive entries, even for $z=\rho$. In the case of $\operatorname{adj}(\bs a)$, only nonnegativity of the entries hold in general. As a consequence, $\operatorname{adj}\bs D\ge0$.

The derivatives also decompose by block: for $z\leq\rho$ and $u\in\{z,Y_1,\dots,Y_m\}$, 
\[\frac{\partial\det\bs D}{\partial u}=
\det(\bs d)\operatorname{tr}(\operatorname{adj}(\bs a)\frac{\partial\bs a}{\partial u })+
\det(\bs a)\operatorname{tr}(\operatorname{adj}(\bs d)\frac{\partial\bs d}{\partial u}).\]
It follows form the previous discussion on signs that both terms are $\le 0$. 

Thus again, both terms of \cref{eq:detJac} are nonpositive. When $z\rightarrow\rho$, the determinant of $\bs D$ in the first one tends to~0. In the second term, all factors have constant sign or are~0. Focusing on the last block of variables, as $\operatorname{adj}\bs d>0$ and at least one of the derivatives $\partial\bs d/\partial Y_j$ for $i+1\le j\le m$ is nonzero, it shows that this factor remains negative there, proving the invertibility of $\bs J_{\operatorname{char}}$.
\end{proof}

\subsection{Newton Iteration for Other Singularities}\label{subsec:newton-u-eq-1}
Systems of generating function equations coming from constructive systems may contain equations of the form $y(z)=1/(1-u(z))$ or $y(z)=\ln(1/(1-u(z))$, where $u(z)$ is a coordinate of a vector defined by an earlier strongly connected component~$\bs U=\bs H(z,\bs U)$. If it has been detected that~$1<u(r)\le\infty$ at its radius of convergence~$r$, then the radius of convergence of~$y(z)$ is the unique positive real solution of~$u(z)=1$. Following the same reasoning as for \cref{thm:newton-radius} shows that this can be solved by Newton iteration as well.
\begin{theorem}\label{thm:newton-u-eq-1}
Let~$\bs U=\bs H_{1:m}(z,\bs U)$ be an irreducible  system, let~$\bs U(z)$ be its generating function solution and $\rho$ its radius of convergence. Let $u$ be a coordinate of~$\bs U$ such that~$u(\rho)>1$. Then there exists~$\bs \zeta=(t,\bs B)\in(0,\rho)\times \mathbb R_{\ge0}^m$ such that $u(t)=1,\bs B=\bs U(t)$ and a neighborhood of~$\bs\zeta$ such Newton's iteration for the system~$\{\bs U=\bs H(z,\bs U),u=1\}$ converges quadratically to~$\bs\zeta$ when started in this neighborhood.
\end{theorem}
Note again that if Newton's iteration for the system converges to a real positive value~$(a,\bs B)$, then computing $\Lambda_{\bs H}(a,\bs B)$ can be used to certify that~$(a,\bs B)=\bs\zeta$, by \cref{thm:Lambda}. 
\begin{proof}The proof is a simpler version of the previous one. Without loss of generality, assume that~$u$ is the first coordinate of~$\bs U$. The Jacobian matrix of the system has the block structure
\begin{equation}\label{eq:jac-u-eq-1}
\begin{pmatrix}
    \bs D&-\frac{\partial\bs H}{\partial z}\\
    1\,0\dots0&0
\end{pmatrix},
\end{equation}
where $\bs D=\Id-\frac{\bpartial \bs H}{\bpartial\bs U}$.
By Schur's formula, the determinant of \ref{eq:jac-u-eq-1} is
\begin{equation}\label{eq:detJac-u-eq1}
(1\,0\dots 0)\operatorname{adj}(\bs D)\frac{\partial H}{\partial z}.
\end{equation}
As in the proof of \cref{thm:newton-radius},
at~$(z,\bs U(z))$ for~$z<\rho$, $\operatorname{adj}\bs D=\det(\bs D)\bs D^{-1}$ has positive entries, while $\frac{\partial H}{\partial z}$ has nonnegative entries. This proves that \cref{eq:detJac-u-eq1} is positive and therefore that the Jacobian matrix is invertible.
\end{proof}
\noindent The choice of initial point for constructive systems is dealt with as above. Again, this can be used on a system that is not irreducible but is defined by the construction of \cref{lemma:wf-irreducible-components}. 


\section{Radius of Convergence of Constructible Generating Functions using Newton Iteration}\label{sec:radius}

The previous section gives properties of Newton's iteration for general well-founded systems. These properties are now exploited in the case of constructive systems (\cref{def:constructive}), which restricts the possible analytic functions that are involved. 

\Cref{th:radius-is-computable,prop:oracle-gives-radius} show that the radius of convergence of a constructible generating function is computable with the \textsf{Constant Oracle} of \cref{def:oracle}. The algorithm underlying the proof is a dichotomy relying strongly on the oracle and therefore not very efficient. Using Newton's iterations from the previous section leads to an algorithm that still relies on the oracle as those of \cref{subsec:radius}, but it uses numerical approximations to make many of its decisions. Only equalities are delegated to the oracle. For that reason, the algorithms presented in this section are labeled as `Numerical Algorithms'.

The data-structure they use to represent a radius of convergence is made of a system of equations in the radius and possibly other variables, together with an interval for each of its variables, so that the system has only one solution when each variable lies in the corresponding interval. As these intervals isolate one of the possibly many solutions of the system, they are called \emph{isolating intervals}. Together with the algorithms of the previous section, this data-structure allows one to compute the solutions to arbitrary precision, limiting to equality questions the tasks devoted to the oracle.

Our algorithms are designed for constructive systems in normal form (see \cref{def:normal_form}).
The main algorithms use subroutines that take as input a component from the decomposition given by~\cref{eq:triangular}. 
Each of these components can be written as a system in normal form $\bc Y = \bc H(\cal Z, \bc U,\bc Y)$ where $\bc U$ is the (possibly empty) set of its predecessors in the decomposition.
The sub-system $\bc U = \bc G(\cal Z,\bc U)$ defining $\bc U$ is also a constructive system in normal form. 

For such a system, testing that a point belongs to the domain of convergence of the generating function~$\bs H$ of $\bc H$ is easy: since the system is constructive and in normal form, the domain of convergence is defined by the inequalities $|z|<r_U$ for $U\in\bs U$, plus the inequalities $|Y_i|<1$ for any~$i$ such that $\Cyc(\mathcal Y_i)$ or $\Seq(\mathcal Y_i)$ is the right-hand side of an equation of~$\bc H$. 

The main result of this section is the following.
\begin{theorem}\label{thm:algo-radius}
Given a \texttt{Constant Oracle} as in \cref{def:oracle}, Algorithm \textup{\textsf{\ref{algo:radius}}} decides whether the generating function solution has a finite radius of convergence or not; if it is finite, it returns a system of equations defining it and a tuple of isolating intervals (one for each variable of the system) where the solution is unique.
\end{theorem}
We first illustrate the working of the algorithm on a few examples, then deal separately with the basic cases of nonrecursive equations and of irreducible systems, before giving the algorithm that combines them, using the decomposition of systems into their strongly connected components.

\subsection{Examples}\label{sec:algo_big_spec}
\begin{example}\label{ex:radius_sys1}Consider the system~\eqref{eq:gf-colored-forests} satisfied by the generating functions from \cref{ex:colored_forest}. 
The algorithm works in a bottom-up way on the strongly connected components.

\subparagraph{First Component $\mathcal G$.}
The generating function of~$\cal G$ is defined by the single nonlinear recursive equation
\[
S_G=\{G(z) = z + G(z)^2\}.
\]
According to \cref{thm:Lambda}, its radius of convergence is either $+\infty$ or a solution of the characteristic system~$\Sigma=\{G = G^2 + z, 1 - 2G = 0\}$. This system has a unique solution such that~$G$ and $z$ are real positive, for which $z=1/4$. Thus, by \cref{thm:Lambda}, the radius of convergence for this component is\[\rho_G=1/4.\] Since one of the equations is linear, Newton's iteration from \cref{thm:newton-radius} converges in 2~steps to this value. Our algorithm returns the system $\Sigma$ and a pair of intervals around $(1/4,1/2)$ where it has a unique solution~$(\rho_G,G(\rho_G))$.

\subparagraph{Second Component $(\mathcal B,\mathcal R)$.}
The generating functions of $\cal B$ and $\cal R$ are defined by the irreducible system
\[
S_{BR}=\left\{B(z) = \frac{z}{1 - R(z)}, R(z) = z^3 + \frac{z}{1 - B(z)}\right\}.
\]    
Again, one can use Newton's iteration on the characteristic system. Starting with the value of $(B,R)$ at $z=0.2$, one gets $z\in [0.2462661,0.2462662]$. As the largest eigenvalue of the Jacobian matrix at~$0.2462661$ is smaller than~$1$, by \cref{thm:Lambda} again, the radius of convergence~$\rho_{BR}$ satisfies
\[\rho_{BR}\in [0.2462661,0.2462662].\]
Besides the characteristic system, our algorithm returns this interval and the corresponding intervals
$(B(\rho_{BR}),R(\rho_{BR}))\in(0.4919992,\infty)\times(0.4995864,\infty)$ where the characteristic system has a unique solution.

\subparagraph{Third Component $(\cal T_r,\cal T_b,\cal T_g)$.}
The corresponding system of generating function equations is the union of the previous two and
\begin{equation}\label{eq:TrTbTg}
S=\left\{T_r(z)=z + \frac{R(z)}{1 - T_b(z)}, T_b(z)=\frac{B(z)}{1 - T_g(z)}, T_g(z)=\frac{G(z)}{1 - T_r(z)}\right\}.
\end{equation}
Again, a consequence of \cref{thm:Lambda} is that the radius of convergence~$\rho_T$ of this system is either on the boundary of the domain of convergence of the system, that is the minimum of~$\rho_G$ and~$\rho_{RB}$, or the equation~$\Lambda(\rho_T)=1$ is satisfied. 
The numerical Newton iteration for this system shows that $\rho<0.2462661$. Thus there exists a solution with real positive values and $z<0.2462661$ for the corresponding characteristic system $S\cup \{D_S=0\}$, with
$D_S$ the determinant of the Jacobian matrix of the system coming from $S$, 
\begin{equation}\label{eq:DS}
D_S(z)=\frac{B(z)R(z)G(z)}{(1-T_r(z))^2(1-T_b(z))^2(1-T_g(z))^2}-1.
\end{equation}
Starting with the value of $(G,B,R,T_r,T_b,T_g)$ at~$z=.16$ and using the iteration of \cref{thm:newton-rho-irred-extended} gives
\[\rho_{T}\in[0.1703916,0.1703917],\]
{together with}
\[(G,B,R,T_r,T_b,T_g)|_{z=\rho_T}>(0.2178503, 0.2193556, 0.2232174, 0.5883123, 0.4658856, 0.5291642),\]
which gives a domain where the characteristic system has a unique solution. In this purely algebraic setting, one can also use elimination to produce an irreducible polynomial of degree~21 that has~$\rho_T$ as its unique root in the interval above, namely
\begin{multline}\label{pol-deg-21}
z^{21}+2 z^{20}+11 z^{19}+10 z^{18}+39 z^{17}-24 z^{16}+44 z^{15}-216 z^{14}+95 z^{13}-412 z^{12}+465 z^{11}-438 z^{10}\\
+678 z^{9}-600 z^{8}+897 z^{7}+192 z^{6}+1332 z^{5}+176 z^{4}+268 z^{3}-272 z^{2}+60 z -4.
\end{multline} 


\subparagraph{Last Component $\cal F$.} Since the generating function of $\cal F$ is $F(z)=\exp(T_r(z))$, the radius of convergence is $\rho_T\in[0.1703916,0.1703917]$ and it arises from the irreducible component defining $T_r$. Thus it is returned with the system $S\cup\{D_S=0\}$ from that component.
\end{example}

\begin{example}We change the system for the component $(\cal T_r,\cal T_r,\cal T_g)$ in the previous example by adding cardinality constraints on the nodes of the trees, so that the system for its generating functions becomes 
\[
S_2=\left\{T_r(z) = z + R(z)T_b(z)^2, T_b(z) = B(z)(1 + T_g(z) + T_g(z)^2), T_g(z) = G(z)(1 + T_r(z))\right\}.
\]
The first two components $S_G$ and $S_{RB}$ are identical. The computations on the third component are done in the same way as above and yield that $\rho_T\in[0.246067,0.246068]$ and, as before, it is the radius of convergence of the whole system, since $F(z)=\exp(T_r(z))$. But let us change again the system, to consider ordered trees, setting $\cal F = \Seq(\redroot)$. We now have to look for a singularity coming from~$\cal F$ whose generating function is $F(z)=1/(1-T_r(z))$, since, in this case, $T_r(0.246067)\sim 1.39031>1$.
Using the numerical Newton iteration to compute a value of~$z$ such that $T_r(z)=1$ following \cref{thm:newton-u-eq-1} and starting from the value of $(G,B,R,T_r,T_b,T_g)$ at $z=0.21$, we find that the radius of convergence is the unique solution of the system $S_G\cup S_{BR}\cup S_2\cup\{T_r(z)=1\}$ with
\[z\in[0.2457166,0.2457167]\] 
and\\[4pt]
\scalebox{.96}{$(0.434, 0.475, 0.483, 0.999, 1.24, 0.869)<(G,B,R,T_r,T_b,T_g)<(0.435, 0.476, 0.484, 1.01, 1.25, 0.870)$.}
\end{example}

\begin{example}We change again the system  for the component $(\cal T_r,\cal T_r,\cal T_g)$ so that the generating function equations become
\[
\left\{T_r(z) = z + R(z)
T_b(z)^2, T_b(z) = B(z)T_g(z)^2, T_g(z) = G(z)T_r(z)^2\right\},
\]
The Jacobian matrix for this system is
\[
\bs J=\begin{pmatrix}
0&2\,B T_g&0\\ 
0&0&2\,G T_r\\ 2\,R T_b&0&0
\end{pmatrix}.
\]
In this case, the numerical Newton iteration for~$(\cal T_r,\cal T_r,\cal T_g)$ converges at~$z=0.24626615$, where the largest eigenvalue of $\bs J$ is $\approx 0.127<1$, which implies that~$\rho=\rho_{BR}\in[0.2462661,0.2462662]$, from the component~$(\cal B,\cal R)$.
To conclude this last example, recall that $F(z)=\exp(T_r(z))$, thus the radius of convergence~$\rho$ of the whole system is $\rho_{BR}$ itself and it arises from the irreducible component defining~$(\cal B,\cal R)$. 
\end{example}

\subsection{Nonrecursive Equations}\label{sec:nonrec}
\begin{numalgo}[t]
\SetAlgoRefName{NR-Radius}
\caption{Radius of convergence \label{algo:nr-radius}}
\DontPrintSemicolon
\Input{$\cal{Y}=\cal{H}(\cal{Z},\bc{U})$ a constructive component in normal form consisting of one nonrecursive equation;\\
\phantom{{\bf Input:}} $\bc U = \bc G(\cal Z,\bc U)$ a constructive system in normal form defining $\bc U$;\\
\phantom{{\bf Input:}} $r_{\bs U}$ radius of convergence of the generating function $\bs{U}$ of $\bc{U}$.
}
\Output{$\min(\rho,r_{\bs U})$, given by a system and isolating intervals, or $\infty$,\\
\phantom{{\bf Output:}} with $\rho$ the radius of convergence of the generating function ${Y}$ of $\cal{Y}$.
\vspace{10pt}}
    \lIf{$\cal{H}(\cal{Z},\bc{U})\in\{1,\cal{Z}\}$}{\Return{$\infty$}}
    \lElseIf{$\cal{H}(\cal{Z},\bc{U})=\Set(\cal{U})$ with $\cal U\in\bc{U}$}{\Return{$r_{\bs U}$}}
    \ElseIf{$\cal{H}(\cal{Z},\bc{U})\in\{\Seq(\cal{U}),\Cyc(\cal{U})\}$ with $\cal U\in\bc{U}$}{
    \tcc*{use \textsf{NewtonIteration} for $U$ in $[0,r_{\bs U}]$}
        \lIf(\tcp*[f]{decide inequality}){$U(r_{\bs U})\le1$}{
            \Return{$r_{\bs U}$}}\label{line7}
        \Else{
            Find an interval~$I_\rho$ for~$\rho\in(0,r_{\bs U})$ s.t. $U(\rho)=1$\tcp*[r]{
            use \cref{thm:newton-u-eq-1}}\label{line9}
            \Return{$I_\rho\cup\bs U(I_\rho),\{\bs U=\bs G(\rho,\bs U)\}\cup\{U=1\}$}\tcp*[r]{
            for $I_\rho$ and $\bs U(I_\rho)$}}}\label{line8}
    \lElse(\tcp*[f]{$\cal{H}$ is a cartesian product or a union}){\Return{$r_{\bs U}$}}
\end{numalgo}
The simplest situation for the computation of the radius of convergence is the case of a nonrecursive equation. Then Algorithm \textsf{\ref{algo:nr-radius}} is a simple distinction of cases, based on properties of the functions $\exp$, $x\mapsto1/(1-x)$ and $x\mapsto\ln1/(1-x)$. 
This algorithm is a refinement of `Algorithm Radius' in the old work of Flajolet, Salvy, Zimmermann~\cite[p.~77]{FlajoletSalvyZimmermann1991}. The main difference is that now, the equation arises inside a larger system, thus information concerning the generating functions on which it depends (denoted $\bc U$ in the algorithm) and that have been analyzed before, forms part of the input of the algorithm. 

Deciding $U(r_{\bs U})<1$ at line~\ref{line7} is obtained by comparing~1 with the solution of a system 
defining~$U(r_{\bs U})$. Such a system is obtained by starting from the system~$S_{r_{\bs U}}$ defining~$r_{\bs U}$ and adding the equations~$\bs U(r_{\bs U})=\bs G(r_{\bs U},\bs U(r_{\bs U}))$ for the coordinates of $\bs U$ that are not in $S_{r_{\bs U}}$ already. The uniqueness of the solution of $S_{r_{\bs U}}$ in the given intervals implies the uniqueness of the solution of the extended system, whose coordinates contain $r_{\bs U}$ and the values of the $\bs U(r_{\bs U})$, and thus in particular $U(r_{\bs U})$. Comparing the latter value to~1 is achieved by checking whether equality holds, using the \texttt{Constant Oracle}, and refining the interval numerically if they are distinct.

In the pseudo-code, the function \textsf{NewtonIteration} takes as input a combinatorial system and a value of~$x$. It implements the iteration of \cref{thm:Newton}. At each iteration, it checks that: the coordinates are increasing; their values belong to the domain of convergence of the corresponding~$\bs H$; the corresponding value of $\Lambda_{\bs H}(a,\bs y^{[n]})$ is at most~1. It exits and returns \texttt{fail} if one of these conditions is not satisfied.

\subsection{Irreducible Components}\label{sec:irred-comp}

\begin{numalgo}[t]
\SetAlgoRefName{Irr-Radius}
\caption{Radius of convergence \label{algo:irr-radius}}
\DontPrintSemicolon
\Input{$\bc{Y}=\bc{H}(\cal{Z},\bc{U},\bc{Y})$ an irreducible constructive component in normal form;\\
\phantom{{\bf Input:}} $\bc U = \bc G(\cal Z,\bc U)$ a constructive system in normal form defining $\bc U$;\\
\phantom{{\bf Input:}} $r_{\bs U}$ radius of convergence of the generating function $\bs{U}$ of $\bc U$.
}
\Output{$\min(\rho,r_{\bs U})$, given by a system and isolating intervals, or $\infty$,\\
\phantom{{\bf Output:}} with $\rho$ the radius of convergence of the generating function $\bs Y$ of $\bc Y$.
\vspace{10pt}}
    Set $\Sigma:=\{\bs U=\bs G(z,\bs U),\det(\Id-\bpartial\bs H/\bpartial\bs Y )=0\}$\tcp*[r]{Characteristic system}\label{line:irr-radius2bis}
    
    \lIf{$\bs H$ is not linear}{Set $\Sigma:=\Sigma \cup \{\bs Y=\bs H(z,\bs U,\bs Y)\}$}\label{line:irr-radius2ter}
    \smallskip
    \tcc*[r]{recall that $\Lambda_{\bs H}(x,\bs Y)$ is the largest eigenvalue of $\bpartial\bs H/\bpartial\bs Y(x,\bs{U},\bs{Y})$}
    \If(\tcp*[f]{Oracle}){\label{algo:irr-rad:line4}$\Sigma$ has a solution $(a,\bs U,\bs B)\ge\bs 0$ 
    where $\bs H$ is convergent\\
     \qquad and $\Lambda_{\bs H}(a,\bs B)=1$ (necessarily $a\le r_{\bs U}$)\label{algo:irr-rad:line3}}{
     {
     Find an isolating interval~$I_a$ for this~$a$\tcp*[r]{use \cref{thm:Newton,thm:newton-rho-irred-extended}}\label{line:irr-radius5}
     \Return $I_a\cup \bs U(I_a)\cup\bs B(I_a),\Sigma$\tcp*[r]{for $I_a,\bs U(I_a),\bs B(I_a)$}}
     }
    \lElse{\Return $r_{\bs U}$}\label{line:irr-radius4}
\end{numalgo}
Let $\bc Y=\bc H(\cal Z,\bc U,\bc Y)$ be an irreducible component input. Algorithm \textsf{\ref{algo:irr-radius}} puts together some of the results of the previous sections.

If the radius of convergence~$\rho$ of the solution is smaller than the radius~$r_{\bs U}$ of the generating function~$\bs U$, then by  \cref{coro:lambda=1,coro:lambda-irr-comp} either $\bs Y$ is infinite at~$\rho$ (then the system is linear, by \cref{lemma:finiteornot}) or it is finite and~$(\rho,\bs U,\bs Y)$ is a solution of the characteristic system with $\Lambda_H(\rho,\bs U,\bs Y)=1$.

If the system is linear, it is of the form $\bc Y=\bc C(\bc U)\bc Y+\bc D(\bc U)$. The generating function solution is $(\Id-\bs C(\bs U))^{-1}\bs D(\bs U)$. It is singular if and only if either $\bs U$ is singular or $\det(\Id-\bs C(\bs U))=0$ at the solution of $\bs U=\bs G(z,\bs U)$. The latter case is precisely the system $\Sigma$ computed in \cref{line:irr-radius2bis}. The hypothesis that $\bs H$ is convergent at $(a,\bs U)$ amounts to $\bs C$ and $\bs D$ being convergent there, which implies that $\rho<r_{\bs U}$. The extra condition $\Lambda_{\bs H}(a,\bs U)=1$ ensures that $a=\rho$ by \cref{coro:lambda-irr-comp}.

If the system is not linear, 
the system~$\Sigma$ computed in \cref{line:irr-radius2bis,line:irr-radius2ter} is the characteristic system of \cref{eq:characteristic-system}. If $\rho<r_{\bs U}$, it has a solution $(a,\bs U,\bs B)\ge0$ where $\bs H$ is convergent  and $\Lambda_H(a,\bs U,\bs B)=1$. Then \cref{coro:lambda-irr-comp} shows that $a$ is the radius of convergence for this component.

\subsection{General Case}

\begin{numalgo}[ht]
\SetAlgoRefName{Radius}
\caption{Radius of convergence (general case)\label{algo:radius}}
\DontPrintSemicolon
\Input{$\bc{Y}=\bc{H}(\cal{Z},\bc{Y})$ a constructive system in normal form with $m$ equations.}
\Output{Radius of convergence of the generating function $\bs{Y}$ of $\bc{Y}$.\vspace{10pt}}
    Decompose the system into its $k$ strictly connected components (note that $i_k=m$)
    ${}\qquad\bc{Y}_{1:i_1}=\bc{H}_{1:i_1}(\cal{Z},\bc{Y}_{1:i_1}),
\dots,\bc{Y}_{i_{k-1}+1:i_k}=\bc{H}_{i_{k-1}+1:i_k}(\cal{Z},\bc{Y}_{1:i_k}).$

    $r:=\infty$
    
    \For(\tcp*[f]{$i_0=0$}){$j=1$ to $k$}{
        ${\bs S}_{\bc U}:=\bc{U}=\bc{H}_{1:i_{j-1}}(\cal Z,\bc{U})$\\

        \If(\tcp*[f]{at most one nonrecursive eqn}){$\partial\bs{H}_{i_{j-1}+1:i_j}/{\partial\bs{Y}_{i_{j-1}+1:i_j}}=0$}
        {$r:=\textsf{NR-Radius}(\mathcal{Y}=\cal{H}_{i_j}(\mathcal{Z},\bc{U},\cal{Y}),\bs S_{\bc U},r)$}
        \lElse{
        $r:=\textsf{Irr-Radius}(\bc{Y}=\bc{H}_{i_{j-1}+1:i_j}(\cal{Z},\bc{U},\bc{Y}),\bs S_{\bc U},r)$
        }
    }
    \Return{$r$}
\end{numalgo}

Finally, the general case is described in Algorithm \textsf{\ref{algo:radius}}. It first decomposes the system into its strongly connected components and then deals with each of them in order, using either of the previous algorithms, depending on the nature of the component. 

The isolating intervals that are returned can be refined to any precision $\epsilon>0$ using standard numerical solvers.
Note that the cases with infinite radius of convergence are straightforward by \cref{prop:finiteness}. Ultimately, they are detected by Algorithm \textsf{\ref{algo:nr-radius}}.


\section{Dominant Singularities of Irreducible Systems}\label{sec:dom-irred} 

The algorithms of the previous section give access to the radius of convergence of the generating function solutions to a combinatorial system.

If $F(z)=\sum_{n\ge0}{f_nz^n}$ is a generating function, the knowledge of its radius of convergence~$\rho$ gives the exponential order of the growth of the sequence of coefficients in the sense that
\[\limsup_{n\ge0}{|f_n|^{1/n}}=\rho^{-1}.\]
Since the coefficients $f_n$ are nonnegative, by Pringsheim's theorem~\cite{Titchmarsh1939} (see also~\cite{Vivanti1893,Hadamard1954}), $\rho$ itself is a singularity of~$F$. When it is the only singularity on the circle of convergence, then the estimate above can be greatly refined into an asymptotic expansion by a further analysis of $F(z)$ as $z\rightarrow\rho$, as recalled in \cref{sec:transfer}. An oscillatory asymptotic behavior occurs when several singularities lie on the circle of convergence. The singularities whose modulus is the radius of convergence are called \emph{dominant singularities}. When there are finitely many of them, as will be the case in our combinatorial applications, again a local analysis can be performed to capture this oscillatory behavior. We now turn to the determination of these dominant singularities and their arguments.

Sequences of combinatorial origin often have dominant singularities with strongly restricted arguments. For instance, Berstel showed than an ${\mathbb R}_+$-rational function (i.e., one obtained from positive real numbers and a variable $z$ by additions, multiplications and applications of $a\mapsto(1-za)^{-1}$) has all its dominant singularities with arguments that are rational multiples of~$\pi$~\cite{Berstel1971,Soittola1976}. 

Note that positivity by itself is not sufficient for this property to hold. A typical example is the following, related to Chebyshev polynomials~\cite{Berstel1971}.

\begin{example}\label{ex:badsing}For all integers $a$ and $c$ such that $0<a<c$, the rational function 
\[f(z)=\frac{1-az}{1-2az+c^2z^2}+\frac{1}{1-cz}=
2+(a+c)z+2a^2z^2+(a+c)(2a-c)^2z^3
+\dotsb\]
has positive integer coefficients. 
It has three dominant singularities, at $1/c$ and $\exp(\pm i\theta)/c$, with $\theta=\arccos(a/c)$ that is not a rational multiple of~$\pi$ in general. 
\end{example}

This property~---~arguments of dominant singularities are rational multiples of~$\pi$~---~holds for solutions of well-founded systems comprising only one equation, either when they are not recursive and use only the basic species of \cref{tab:sum_esp_sg}~\cite[Prop.~4.4]{FlajoletSalvyZimmermann1991}, or when they are recursive under an extra assumption~\cite[Lemma~26]{BellBurrisYeats2006}. Moreover, the generating function solution in those cases can be written in the form $f(z)=z^ag(z^q)$ and the dominant singularities have complex arguments that are integer multiples of $2\pi/q$; the integers $a,q$ can be computed directly from the combinatorial definition.

For irreducible systems with more than one equation, that same decomposition property actually  holds in the nonlinear case as well. This has been proved by Bell, Burris and Yeats in their study of supports of power series~\cite[Thm.~4.6]{BellBurrisYeats2011}.
\begin{definition}\label{def:support_valuation}
The \emph{support} of a nonzero power series~$f\in R[[z]]$ is the set~$\operatorname{Supp}(f)$ of indices of its nonzero coefficients. The \emph{valuation} $\val(f)$ of~$f$ is $\min\operatorname{Supp}(f)$. Its \emph{period} is $\gcd(\operatorname{Supp}(f)-\val(f))$. By extension the period of a power series reduced to one monomial is~0.
\end{definition}
Regarding generating functions, the definition of the valuation corresponds to the one given for species in \cref{sec:characterization-leading}.
\begin{theorem}\cite[Thm.4.6, Cor.4.7]{BellBurrisYeats2011}\label{thm:periods}
Let $\bs Y=\bs H(z,\bs Y)$ be an irreducible system and let $\bs Y$ be its generating function solution.
All coordinates $Y_i$ have infinite support and are periodic with the same period~$q$. If moreover $\bs H$ is not linear, then they are of the form $Y_i(z)=z^{\val(Y_i)}V_i(z^q)$ with $V_i$ a power series with at most finitely many coefficients equal to~0.
\end{theorem}

\begin{proof}Bell, Burris and Yeats prove this result when $\bpartial\bs H/\bpartial\bs Y(0,\bs 0)=\bs 0$ and $\bs H(0,\bs 0)=\bs 0$. By the reductions of \cref{sec:comparison}, all well-founded systems can be transformed into systems with those properties, but these reductions do not preserve valuations and periods. Still, the arguments used in their proofs apply. We sketch the proof here and refer to their article~\cite{BellBurrisYeats2011} for details. 

For all $i\in\{1,\dots,n\}$, let $v_i$ be the valuation of~$Y_i$ so that $Y_i=z^{v_i}\hat{Y}_i(z)$ with $\hat{Y}_i(0)\neq0$. Note that a Taylor expansion of~$Y_i=H_i(z,\bs Y)$ at~$(z,\bs 0)$ shows that $\val(\partial H_i/\partial Y_j)(z,\bs 0)\ge v_i-v_j$ for all $(i,j)$. The system 
\[\hat Y_i=H_i(z,z^{v_1}\hat{Y}_1,\dots,z^{v_n}\hat{Y}_n)/z^{v_i}=:\hat{H}_i(z,\hat{\bs Y}),\qquad i=1,\dots,n\] is well founded: that $\hat{H}_i$ is a power series follows from the valuations above; the coefficients of~$\hat{\bs H}$ are nonnegative since they are formed by positive combinations of coefficients of~$\bs H$; the~$\bs U_k$ are well defined because~$\bs H$ itself is well founded; $\bpartial\hat{\bs H}/\bpartial\hat{\bs Y}(0,\hat{\bs Y}(\bs 0))=\bpartial\bs H/\bpartial\bs Y(0,\bs Y(\bs 0))$ is nilpotent. 

By \cref{lemma:NoZeroInJac}, we can assume that $\partial\hat H_i/\partial\hat Y_j\neq0$ for all~$(i,j)$. Then, the Taylor expansion of~$\hat H_i$ contains a nonzero monomial $cz^p\hat{\bs Y}{}^{\bs k}$ with $k_j>0$. Since~$\hat{Y}_\ell(0)\neq0$ for all~$\ell$, $\hat Y_i=\hat H_i(z,\hat{\bs Y})$ implies that $p\in\operatorname{Supp}\hat Y_i$ and $p+\operatorname{Supp}(\hat Y_j)\subset\operatorname{Supp}(\hat Y_i)$, so that the period of~$\hat Y_i$, which is that of $Y_i$, divides that of~$\hat Y_j$, which is that of $Y_j$. Since this holds for all~$(i,j)$, they are equal. Applying this argument with~$i=j$, one can further assume that there is an index~$j$ such that~$p>0$, since otherwise the Jacobian matrix would not be nilpotent. For such a~$j$, it follows that $p+\operatorname{Supp}(\hat Y_j)\subset\operatorname{Supp}(\hat Y_j)$, which shows that this support is infinite and then so are those of the other~$\hat Y_i$ by the inclusion above.

Let~$q$ be the common period of the~$\hat{Y}_i$. When~$\bs H$ is not linear, a similar reasoning shows that~$p+\operatorname{Supp}(\hat Y_j^k)\subset\operatorname{Supp}(\hat Y_j)$ for some $k>1$ and $p>0$. The conclusion then follows from~\cite[Lemma~2.13]{BellBurrisYeats2011}. 
\end{proof}

As shown by \cref{ex:badsing}, the theorem above is not sufficient to conclude that the dominant singularities have arguments that are rational multiples of~$\pi$. A result of that type is announced for nonlinear irreducible systems by Bell, Burris, Yeats~\cite[Prop. 7.2(i)]{BellBurrisYeats2011}, with a one-line proof from which we could not reconstruct their intent and that seems to be missing a hypothesis. A proof based on Wielandt's theorem is given below.

\begin{theorem}\label{thm:BBY72i} In the conditions of \cref{thm:periods}, let $q$ be the period of the coordinates of the generating function solution~$\bs Y$. If the system is \emph{nonlinear}, then the dominant singularities are of the form~$\rho\exp(2i\pi k/q)$, $k\in\{0,\dots,q-1\}$.
\end{theorem}
\begin{theorem}\label{thm:period-linear}
Let $\bs Y=\bs A(z)+\bs B(z)\bs Y$ be an irreducible \emph{linear} well-founded system. Let~$\rho_A$, $\rho_B$ and $\rho_Y$ be the radii of convergence of~$\bs A,\bs B,\bs Y$ and let $\Sigma_A$, $\Sigma_B$ and $\Sigma_Y$ be their sets of dominant singularities.

If the largest real eigenvalue of $\bs B(\rho_B)$ is smaller than~1, then let $\rho_C:=+\infty$ and $\Sigma_C:=\emptyset$. 
Otherwise,
there exists a unique $\rho_C\in(0,\rho_B]$ such that the largest real eigenvalue of~$\bs B(\rho_C)$ is~1. Let~$q_B$ be the gcd of the periods of the entries of~$\bs B$ if $q_B\neq0$ and the gcd of the periods of the entries of the diagonal of~$(\Id-\bs B)^{-1}$ otherwise and let
\[\Sigma_C:=\left\{\rho_C\exp(2i\pi k/q_B)\mid k\in\{0,\dots,q_B-1\}\right\}.\]
Then, $\rho_Y=\min(\rho_A,\rho_B,\rho_C)$ and
\[
\Sigma_Y\subset\bigcup_{\substack{X\in\{A,B,C\}\\ \rho_{\!X}=\rho_{Y}}}\Sigma_X.
\]
\end{theorem}

\noindent The rest of this section is devoted to a proof of these theorems.

\subsection{Periods and Dominant Singularities}
Nonnegativity of the coefficients of combinatorial generating functions has a strong impact on their growth. The first one is a consequence of the triangular inequality.
\begin{lemma}\cite[``Daffodil Lemma'']{FlajoletSedgewick2009}\label{lemma:positivity} Let $F(z)=\sum_{n\ge0}{f_nz^n}$ be a power series with nonnegative coefficients $f_n$, two of which at least being nonzero.
Let~$u>0$ be such that $F$ converges at~$z=u$. If $\phi\in[0,2\pi)$ is such that $|F(u e^{i\phi})|=F(u)$, then there exist relatively prime integers $(p,q)$ with~$q\neq0$ such that $\phi=\pi p/q$ and $F$ has period~$q$.
\end{lemma}
The next statement is purely at the level of formal power series. See \cite{BellBurrisYeats2011} for similar results.
\begin{lemma}\label{lemma:compose-period}
Let $\phi_1(z),\dots,\phi_k(z)$ (resp. $F(z,y_1,\dots,y_k)$) be nonzero power series in one (resp. $k+1$) variables with nonnegative coefficients and satisfy $\partial F/\partial y_i\neq0$ for $i=1,\dots,k$. If $F(z,\phi_1(z),\dots,\phi_k(z))$ has period $q$, then the periods of~$\phi_1,\dots,\phi_k$ are all multiples of~$q$. 
\end{lemma}
\begin{proof} Let $m_1,\dots,m_k$ be the valuations of the power series $\phi_i$ and $q_1,\dots,q_k$ their periods.
Three basic cases are multiplication by a power of~$z$, which follows from the definition of period, sum and product. Without loss of generality assume $m_1\le m_2$. Then $m_1=\val(\phi_1+\phi_2)$ and 
\[\operatorname{Supp}(\phi_1(z)+\phi_2(z))-m_1=
\left(\operatorname{Supp}(\phi_1(z))-m_1\right)\cup
\left((\operatorname{Supp}(\phi_2(z))-m_2)+(m_2-m_1)\right).\]
If $\phi_1+\phi_2$ has period~$q$, the gcd of the elements of this set is~$q$. In particular, $q_1$, as gcd of the elements of the first set in the right-hand side, is a multiple of~$q$. Then so is 
$m_2-m_1$ since $m_2\in\operatorname{Supp}(\phi_2)$ and finally so is $q_2$.

Similarly, the product $\phi_1\phi_2$ has valuation $m_1+m_2$ and
\[\operatorname{Supp}(\phi_1(z)\phi_2(z))-m_1-m_2=
\left(\operatorname{Supp}(\phi_1(z))-m_1\right)+
\left(\operatorname{Supp}(\phi_2(z))-m_2\right).\]
Again, since both sets contain~0, if this is a subset of~$q\mathbb N$, then each of them must be, implying that the gcds of their elements are both multiples of~$q$.

By induction on the number of variables, the lemma then holds for each monomial in the variables~$y_i$ and then for any sum of these monomials and thus for~$F$. The property applies to all~$\phi_i$ because $\partial F/\partial y_i\neq0$.
\end{proof}
\subsection{Dominant Eigenvalues}
The study of dominant singularities that are not real is in turn related to eigenvalues of the Jacobian matrix outside the positive real axis. The notation~$\Lambda_{\bs H}(z,\bs y)$ is now extended to denote the maximal modulus of an eigenvalue of~$\bpartial\bs H/\bpartial\bs Y$ at $(z,\bs y)\in\mathbb C^{m+1}$. A basic tool is Wielandt's theorem~\cite[\S XIII.2.3, Lemma~2]{Gantmacher1959}.
\begin{lemma}[Wielandt's Theorem]\label{lemma:Wielandt} If $A=(a_{ij})$ is a nonnegative, irreducible square matrix and $C=(c_{ij})$ is a matrix with complex entries of the same dimension such that $|c_{i,j}|\le a_{i,j}$ for all $i,j$, then any eigenvalue~$\lambda$ of $C$ has modulus bounded by that of the dominant eigenvalue~$\Lambda$ of~$A$. If for some $\lambda$ one has $|\lambda|=\Lambda$, then $C$ can be written~$e^{i\phi}DAD^{-1}$ with $\phi=\arg\lambda$ and $D$ a diagonal matrix with entries of modulus~1.
\end{lemma}
This leads to the following result for irreducible  systems.
\begin{lemma}\label{lemma:entriesJ}
Let $\bs Y=\bs H(z,\bs Y)$ be an irreducible system. Let~$\rho$ be the radius of convergence of the generating function solution~$\bs Y(z)$ and let $\bs J(z):=\bpartial\bs H/\bpartial\bs Y(z,\bs Y(z))$ be the Jacobian matrix of $\bs H$ evaluated at $(z,\bs Y(z))$. If  $\alpha$ is a singularity of $\bs Y(z)$ with $|\alpha|=\rho$,                                                                              then $\Lambda_{\bs H}(\alpha)\le\Lambda_{\bs H}(\rho)\le 1$. If $\Lambda_{\bs H}(\alpha)=1$, each entry of $\bs J(\alpha)$ has modulus equal to that of the corresponding entry in~$\bs J(\rho)$. 
\end{lemma}
\begin{proof}
First, for any nonnegative real $a\le\rho$ and real $\theta$, if $\bs H$ converges at $(a,\bs Y(a))$, then it also converges at $(ae^{i\theta},\bs Y(ae^{i\theta}))$ by absolute convergence. 

The Jacobian matrix $\bs J(z)$ also has a Taylor expansion at~0 with nonnegative coefficients. Thus, again by absolute convergence, for positive real~$a<\rho$ and real $\theta$, all the entries of $\bs J(ae^{i\theta})$ have modulus bounded by the corresponding entry in $\bs J(a)$. This implies that the spectral radius of $\bs J(ae^{i\theta})$ is at most that of~$\bs J(a)$ by Wielandt's theorem. By continuity of~$\bs J$ and~$\Lambda_{\bs H}$, this holds when~$a=\rho$ and $ae^{i\theta}=\alpha$: $\Lambda_{\bs H}(\alpha)\le\Lambda_{\bs H}(\rho)\le 1$, where the last inequality is provided by \cref{coro:lambda_at_most_1}.

If $\Lambda_{\bs H}(\alpha)=1$, then the previous inequality implies that $\Lambda_{\bs H}(\alpha)=\Lambda_{\bs H}(\rho)=1$ and each entry of $\bs J(\alpha)$ has modulus \emph{equal} to the corresponding entry in~$\bs J(\rho)$, again by Wielandt's theorem.
\end{proof}

\subsection{Proof of Theorem \ref{thm:BBY72i}} 

By \cref{lemma:finiteornot}, all coordinates~$Y_i(\rho)$ are finite and then so are the~$|Y_i(\alpha)|\le Y_i(\rho)$. 

If for all~$i$ these inequalities are equalities, then the ``daffodil lemma''~\ref{lemma:positivity} concludes.
Otherwise, one of these inequalities is strict and then by absolute convergence, $(\alpha,\bs Y(\alpha))$ lies inside the domain of convergence of~$\bs H$. Then $\alpha$ being a singularity implies that $\Lambda_{\bs H}(\alpha)=1$ by the implicit function theorem. \Cref{lemma:entriesJ} then gives that each entry of~$\bs J(\alpha)$ has modulus equal to that of the corresponding entry in~$\bs J(\rho)$. 
Since $\bs H$ is zero-free (by \cref{def-irred}) and not linear, one of the entries~$J_{ij}$ of $\bs J$ depends on a coordinate~$Y_k$ that has infinite support by \cref{thm:periods}. Then \cref{lemma:positivity} applies and $\arg\alpha=\pi p/q'$ for $q'$ a period of that entry.
Next, \cref{lemma:compose-period} implies that that coordinate $Y_k(z)$  itself has period a multiple~$q$ of~$q'$ and then so do all the coordinates~$Y_j$ by \cref{thm:periods}. 

\subsection{Proof of Theorem \ref{thm:period-linear}}
By the last point of the \cref{def:wf-analytic} of well-founded systems, the matrix $\Id-\bs B(z)$ is invertible in a neighborhood of $z=0$. It follows that 
the generating function solution is
\[(\Id-\bs B(z))^{-1}\bs A(z),\]
and it is sufficient to consider the singularities of both terms in this product and keep those that have minimum modulus. The second term is~$\bs A$, with singularities~$\Sigma_A$. 
If the largest eigenvalue of $\bs B(\rho_B)$ is smaller than~1, then $\Id-\bs B(z)$ is invertible on $[0,\rho_B)$ and its dominant singularities lie in~$\Sigma_B$.
Otherwise, since the system is well founded, the largest eigenvalue of~$\bs B(z)$ is~0 at~$z=0$ and larger or equal to~1 at~$z=\rho_B$, thus by monotonicity and continuity there is a unique~$\rho_C$ where the largest eigenvalue is~1.

If $\alpha$ is a dominant singularity of $(\Id-\bs B(z))^{-1}$ that is not a singularity of~$\bs B$, then $\bs B(\alpha)$ has~1 for dominant eigenvalue.
Then, \cref{lemma:entriesJ} applies to the system $\bs Y=\bs K(z,\bs Y)$ where $\bs K:=\Id+\bs B\bs Y$ in a situation where $\Lambda_{\bs K}(\alpha)=1$. Thus $\Lambda_{\bs K}(|\alpha|)=1$ and $|\alpha| = \rho_C$.
This shows that the entries of $\bs B(\alpha)$ have modulus equal to the corresponding entries in~$\bs B(\rho_C)$. 
Let~$q$ be the gcd of the periods of the entries of~$\bs B$. If $q\neq0$, then $q_B=q $ and one of the entries of~$\bs B$ has a nonzero period~$q$ and \cref{lemma:positivity} then implies that $\arg\alpha=\pi p/q_B$ for some integer~$p$.

If $q=0$ then by definition of the periods, the entry $(i,j)$ of~$\bs B$ is a monomial $c_{ij}z^{m_{ij}}$ for all $(i,j)$, so that $\rho_B=\infty$. Since $\bs B(|\alpha|)$ has~1 for dominant eigenvalue, Wielandt's theorem implies that there exists a matrix~$\bs D=\operatorname{diag}(e^{i\phi_1},\dots,e^{i\phi_m})$ such that $\bs B(\alpha)=\bs D\bs B(|\alpha|)\bs D^{-1}$. If $\theta=\arg\alpha$, extracting entries of this matrix identity yields
\[c_{i,j}\neq0\quad\Longrightarrow \quad m_{i,j}\theta=\phi_i-\phi_j\bmod 2\pi.\]
A circuit in the dependency graph of $\bs Y=\Id+\bs B\bs Y$ corresponds to a sequence of edges $(i_1,i_2),(i_2,i_3),\dots,(i_k,i_1)$ such that each edge starts from the endpoint of the previous one, the last edge ends at the starting point of the first one and each of the coefficients $c_{i_1,i_2},c_{i_2,i_3},\dotsc$ is nonzero. 
The system being irreducible, it contains at least one such circuit. Summing the identity above along a circuit, the right-hand side telescopes, which shows that $\theta$ is an integer multiple of $2\pi/\sum{m_{i,j}}$. Taking the gcd over all circuits concludes the proof.
The sum over all circuits starting at vertex~$k$ of the products of $c_{i,j}z^{m_{i,j}}$ over the edges of the circuit is exactly the $k$th entry in the diagonal of~$(\Id-\bs B(z))^{-1}$ (by the transfer-matrix method~\cite[\S4.7]{Stanley1986}). 

The statement of the theorem is obtained by taking the union of all these sets of potential singularities on the circle~$|z|=\rho_Y$.


\section{Dominant Singularities of Constructible Generating Functions}\label{sec:dominant} 
Recall that constructive systems (\cref{def:constructive}) are well-founded systems defined using the basic constructions $\Set,\Cyc,\Seq$. This specific structure, coupled with the results of the previous section, lets us design algorithms that determine the argument of the dominant singularities other than the radius of convergence.

In the case of nonlinear irreducible systems, \cref{thm:BBY72i} shows that the arguments of dominant singularities are of the form $\pi p/q$ where $q$ is the period of the generating functions. We first show how to compute this period. The general case of systems that are not necessarily irreducible then follows from the block-decomposition of systems into irreducible components and single nonrecursive equations, with some care for linear cases. This leads to the following analogue of \cref{thm:algo-radius} for dominant singularities.

\begin{theorem}
Given a \texttt{Constant Oracle} as in \cref{sec:constant-oracle}, 
Algorithm \textup{\textsf{\ref{algo:dominant}}}
computes a finite superset of the dominant singularities of the generating function solution of a constructive system in normal form.
\end{theorem}

The proof of this theorem splits into several cases depending on the type of the system, detailed in the following subsections.

\subsection{Valuations and Periods}

The computation of the periods is intertwined with the computation of valuations that was done in Algorithm \textsf{\ref{algo:wellfoundedandleadingterm}} on page~\pageref{algo:wellfoundedandleadingterm}.

Denoting by~$\Pi(y)$ the period of the power series~$y$, the following properties are simple consequences of the definitions~\cite[Prop.~2.9]{BellBurrisYeats2011}:
\begin{gather*}
\Pi(1)=\Pi(z)=0,\quad
\Pi(e^{y})=\Pi\left(\frac1{1-y}\right)=\Pi\left(\ln\frac1{1-y}\right)=
\operatorname{gcd}\left(\Pi(y),\val(y)\right),\\
\Pi(y_1+\dots+y_n)=\gcd\left(\Pi(y_1),\dots,\Pi(y_n),
\val(y_2)-\val(y_1),\dots,
\val(y_n)-\val(y_1)\right),\\
\Pi(y_1\dotsm y_n)=\operatorname{gcd}\left(\Pi(y_1),\dots,\Pi(y_n)\right).
\end{gather*}

\begin{algorithm}
\SetKwFunction{Periods}{Periods}
\SetKwFunction{gcd}{gcd}
\SetAlgoRefName{Periods}
\SetKwFunction{WFLT}{WellFoundedAndLeadingTerms}
\SetFuncSty{xAlCapSlimSty}
\caption{\label{algo:periods}Periods of the generating functions (component)}
\Input{$\bc{Y}=\bc{H}(\cal{Z},\bc{U},\bc{Y})$ a constructive component in normal form with $m$ equations\\
\phantom{{\bf Input: }}$\bc U = \bc G(\cal Z,\bc U)$ the sub-system defining $\bc U$
}
\Output{Vector of periods of the generating functions $\bs{Y}$ of $\bc{Y}$\vspace{10pt}}
    Let $\bc T = (\bc U,\bc Y)$ and $\bc F = (\bc G(\cal{Z},\bc U), \bc H(\cal{Z},\bc U, \bc Y))$\;
    $v_{1:\ell} :=$ valuations of $\bc T$ from \WFLT{$\bc T = \bc F(\cal{Z},\bc T)$}\;
    $q_{1:\ell}:=(0,\dots,0)$
    
    \Repeat{$\bs q=\bs \pi$}{
        $\bs \pi := \bs q$\;
        \For{$i=1$ to $\ell$}{
            \lIf{$\cal F_i=\mathcal{E}$ or $\cal F_i=\mathcal Z$}{$q_i:=0$}
            \lIf{$\cal F_i=\Phi(\cal T_j)$ with $\Phi\in\{\Set,\Cyc,\Seq\}$}{
                $q_i:=\gcd(\pi_j,v_j)$}
            \lIf{$\cal F_i=\cal T_{k_1}+\dots+\cal T_{k_n}$}{
                $q_i:=\gcd(\pi_{k_1},\dots,\pi_{k_n},v_{k_2}-v_{k_1},\dots,
                v_{k_n}-v_{k_1})$}
            \lIf{$\cal F_i=\cal T_{k_1}\times\dotsm\times\cal T_{k_n}$}{
                $q_i:=\gcd(\pi_{k_1},\dots,\pi_{k_n})$
            }
        }
        
    }
    \Return{the last $m$ entries of $\bs q$}
\end{algorithm}

Thus, following the same template as for the valuations,  Algorithm \textsf{\ref{algo:periods}} computes the periods. Termination is clear: as the system is constructive, it is zero-free and thus the valuations are all nonnegative integers. Then the iteration taking gcds can only reduce coordinates a finite number of times. This algorithm is slightly unnecessarily technical by taking as input a component and a sub-system, but it is convenient to recover only the periods of the component in other algorithms below.

Note that both algorithms apply to systems defining vectors as well as matrices. This is used in Algorithm \textsf{\ref{algo:dominant-irr}}.

\begin{example}Consider the system
\[\cal A=\cZ^2+\cZ^2\times \cal B^2,\quad \cal B=\cZ^3+ \cZ^3\times\cal A^2.\]
In order to simplify the presentation, we do not rewrite this system in normal form.
The valuations computed by Algorithm \textsf{\ref{algo:wellfoundedandleadingterm}} are $(v_A,v_B)=(2,3)$.
From there, the periods are given as solutions to the system
\[q_A=\gcd(6,q_B),\quad q_B=\gcd(4,q_A).\]
Fixed point iteration gives $(0,0)\mapsto(6,4)\mapsto(2,2)\mapsto(2,2).$ Thus both series have period~2.
\end{example}

\subsection{Nonrecursive Equations}

\begin{algorithm}[ht]
\SetFuncSty{xAlCapSlimSty}
\SetAlgoRefName{NR-DS}
\caption{Dominant singularities of the generating function \label{algo:dominant-nr}}
\Input{$\cal{Y}=\cal{H}(\cal{Z},\bc{U})$ a constructive component in normal form of one nonrecursive eqn;\\
\phantom{{\bf Input:}} $R$ radius of convergence of $Y$;\\
\phantom{{\bf Input:}} $\bs S_{\bc U}:=\bc U = \bc G(\cal Z,\bc U)$\hfill\texttt{\color{darkgray}\small// $\bs{U}$ is the generating function of its solution}\\
\phantom{{\bf Input:}} $\Delta_{\bc U} := \bigcup_{\cal{U}\in\bc{U}}\Delta_{\cal{U}}(R)$ with $\Delta_{\cal{U}}(R)\supseteq\{\alpha\mid Re^{i\alpha} \text{ singularity of }{U}\}$
}
\medskip
\Output{$\mathcal{S}\supseteq\{\alpha\mid Re^{i\alpha} \text{ singularity of }Y\}$\vspace{10pt}}
\lIf{$\cal{H}(\cal{Z},\bc{U})\in\{\cal{Z},1\}$}{\Return{$\emptyset$}}
\lElseIf{$\cal{H}(\cal{Z},\bc{U})=\Set(\cal{U})$, with $\cal{U}\in \bc{U}$}{\Return{$\Delta_{\cal{U}}(R)$}}
\ElseIf{$\cal{H}(\cal{Z},\bc{U})\in\{\Cyc(\cal{U}),\Seq(\cal{U})\}$, with $\cal{U}\in \bc{U}$}
{
$v :=$ valuation of $\cal U$ from \WFLT{$\bs S_{\bc U}$}\;\label{algo:dom:line6}
$q :=$ period of $\cal U$ from  \Periods{$\bs S_{\bc U},\emptyset$}\;\label{algo:dom:line7}
$p:=\gcd(v,q)$\;\label{algo:dom:line8}
\Return{$\{2k\pi/p\mid k=0,\dots,p-1\}$}
}
\lElse(\tcp*[f]{$\cal{H}(\cal{Z},\bc{U})$ is a product or a union  of $\bc{U}$}){\Return{$\Delta_{\bc U}$}}
\end{algorithm}

In the case of a nonrecursive equation in normal form,
the species~$\cal H(\cal{Z},\bc U)$ is one of
\[\{1,\mathcal Z,\cal{U}_1+\dots+\cal{U}_k,\cal{U}_1\times\dots\times\cal{U}_k,\Set(\cal{U}),\Seq(\cal{U}),\Cyc(\cal{U})\}\]
by \cref{def:normal_form}. In the first two cases, there are no singularities. In the next two cases (union and cartesian product), the set of singularities is a subset of the union of the sets of singularities of the~$U_i$.\footnote{It can be a strict subset. Consider for instance $\{\mathcal Y=\cal{U}_1\times\cal{U}_2, \cal{U}_1=1+\mathcal Z, \cal{U}_ 2 = \Seq(\mathcal Z\times\mathcal Z)\}$. The sets of singularities of the generating functions of $\cal{U}_1$ and $\cal{U}_2$ are $\emptyset$ and $\{-1,1\}$, but the only singularity of the generating function~$Y$ is~$1$.} In the case of $\Set$, the singularities are those of~$\bs U$. Finally, for the last two cases, dominant singularities occur when the generating function $U$ reaches~1 on the circle of radius~$R$. By \cref{lemma:positivity}, the argument of the dominant singularities  is of the form $\pi p/q$ with $q$ a divisor of the period of~$U$, which is computed in \crefrange{algo:dom:line6}{algo:dom:line8} of Algorithm \textsf{\ref{algo:dominant-nr}}.

\subsection{Irreducible Systems}
\begin{algorithm}[t]
\SetFuncSty{xAlCapSlimSty}
\SetAlgoRefName{Irr-DS}
\caption{Dominant singularities of the generating function \label{algo:dominant-irr}}
\Input{$\bc{Y}=\bc{H}(\cal{Z},\bc{U},\bc{Y})$ an irreducible constructive component in normal form with $\partial\bc{H}/\partial\cal{U}\neq0$;\;
\phantom{{\bf Input:}} $R$ radius of convergence of its generating function solution $\bs{Y}$;\\
\phantom{{\bf Input:}} $\bs S_{\bc U}:=\bc U = \bc G(\cal Z,\bc U)$\hfill\texttt{\color{darkgray}\small// $\bs{U}$ is the generating function of its solution}\\
\phantom{{\bf Input:}} $\Delta_{\bc U} := \bigcup_{\cal{U}\in\bc{U}}\Delta_{\cal{U}}(R)$ with $\Delta_{\cal{U}}(R)\supseteq\{\alpha\mid Re^{i\alpha} \text{ singularity of }{U}\}$
}
\medskip
\Output{$\mathcal{S}\supseteq\{\alpha\mid Re^{i\alpha} \text{ singularity of }\bs{Y}\}$\vspace{10pt}}
\If(\tcp*[f]{nonlinear case}){$\bpartial^2\bc{H}/\bpartial\bc{Y}^2\neq0$}{
$\bs{q}:=$ \Periods{$\bc{Y}=\bc{H}(\cal{Z},\bc{U},\bc{Y}),\bs S_{\bc U}$}\;\label{algo:dom:line13}
\Return{$\{2k\pi/q_1\mid k=0,\dots,q_1-1\}$}\tcp*[r]{\cref{thm:periods,thm:BBY72i}}}
\smallskip
\Else(\tcp*[f]{linear case}\label{algo:dom:line14}){
\lIf(\tcp*[f]{inequality}\label{line16}){$\Lambda_{\bs H}(R)<1$}{\Return $\Delta_{\bc U}$}
\Else{
Define $\bc{A}(\cal Z,\bc U)$, $\bc{B}(\cal Z,\bc U)$ by $\bc{H}=\bc{A}+\bc{B}\bc{Y}$\tcp*{$\bc{Y}=(\Id-\bc{B})^{-1}\bc{A}$}
$q:= \gcd( $\Periods{$\bc V = \bc{B}(\cal Z,\bc U)$, $\bs S_{\bc U}$})\label{algo:dom:lineMatrix}\;
\If{$q=0$}{
$\bs Q:= $ \Periods{$\bc Y=\Id+\bc B(\cal{Z},\bc{U})\bc Y$, $\bs S_{\bc U}$}\; 
$q:=\gcd(\operatorname{diagonal}(\bs Q))$}
\Return $\Delta_{\bc U} \cup \{2k\pi/q\mid k=0,\dots,q-1\}$
\tcp*[r]{\cref{thm:period-linear}\label{algo:dom:line23}}}
}
\end{algorithm}

\begin{algorithm}
\SetAlgoRefName{DominantSingularities}
\caption{Dominant singularities of the generating function (general case)\label{algo:dominant}}
\DontPrintSemicolon
\Input{$\bc{Y}=\bc{H}(\cal{Z},\bc{Y})$ a constructive system in normal form with $m$ equations }
\Output{Supersets of the dominant singularities of its generating function solution $\bs Y$.\vspace{10pt}}
    {$R:=\textsf{Radius}(\bc{Y}=\bc{H}(\cal{Z},\bc{Y}))$}\\
    Decompose the system into its $k$ strictly connected components (note that $i_k=m$)
    ${}\qquad\bc{Y}_{1:i_1}=\bc{H}_{1:i_1}(\cal{Z},\bc{Y}_{1:i_1}),
\dots,\bc{Y}_{i_{k-1}+1:i_k}=\bc{H}_{i_{k-1}+1:i_k}(\cal{Z},\bc{Y}_{1:i_k})$\;
    $\Delta:=\emptyset$\;
    \For{$j=1$ to $k$}{
        ${\bs S}_{\bc U}:=\bc{U}=\bc{H}_{1:i_{j-1}}(\cal Z,\bc{U})$\tcp*[f]{$i_0=0$}\\
        \If(\tcp*[f]{comparison already computed along with $R$}){\em\textsf{Radius}$({\bs S}_{\bc U})=R$}
        {

        \If(\tcp*[f]{at most one nonrecursive eqn}){$\partial\bs{H}_{i_{j-1}+1:i_j}/{\partial\bs{Y}_{i_{j-1}+1:i_j}}=0$}
        {$\Delta:=\textsf{NR-DS}(\mathcal{Y}=\cal{H}_{i_j}(\mathcal{Z},\bc{U},\cal{Y}), R, {\bs S}_{\bc U}, \Delta)$}
        \lElse{
        $\Delta:=\textsf{Irr-DS}(\bc{Y}=\bc{H}_{i_{j-1}+1:i_j}(\cal{Z},\bc{U}, \bc{Y}), R, {\bs S}_{\bc U}, \Delta)$
        }}
    }
    \Return{$R,\Delta$}\tcp*[f]{$\Delta\supseteq\{\alpha\mid Re^{i\alpha} \text{ singularity of }{\bs Y}\}$}
\end{algorithm}

In the case of irreducible \emph{nonlinear} systems the computation of periods above is sufficient to locate a finite set containing the dominant singularities of the generating functions by \cref{thm:BBY72i}. 

In the linear case, \cref{thm:period-linear} leads to considering another set of periods.

\begin{example}The simple variation 
\[\cal A=\cZ^2+\cZ^2\times \cal B,\quad \cal B=\cZ^3+ \cZ^3\times\cal A.\]
of the previous example can be treated in the same way. The valuations are again $(2,3)$. The system of equations for the periods is
\[q_A=\gcd(3,q_B),\quad q_B=\gcd(2,q_A),\]
from where an iteration shows that the periods are $(1,1)$. This is correct, but the generating function solution is
\[A(z)=z^2\frac{1+z^3}{1-z^5},\quad B(z)=z^3\frac{1+z^2}{1-z^5}\]
with denominators indicating dominant singularities at $\exp(2ik\pi/5)$ for $k=0,\dots,4$.

Indeed, the algorithm rewrites the system as $\bs{Y}=\bs S(z) +\bs J(z)\bs Y$, with
$$\bs S=\begin{pmatrix}z^2\\z^3\end{pmatrix} \mbox{ and }\bs J=\begin{pmatrix}0&z^2\\z^3&0\end{pmatrix},$$
and the periods of $(\Id-\bs J)^{-1}$ are those of the solutions of the system $\bs M=\Id+\bs J\bs M$. 
\cref{thm:algo-leading} shows that the valuations of~$\bs M$ satisfy
\[V=\begin{pmatrix}v_{11}&v_{12}\\ v_{21}&v_{22}\end{pmatrix}=
\begin{pmatrix}\min(0,2+v_{21})&2+v_{22}\\ 3+v_{11}&\min(0,3+v_{12})\end{pmatrix}
\]
and thus by iteration
\[V=\begin{pmatrix}0&2\\ 3&0\end{pmatrix}.\]
Next, the matrix of periods satisfies
\[\Pi=\begin{pmatrix}\pi_{11}&\pi_{12}\\ \pi_{21}&\pi_{22}\end{pmatrix}=
\begin{pmatrix}\gcd(\pi_{21},2+v_{21})&\pi_{22}\\ \pi_{11}&\gcd(\pi_{12},3+v_{12})\end{pmatrix}.
\]
Algorithm \textsf{\ref{algo:periods}} starts from the 0~matrix and computes the fixed point
\[\Pi=\begin{pmatrix}5&5\\ 5&5\end{pmatrix}\]
in two iterations. The gcd of the periods on the diagonal is~5 and the location of possible dominant singularities follows.

Note that in this example, the periods and the valuations (and therefore the leading terms) are computed for a matrix. This is done by flattening the entries of the matrix into a vector in order to use Algorithms~\textsf{\ref{algo:periods}} and \textsf{\ref{algo:wellfoundedandleadingterm}} with appropriate input.
\end{example}

\begin{proposition}\label{dom-irred}
Given a \texttt{Constant Oracle} as in \cref{sec:constant-oracle}, Algorithm \textup{\textsf{\ref{algo:dominant-irr}}} computes a finite superset of the dominant singularities of the generating function solution of a strongly connected component of a constructive system in normal form.
\end{proposition}
\begin{proof}
When the system is irreducible and nonlinear, \cref{thm:BBY72i,thm:periods} show that the singularities are deduced from the common period of the coordinates of~$\bs Y$, 
computed in \cref{algo:dom:line13}.

When the system is irreducible and linear, the location of the dominant singularities is given by \cref{thm:period-linear}, implemented in \crefrange{algo:dom:line14}{algo:dom:line23}. The decision on \cref{line16} is possible with the oracles at hand, since it was already done in the computation of the radius of convergence (\cref{algo:irr-rad:line4} of Algorithm \textsf{\ref{algo:radius}}). 
\end{proof}

\subsection{General Case}

As in the case of Algorithm \textsf{\ref{algo:radius}}, the computation of dominant singularities proceeds by finding them for each of the strongly connected components of the system in order, restricting to those whose generating function have the same radius of convergence as the full system. The result of each strongly connected component is used as input in the next one. This is done by Algorithm \textsf{\ref{algo:dominant}}, which relies on the Algorithms \textsf{\ref{algo:dominant-nr}} and \textsf{\ref{algo:dominant-irr}} depending on the nature of the component.  
The correctness of this approach follows from \cref{dom-irred}.

\begin{example}The generating function~$1/(1-z)+1/(1-2z^4)$ has periodicity~1 but~4 dominant singularities. It is the generating function of the last component of the following system in normal form
\[\{\cal A=\cZ^4,\,\cal B=\cal A+\cal A,\,\cal C=\Seq(\cal B),\,\cal D=\Seq(\mathcal Z),\,\cal E=\cal C+\cal D\}.\]
The algorithm starts by finding that the radius of convergence of the generating function is the positive root of~$1-2z^4=0$ and that it is the radius of convergence of~$C(z)$ and~$E(z)$. Then \textsf{NR-DS} is called with each of the two equations defining~$\cal C$ and~$\cal E$. For the first one, the valuation of~$\cal B$ is found to be~4 so that $S=\{0,2\pi/4,\pi/2,3\pi/4\}$ is returned. The second one is a union, resulting in the output~$S$.
\end{example}

\section{Singular Behavior}\label{sec:singular-behavior}
\begin{algorithm}[t]
\SetAlgoRefName{SingularExpansions}
\caption{Singular expansions of the generating functions\label{algo:singularexpansion}}
\Input{$\bc{Y}=\bc{H}(\cal{Z},\bc{Y})$ a constructive system in normal form with $m$ equations;\\
\phantom{{\bf Input:}} $p>1$ precision used in truncated series expansions.}
\Output{A superset of $\{(\sigma,\text{local expansion of $\bs Y$ at $\sigma$})\mid\sigma \text{ dominant singularity}\}$\vspace{10pt}}
    $R,\Delta:=\textsf{DominantSingularities}(\bc{Y}=\bc{H}(\cal{Z},\bc{Y}))$\\
    \lFor{$\alpha\in\Delta$}{$\bs S_\alpha:=\textsf{Localbehavior}(\bc{Y}=\bc{H}(\cal{Z},\bc{Y}),Re^{i\alpha},p)$}
    \Return{$\{(Re^{i\alpha},\bs S_\alpha)\mid \alpha\in\Delta\}$}
\end{algorithm}

While the radius of convergence of a power series conveys information on the exponential growth of its coefficients, singularity analysis obtains a more precise asymptotic behavior by first performing a more detailed analysis of the function in a neighborhood of its dominant singularities. 

In their book, Flajolet and Sedgewick ask for an automatic way to do this~\cite[p.~493]{FlajoletSedgewick2009}:
\begin{quote}\em
It would at least be desirable to determine directly, from
a positive (but reducible) system, the type of singular behavior of the solution, but
the systematic research involved in such a programme is yet to be carried out.
\end{quote}

An answer was provided by Banderier and Drmota~\cite{BanderierDrmota2015} for the case of positive systems of polynomial or entire functions under mild assumptions that prevent essential dominant singularities~\cite{BanderierDrmota2015}. The class of constructive systems is both smaller (not all positive entire functions occur) and larger (cycles and sequences are allowed). In this section, their results are extended to constructive systems. As exponentials and logarithms can now occur in the systems, there is a larger variety of possible singular behaviors that we classify, focusing on algebraic-logarithmic singularities.

The basis for the analysis of Banderier and Drmota and for ours is the classical Drmota-Lalley-Woods theorem~\cite{Drmota1997};\cite{Lalley1993};\cite{Woods1997};\cite[Thm.~VII.6]{FlajoletSedgewick2009}, stating that the dominant singularities of `nice' irreducible systems behave in a very predictable way: only square-root-type singularities occur. \Cref{sec:DLW} gives a version of this result in the case of nonlinear irreducible components of well-founded systems, with few extra hypotheses. This is then used to classify the possible behaviors of all systems at their dominant singularity on the positive real axis. In particular, we recover a \emph{gap} phenomenon that is known in the case of nonrecursive systems~\cite[Prop.~4.8]{FlajoletSalvyZimmermann1991}:
the behavior is either of the algebraic-logarithmic type (see \cref{def:alg-log})
or at least as large as $\exp(c\ln^2(1-z/\rho))$ for some $c>0$. Using this information on the dominant \emph{real} positive singularity, the possibilities at other dominant singularities can then be classified. There, nonreal exponents can appear, in contrast with the case of entire functions considered by Banderier and Drmota.

\subsection{Drmota-Lalley-Woods Theorem}\label{sec:DLW}
 We give a version of the Drmota-Lalley-Woods theorem suited to our context with an irreducible component inside a wider system, as an intermediate step to results for reducible systems. This is the source of the occurrence of Puiseux expansions in the system itself (coming from the bottom-up exploration of the DAG) rather than only in its solution.

\begin{definition}If $S(z)\in\mathbb C((z^{1/r}))$ is a power series
\[S(z)=\sum_{i\ge i_0\in\mathbb Z}c_iz^{i/r},\qquad c_{i_0}\neq0,\]
then by extension of \cref{def:support_valuation}, the rational number~$i_0/r$ is called the valuation of~$S$ and denoted~$\val_0S$. Similarly, we call valuation at~$\sigma$ and write~$\val_\sigma S$ for the valuation of $S(1-z/\sigma)\in\mathbb C(((1-z/\sigma)^{1/r}))$.
\end{definition}

\begin{restatable}{proposition}{DLW}\label{prop:DLW}
Let $\bs Y=\bs H(\mathcal Z,\bs U,\bs Y)$ be a nonlinear irreducible component of a constructive system, with $\bs H$ depending polynomially on $\bs U$. Let $\sigma\in\mathbb C$ be a dominant singularity of the generating function solution~$\bs Y$ and~${\rho=|\sigma|}$. 
If $\bpartial\bs H/\bpartial\bs U=\bs 0$, set~$r=1$; otherwise, assume that all the generating functions~$U_j$ such that ${\partial\bs H/\partial {U}_j \neq 0}$ have convergent Puiseux expansions in $\mathbb C\{(1-z/\sigma)^{1/r}\}$, and that the coefficient $\bs T$ of $(1-z/\sigma)^{1/r}$ in $\bs U(z)$ is nonzero with coordinates in $\mathbb R_{\le 0}$.

If $\Lambda_{\bs H}(\rho)<1$, then all coordinates of $\bs Y$ have an expansion in $\mathbb C\{(1-z/\sigma)^{1/r}\}$.
Otherwise, $\Lambda_{\bs H}(\rho)=1$ and each coordinate $Y_i$ of $\bs Y$ has an expansion in 
$\mathbb C\{(1-z/\sigma)^{1/2r}\}$ satisfying
\[Y_i(z)=\left(\frac\sigma\rho\right)^{\!\val_0(Y_i)}\left(Y_i(\rho)-C_i(1-z/\sigma)^{1/2r}+O((1-z/\sigma)^{2/r})\right),\]
with all coordinates of $\bs Y(\rho)$ and $\bs C$ in $\mathbb{R}_{>0}$. If $\bs v$ and $\bs w$ are respectively a left and a right eigenvector of~$\bpartial\bs H/\bpartial\bs Y(\rho,\bs Y(\rho),\bs U(\rho))$ for the eigenvalue~1 with positive coordinates then 
\begin{align}\label{eq:CandW}
\bs C&=\sqrt{\frac{2\bs v\cdot \bs W}{
\bs v\cdot\frac{\bpartial^2\bs H}{\bpartial\bs Y^2}(\rho,\bs Y(\rho),\bs U(\rho))(\bs w,\bs w)}}\,\bs w\quad\text{with}\\
\bs W&=\frac{\bpartial\bs H}{\bpartial\bs U}(\rho,\bs Y(\rho),\bs U(\rho))\cdot (-\bs T)+
\begin{cases}\rho\frac{\partial \bs H}{\partial z}(\rho,\bs Y(\rho),\bs U(\rho))&\text{if $r=1$},\\
0,\qquad &\text{otherwise.}
\label{eq:CandW2}
\end{cases}
\end{align}
\end{restatable}
Note that the condition that $\bs H$ depends polynomially on $\bs U$ is automatically satisfied for a constructive system in normal form (see \cref{def:normal_form}).

\noindent The proof is given in \cref{appendix:proofpropDLW}.

\subsection{Algebraic-Logarithmic Singularities}
The `transfer theorems' of Flajolet and Odlyzko~\cite{FlajoletOdlyzko1990a} give precise information on the asymptotic behavior of the coefficients of a generating function given its behavior in the neighborhood of its dominant singularities (see \cref{sec:transfer}).
Of particular importance is the class of functions with algebraic-logarithmic expansions in the following sense.
\begin{definition}\label{def:alg-log}
A function $f$ analytic in a slit neighborhood $\Delta$ of~$\sigma\in\mathbb C\setminus\{0\}$ is said to have an \emph{algebraic-logarithmic} singularity at~$\sigma$ if it satisfies
\[f\sim C\,(1-z/\sigma)^\alpha\ln^k\!\left(\frac{1}{1-z/\sigma}\right),\quad z\rightarrow \sigma\]
with $C$ and $\alpha$ in $\mathbb C$, $k\in\mathbb N$.
\end{definition}
The results of Flajolet and Odlyzko also hold for exponents~$k$ that are not natural integers and for iterated logarithms, but these will be proved not to appear in expansions of constructible generating functions (\cref{thm:alg-log-real,thm:alg-log-dom}).

Constructible generating functions with finite radius of convergence satisfy a \emph{gap} property: either they have algebraic-logarithmic expansions at their dominant singularities or they have superpolynomial growth there, in the following sense.
\begin{definition}
The function $f$ is said to have \emph{superpolynomial growth} at $\sigma$ if
\[\exists c>0 \text{ s.t. } |\ln f|>c\ln^2\!\left(|1-z/\sigma|\right)\text{ for $z$ sufficiently close to~$\sigma$}.\]
\end{definition}
This gap property is known  for the special case of nonrecursive equations~\cite[Prop.~4.8]{FlajoletSalvyZimmermann1991}; \cref{thm:alg-log-real,thm:alg-log-dom} show that it holds for systems, while giving more information on the full local expansion. 

\subsubsection{Real Positive Dominant Singularity}

The positivity of the coefficients of generating functions at the origin induces strong constraints on their behavior at their dominant real positive singularity.
\begin{theorem}\label{thm:alg-log-real}
Let $\bc Y=\bc H(\mathcal Z,\bc Y)$ be a 
constructive system and let $\rho>0$ be the real dominant singularity of the generating function~$\bs Y$. Let $Z=1-z/\rho$ and $L=\ln(1/Z)$. As $z\rightarrow\rho-$, each coordinate $Y_i$ of $\bs Y$ is either of superpolynomial growth at~$\rho$ or can be expressed as a finite sum of expansions of one of the types
\begin{align*} 
L^k\times P,&\quad k\in\mathbb N,P\in\mathbb R\{Z^{1/r}\},\val P=0;\\
Z^\alpha L^k\sum_{i\ge0}p_i(L)Z^{i/r},&\quad \alpha\in\mathbb R_{<0},k\in\mathbb N,p_i\in\mathbb R[L],\deg p_i\le i,p_0\neq0.
\end{align*}
In both cases, $r$ is of the form~$2^m$ with $m\in\mathbb N$.

Moreover, when the system does not use the $\Cyc$ construction, then $L$ does not appear in the expansions ($k=0$ and $\deg p_i=0$), the exponents~$\alpha$ that occur in the second type of expansion are rational with denominator~$r$ and the corresponding sum belongs to~$\mathbb R\{Z^{1/r}\}$.
\end{theorem}
The second part of the theorem recovers the result of Banderier and Drmota~\cite[Thm.~5.3]{BanderierDrmota2015} for constructible generating functions. Examples where $r$ is 4~or~8 are given by supertrees and supersupertrees~\cite[p.~412--414]{FlajoletSedgewick2009}.

\begin{example}\label{ex:transcendental-exponent}This is a simplified variant of~\cite[Ex.~4.9]{FlajoletSalvyZimmermann1991}. It gives an example of an expansion of the second type, with an exponent that is not algebraic.
The species
\[\mathcal Y=\Set\!\left(\mathcal Z\times\Cyc\!\left(\cZ^2+\cZ^2\right)\times\Set\!\left(\Seq(\cZ^3)\right)\right)\]
has for generating function
\[Y=\exp\!\left(z\ln\!\left(\frac1{1-2z^2}\right)\exp\!\left(\frac1{1-z^3}\right)\right).\]
At the dominant singularity~$z=1/\sqrt2$, its behavior is
\[Y(z)=\frac{2^\alpha}{(1-z\sqrt2)^\alpha}(1+O(1-z\sqrt2)),\quad\alpha=\frac{\sqrt2}2\exp\!\left(\frac{8+2\sqrt2}{7}\right).\]
The exponent $\alpha$ is transcendental by Lindemann's theorem.
\end{example}

\Cref{thm:alg-log-real} is a classification theorem. Its main point is that no other behavior is possible.

\begin{notation}
The set of functions with superpolynomial growth at~$\rho>0$ is denoted $\mathcal E^{\exp}_{\rho}$. Similarly, $\mathcal E_{\rho,k}^{\ln}$ denote the sets of functions whose behavior at $\rho$ is a finite sum of expansions of the first type in the theorem, with $k$ the maximal value in these expansions. (Functions that are analytic at~$\rho$ belong to~$\mathcal E_{\rho,0}^{\ln}$ with $r=1$.) Finally, $\mathcal E^{\al}_{\rho,\alpha,k}$ denotes the set of functions with an algebraic-logarithmic singularity at~$\rho$ with $\alpha\neq0$ whose behavior at $\rho$ is a finite sum of expansions of the first and second type, with $\alpha$ the minimal value and $k$ the maximal one for this value of~$\alpha$.

Also, we write
$\mathcal E^{\ln}_\rho=\cup_k\mathcal E^{\ln}_{\rho,k}$ and 
$\mathcal E^{\al}_\rho=\cup_{\alpha,k}\mathcal E^{\al}_{\rho,\alpha,k}$.   
\end{notation}

\begin{proof}[Proof of \cref{thm:alg-log-real}] Without loss of generality, we consider the system to be in normal form: the coordinates~$Y_i$ of the original system form a subset of those of the normalized system. 

The proof is by induction on the strongly connected components~$\bc Y=\bc H(\mathcal Z,\bc U,\bc Y)$ of the decomposition of the system. The induction hypothesis is that $\bs U\in\mathcal E_\rho^{\exp}\cup\mathcal E_\rho^{\ln}\cup\mathcal E_\rho^{\al}$ and moreover, if a nonconstant component $U_i$ has an expansion in~$\mathcal E^{\ln}_{\rho,0}$, then the coefficient of $Z^{1/r}$ in that expansion is negative (and in particular nonzero).

\medskip\noindent\emph{Irreducible component.}
In the nonlinear case, the result is a consequence of \cref{prop:DLW}. 

Otherwise, the component is of the form~$\bc Y=\bc A+\bc B\bc Y$, with $\bs A,\bs B$ having entries with expansions in~$\mathbb R\{\{Z^{1/r}\}\}$, which is a field. Therefore the coordinates~$Y_i$ of the solution also belong to~$\mathbb R\{\{Z^{1/r}\}\}$. By \cref{lemma:finiteornot}, either they all tend to infinity or they all have a finite limit. In the first case, they belong to~$\mathcal E^{\al}_{\rho,\alpha,0}$ with $\alpha$ rational and with denominator~$r$. The second case is when $\Lambda(\bs B)<1$. The expansion is given by
\[\bs Y=(\bs I-\bs B)^{-1}\bs A
=\left(\bs I+\bs C\bs D+O((\bs C\bs D)^2)\right) \bs C\bs A
\]
with $\bs C=(\bs I-\bs B(\rho))^{-1}$, $\bs D=\bs B-\bs B(\rho)$. By irreducibility, the entries of $\bs C$ are all positive; by the induction hypothesis at least one of the entries of $\bs A$ and $\bs D$ has a nonzero coefficient of~$Z^{1/r}$, which is moreover negative, and thus so do all the entries in the product.

\medskip\noindent\emph{Nonrecursive component.}
The system is $\cal Y=\mathcal H(\mathcal Z,\mathcal U)$.

If $\mathcal H\in\{1,\mathcal Z\}$, then the expansions are~$1$ and~$\rho(1-Z)$ in $\mathcal E^{\ln}_{\rho,0}$, that satisfy the induction hypothesis.

\smallskip{\em Sums and products.}
Let $f,g$ be analytic functions in a slit neighborhood of~$\rho$. Then the cases when $\mathcal H$ is a sum or a product follow from
\begin{align*}
f\in\mathcal E^{\exp}_\rho,g\in \mathcal E_\rho^{\exp}\cup\mathcal E_\rho^{\ln}\cup\mathcal E_\rho^{\al}&\Rightarrow \{fg,f+g\}\subset\mathcal E^{\exp}_\rho;\\
f\in\mathcal E^{\al}_{\rho},g\in\mathcal E_{\rho}^{\ln}\cup\mathcal E_{\rho}^{\al}&\Rightarrow \{fg,f+g\}\subset\mathcal E^{\al}_\rho;\\
f\in\mathcal E^{\ln}_{\rho},g\in \mathcal E_{\rho}^{\ln}&\Rightarrow\{fg,f+g\}\subset\mathcal E^{\ln}_{\rho}.
\end{align*}
In the last case, if both $f,g$ belong to~$\mathcal E^{\ln}_{\rho,0}$ with a negative coefficient of~$Z^{1/r}$, then so does their product and sum.

\smallskip{\em Sets.} If $\mathcal H=\Set$ with $U\in \mathcal E_\rho^{\exp}\cup\mathcal E_\rho^{\ln}$ or $U\in\mathcal E^{\ln}_{\rho,k}\text{ and }k\ge2$, then $Y=\exp(U)\in\mathcal E^{\exp}_\rho$.
If $U\in\mathcal E^{\ln}_{\rho,0}$, then $\textbf{}\exp(U)\in\mathcal E^{\ln}_{\rho,0}$ by composition with a convergent Taylor expansion. The coefficient of $Z^{1/r}$ being negative is preserved by the exponential.

The remaining case for $\Set(\mathcal U)$ is when $U\in\mathcal E^{\ln}_{\rho,1}$ (induced by a $\Cyc$ in a previous component). Then we have
\[
U(z)=aL+LP_{>0} + Q,\quad z\rightarrow\rho,
\]
with $a>0$ and $P_{>0}\in\mathbb R\{Z^{1/r}\}$, $\val P_{>0}>0$ and a possibly nonzero~$Q\in\mathbb R\{Z^{1/r}\}$ with $\val_\rho Q=0$. 
The exponential of $U$ is $Z^{-a}\exp(LP_{>0})\exp(Q)$ and thus belongs to~$\mathcal E^{\al}_{\rho}$.

\smallskip{\em Sequences and cycles.}
When $\mathcal Y=\mathcal H(\mathcal U)$ with $\mathcal H\in\{\Seq,\Cyc\}$, the generating function~$U$ belongs to~$\mathcal E_{\rho,0}^{\ln}$ and its limit at $\rho$ is at most~1. If $\lim U<1$, then by composition with a Taylor series, $Y\in\mathcal E_{\rho,0}^{\ln}$ and the coefficient of~$Z^{1/r}$ is negative. Otherwise,
$U$ has an expansion
\[
U(z)=1-cZ^{1/r}P,\quad c>0,\val_\rho P=0,\qquad z\rightarrow\rho
\]
and the leading coefficient of~$P$ is~1.
Then, for $\Seq(\mathcal U)$, 
\[
\frac1{1-U(z)}=\frac1cZ^{-1/r}\frac1P
\]
and $Q=1/P\in\mathbb R\{Z^{1/r}\}$ has $\val_\rho Q=0$ so that this is an expansion in $\mathcal E^{\al}_{\rho,0}$. In the case of $\Cyc(\mathcal U)$ one obtains
\[
\ln\frac1{1-U(z)}=\frac 1rL-\ln c+\ln\frac1P
\]
which is the sum of the expansion $(1/r)L\in\mathcal E^{\ln}_\rho$ and the expansion $-\ln c-\ln(P)\in\mathbb R\{Z^{1/r}\}$, concluding the proof.
\end{proof}

\subsubsection{Other Dominant Singularities}
In view of the application of transfer theorems to the local expansions in the next section, we focus here on generating functions whose behavior along the positive real axis is algebraic-logarithmic. This behavior persists at the other dominant singularities, i.e., they do not have superpolynomial growth there. Still cancellations and complex coefficients can occur in the expansions and in the exponents at nonreal singularities, leading to more diversity in the possible local behaviors. 
\begin{theorem}\label{thm:alg-log-dom}
Let $\bc Y=\bc H(\mathcal Z,\bc Y)$ be a 
constructive system and let $\sigma$ be a dominant singularity of the generating function~$\bs Y$ with $\rho=|\sigma|$. Let $Z=1-z/\sigma$ and $L=\ln(1/Z)$.
Let $f$ be a coordinate of $\bs Y$ that has an algebraic-logarithmic singularity at~$\rho$. Then as $z\rightarrow\sigma$, its behavior can be expressed as a finite sum of expansions of one of the types
\begin{align}
L^{k}\times P,&\quad k\in\mathbb N,P\in\mathbb C\{Z^{1/r}\},\val_\sigma P\ge0;\label{eq:type1}\\
Z^\alpha L^k\sum_{i\ge 0}p_i(L)Z^{i/r},&\quad\alpha\in\mathbb C, k\in \mathbb N,p_i\in\mathbb C[L],\deg p_i\le i.\label{eq:type2}
\end{align}
In both cases, $r$ is of the form~$2^m$ with $m\in\mathbb N$.

When the system does not use the $\Cyc$ construction, then $L$ does not appear in the expansions, the exponents~$\alpha$ that occur in the second type of expansion are rational with denominator~$r$ and the corresponding sum belongs to~$\mathbb C\{Z^{1/r}\}$.
\end{theorem}
The main differences with the real positive case of \cref{thm:alg-log-real} is that in the first type of expansions, the leading term may have~$Z$ to a positive exponent and in the second type, $\alpha$ can be a nonreal complex number and the coefficient $p_0$ may be 0.

\begin{example}\label{ex:imaginary-exponent}Consider the species $\mathcal Y=\Set(\cal Z\times\Cyc(\cZ^4))$. Its generating function $Y(z)$ has a complex exponent at its singularity~$i(=\sqrt{-1})$:
\[Y(z)=\exp\!\left(z\ln\!\left(\frac1{1-z^4}\right)\right)\sim 4^{-i}(1-z/i)^{-i},\quad z\rightarrow i.\]
This imaginary exponent is the basis for a term in $\cos(\ln n+\phi)$
in the asymptotic behavior given in \cref{ex:ex27} below.
\end{example}

\begin{notation}
The set of functions whose behavior at~$\sigma$ is a finite sum of expansions of the first type is denoted~$\hat{\mathcal E}^{\ln}_{\sigma,k}$, with $k$ the maximal value reached in these expansions. Those whose behavior is a finite sum of expansions of the first and second type is denoted~$\hat{\mathcal E}^{\al}_{\sigma,\alpha,k}$ with $\alpha$ and~$k$ as above. Again, $\hat{\mathcal E}^{\ln}_\sigma=\cup_k\hat{\mathcal E}^{\ln}_{\sigma,k}$ and $\hat{\mathcal E}^{\al}_\sigma=\cup_{\alpha,k}\hat{\mathcal E}^{\al}_{\sigma,\alpha,k}$.
\end{notation}
\begin{proof}[Proof of \cref{thm:alg-log-dom}]
The proof is similar to the previous one, by an induction on the strongly connected components~$\bc Y=\bc H(\mathcal Z,\bc U,\bc Y)$ of the decomposition of the system. The new ingredient is that the behavior of~$f$ at~$\sigma$ is dominated by its behavior at~$\rho$, given by \cref{thm:alg-log-real}.
In particular if~$f$ is analytic at~$\rho$, it is also analytic at~$\sigma$.
The induction hypothesis is that $\bs U\in\hat{\mathcal E}_\sigma^{\ln}\cup\hat{\mathcal E}_\sigma^{\al}$.

\medskip\noindent\emph{Irreducible component.}
The reasoning is the same as in the real case.
If the component is linear, the results follows from the fact that~$\mathbb C\{\{Z^{1/r}\}\}$ is a field. 
If the component is not linear, then if $\Lambda_\bs H(\sigma)<1$, the result is a consequence of \cref{prop:DLW}. Otherwise, $\Lambda_\bs H(\sigma)=\Lambda_\bs H(\rho)=1$ and $\sigma$ is a singularity of~$\bs Y$. By \cref{thm:periods}, all coordinates $Y_i$ are of the form $Y_i(z)=z^{\val Y_i}V_i(z^q)$ with $q$ such that $\arg\sigma=2\pi k/q$ for some~$k\in\mathbb N$. Thus, their behavior as~$z\rightarrow\sigma$ is given by their behavior as~$z\rightarrow\rho$; \cref{thm:alg-log-real} then implies that $\bs Y\in\hat{\mathcal E}_{\sigma}^{\ln}$.

\medskip
\noindent\emph{Nonrecursive component.}
The system is $\cal Y=\mathcal H(\mathcal Z,\mathcal U)$.

If $\mathcal H\in\{1,\mathcal Z\}$, then the expansions are~$1$ and~$\rho(1-Z)$ in $\hat{\mathcal E}^{\ln}_{\sigma,0}$.

\smallskip{\em Sums and products.}
The cases when $\mathcal H$ is a sum or product follow from
\begin{align*}
f\in\hat{\mathcal E}^{\al}_{\sigma},g\in\hat{\mathcal E}_{\sigma}^{\ln}\cup\hat{\mathcal E}_{\sigma}^{\al}&\Rightarrow \{fg,f+g\}\subset\hat{\mathcal E}^{\al}_\sigma;\\
f\in\hat{\mathcal E}^{\ln}_{\sigma},g\in\hat{\mathcal E}_{\sigma}^{\ln}&\Rightarrow\{fg,f+g\}\subset\hat{\mathcal E}^{\ln}_{\sigma}.
\end{align*}

\smallskip{\em Sets.}
If $\mathcal H=\Set$ and $f\in\mathcal E^{\ln}_\rho$, then $U\in\mathcal E_{\rho,0}^{\ln}$. Following the different cases in the proof of \cref{thm:alg-log-dom}, this implies that $U$ belongs to $\hat{\mathcal E}_{\sigma,0}^{\ln}$ and therefore so does $Y=\exp(U)$ by composition with an analytic function.

The remaining case for $\Set(\mathcal U)$ is when $f\in\mathcal E_\rho^{\al}$. Then $U\in\mathcal E_{\rho,1}^{\ln}$. By the induction hypothesis, this implies $U\in\hat{\mathcal E}_{\sigma,1}^{\ln}\cup\hat{\mathcal E}_{\sigma,0}^{\ln}$. If $U\in\hat{\mathcal E}_{\sigma,0}^{\ln}$ then by composition with $\exp$ also $f\in\hat{\mathcal E}_{\sigma,0}^{\ln}$. Otherwise, 
\[U(z)=aL+LP_{>0}\,(+Q),\quad z\rightarrow\rho,\]
with $a\in\mathbb C$ and $P_{>0}\in\mathbb C\{Z^{1/r}\}$, $\val P_{>0}>0$ and a possibly nonzero~$Q\in\mathbb C\{Z^{1/r}\}$ with $\val_\rho Q\ge0$. 
The exponential of $U$ is $Z^{-a}\exp(LP_{>0})\exp(Q)$ and thus belongs to~$\mathcal E^{\al}_{\sigma}$. This is where exponents that are not rational with denominator a power of~2 can occur when~$\Cyc$ is used.

\smallskip{\em Sequences and Cycles.} 
If $\mathcal Y=\mathcal H(\mathcal U)$ with $\mathcal H\in\{\Seq,\Cyc\}$ and $Y$ is not analytic at~$\sigma$ while it belongs to~$\mathcal E^{\al}_\rho$, then it has limit~1 at~$\rho$ and at~$\sigma$. By \cref{thm:BBY72i,thm:periods}, $\sigma=\rho\exp(2i\pi k/q)$ with $k\in\{0,\dots,q-1\}$ and $U(z)=z^{\val_0 U}V(z^q)$ with~$V$ a power series at~0, which is then analytic for $|z|\le\rho^q$. Since~$\lim U=1$ at~$z=\sigma$, it follows that $\arg\sigma^{\val_0U}=0$, which implies that~$\val_0U=0$. Thus the behavior of~$U$ at~$\sigma$ is given by its behavior at~$\rho$ and by \cref{thm:alg-log-real},
\[U(z)=1-cZ^{1/r}P,\quad c>0,P(0)=1\]
with $Z=1-z/\sigma$, so that $1/(1-U(z))\in\hat{\mathcal E}_\sigma^{\al}$ and $\ln(1/(1-U(z)))\in\hat{\mathcal E}_\sigma^{\ln}$. 
\end{proof}

\subsection{Algorithm}
\begin{algorithm}[H]
\SetAlgoRefName{Irr-Localbehavior}
\caption{Local expansions of the generating functions\label{algo:localbehavior-irred}}

\Input{$\bc{Y}=\bc{H}_{1:m}(\cal{Z},\bc{U},\bc{Y})$ a constructive component in normal form with $\partial\bc{H}/\partial\cal{U}\neq0$ for all $\cal{U}\in\bc{U}$;\\
\phantom{{\bf Input:}} $\sigma\in\mathbb C$ with $|\sigma|\le\rho=$radius of convergence of its generating function solution $\bs Y$;\\
\phantom{{\bf Input:}} $\{S_{\cal{U}}\mid \cal{U}\in \bc U\}$ local expansion of $\bs U$ at~$\sigma$\\
\phantom{{\bf Input:}} given as $\Exp$ or linear combinations of expansions of types~\eqref{eq:type1},\eqref{eq:type2};\\ 
\phantom{{\bf Input:}} $p>1$ precision used in truncated series expansions.
}
\Output{$\{S_{\cal Y}\mid\cal Y\in\bc Y\} $ truncated local expansion of $\bs{Y}$ at $\sigma$,\\
\phantom{{\bf Input:}} or {\sc Exp} if $\bs Y$ has a superpolynomial behavior at~$\sigma$\vspace{10pt} 
}
\tcc*[r]{$r$ is given by the local expansion of $\bs U$ in~$\mathbb C\{(1-z/\sigma)^{1/r}\}$; $r=1$ when $\bs U=\emptyset$}
\If(\tcp*[f]{linear case}){$\bpartial^2\bc{H}/\bpartial\bc{Y}^2=0$}{
Define $\bc{A},\bc{B}$ by $\bc{H}=\bc{A}+\bc{B}\bc{Y}$\tcp*{$\bc{Y}=(\Id-\bc{B})^{-1}\bc{A}$}
\Return{$(\Id-S_{\bc B}+O((1-z/\sigma)^{p/r}))^{-1}S_{\bc A}$}\tcp*[f]{Oracle}}
\Else(\tcc*[f]{$\Lambda_{\bs H}(\sigma):=\max|\lambda|\text{\ s.t.\ }\det(\bpartial\bs H/\bpartial\bs Y(\sigma,\bs{U}(\sigma),\bs{Y}(\sigma))-\lambda\Id)=0$}){
\If(\tcp*[f]{non singular, nonlinear case, Oracle, $\bs Y$ finite at~$\sigma$}\label{algline:nonsing-nonlin}){$\Lambda_{\bs H}(\sigma)<1$}
{
Compute expansion of $\bs Y$ by undetermined coefficients or 
Newton iteration\\
\Return{$\bs Y(\sigma)+\bs c_1(1-z/\sigma)^{\alpha_1/r}+\dots+O((1-z/\sigma)^{p/r})$}
}
\Else(\tcp*[f]{singular, nonlinear case, $\Lambda_\bs H(\sigma)=1,|\sigma|=\rho$}\label{algline:sing-nonlin}){
\tcc*{Halving of the exponent, following the proof of \cref{prop:DLW}}
Compute the value $\bs Y(\sigma)$ by \textsf{Algorithm \ref{algo:irr-radius}};\\
Let $\bs T_{\bc U}(v)=\bs S_{\bc U}(\rho(1-v^2))$;\tcp*[f]{$v$ stands for $(1-z/\rho)^{1/(2r)}$}\\
Let $\tilde{\bs H}_{2:m}(\bs Y_{2:m},\delta,v)=\bs H_{2:m}(\rho(1-v^2),\bs T_{\bc U}(v),Y_1(\rho)+\delta v,\bs Y_{2:m})$;\\
\smallskip
\tcc*{for instance by Newton iteration over power series}
Compute $\bs S_{2:m}(v,\delta)$, the unique power series in $\mathbb C[\delta][[v]]$ such that\\
\qquad $\bs S_{2:m}(v,\delta)=\tilde{\bs H}_{2:m}(\bs S_{2:m}(v,\delta),\delta,v)$ and $\bs S_{2:m}(0)=\bs Y_{2:m}(\rho)$, \\
\smallskip
\tcc*{$Y_1(\rho)=H_1(\rho,\bs T_{\bc U}(0),\bs Y(\rho))$, $\partial H_1/\partial Y_1(\rho,\bs T_{\bc U}(0),\bs Y(\rho))=1$}
\tcc*{and $a_2$ of the form $c-d\delta^2$ by \cref{eq:Schur}}{Compute the expansion $E(v,\delta)=a_2v^2+a_3v^3+\dotsb$ of}\\
{\qquad $H_1(\rho(1-v^2),\bs T_{\bc U}(v),Y_1(\rho)+\delta v,S_{2:m}(v,\delta))-Y_m(\rho)-\delta v;$}\\
\smallskip
\tcc*{for instance by Newton iteration over power series}
Compute the unique power series~$\Delta(v)$ with $\Delta(0)=-\sqrt{c/d}$ and $E(v,\Delta(v))=0$,\\
Set $\Sigma_1(v)=Y_1(\rho)+v\Delta(v)$; $\Sigma_{2:m}(v)=S_{2:m}(v,\Delta(v))$ expanded in powers of~$v$;\\
\Return{$\{(\sigma/\rho)^{\val_0 Y_i}\Sigma_i((1-z/\sigma)^{1/(2r)})\mid i=1,\dots,m\}$}
}
}
\end{algorithm}

\Cref{thm:alg-log-real,thm:alg-log-dom} become effective in Algorithms \textsf{\ref{algo:localbehavior-irred}}, \textsf{\ref{algo:localbehavior-nonrecursive}} and \textsf{\ref{algo:localbehavior}}, whose properties are summarized in the following.
\begin{theorem}
Let $\bc Y=\bc H(\cal Z,\bc Y)$ be a 
constructive system in normal form (\cref{def:irreducible-component-com}) and assume that $\sigma\in\mathbb C$ lies inside or on the boundary of the domain of convergence of the generating function~$\bs Y$.
Algorithm \textsf{\ref{algo:localbehavior}} computes truncations of the expansions of all coordinates that are analytic or have an algorithmic-logarithmic singularity at~$\sigma$ (\cref{def:alg-log}). It returns {\sc Exp} for the other coordinates.
\end{theorem}
\begin{proof}
The algorithm first decomposes the system into its irreducible components. It then takes each component in turn and deals with systems where $\bc Y$ may depend on auxiliary variables $\bc U$ coming from earlier components. The proof of correctness is by induction and the induction hypothesis is that when $\bc U$ is not empty, all its coordinates have expansions at~$\sigma$ given by \cref{thm:alg-log-real,thm:alg-log-dom}, for which truncations have been computed recursively. 

The nonrecursive cases are dealt with by Algorithm \textsf{\ref{algo:localbehavior-nonrecursive}}, the recursive ones by Algorithm \textsf{\ref{algo:localbehavior-irred}}. In both cases, the algorithm follows the steps of the proofs of \cref{prop:DLW,thm:alg-log-real,thm:alg-log-dom}.
\end{proof}
Note that this result does not require any help from an oracle: all difficult decisions have been made during the computation of the dominant singularities.

\begin{algorithm}[t]

\SetAlgoRefName{NR-Localbehavior}
\SetAlgoCaptionLayout{textbf}
\caption{Local expansions of the generating functions \label{algo:localbehavior-nonrecursive}}
\Input{$\cal{Y}=\cal{H}(\cal{Z},\bc{U})$ a constructive component in normal form of one nonrecursive equation with $\partial\bc{H}/\partial\cal{U}\neq0$ for all $\cal{U}\in\bc{U}$;\\
\phantom{{\bf Input:}} $\sigma\in\mathbb C$ with $|\sigma|\le\rho=$radius of convergence of its generating function solution $\bs Y$;\\
\phantom{{\bf Input:}} $\{S_{\cal{U}}\mid \cal{U}\in \bc U\}$ local expansion of $\bs U$ at~$\sigma$ \\
\phantom{{\bf Input:}} given as $\Exp$ or linear combinations of expansions of types~\eqref{eq:type1},\eqref{eq:type2};\\ 
\phantom{{\bf Input:}} $p>1$ precision used in truncated series expansions.
}
\Output{$S_{\cal Y} $ truncated local expansion of $Y$ at $\sigma$,\\
\phantom{{\bf Input:}} or {\sc Exp} if $Y$ has a superpolynomial behavior at~$\sigma$\vspace{10pt}}

{
\tcc*[r]{$r$ is given by the local expansion of $\bs U$ in~$\mathbb C\{(1-z/\sigma)^{1/r}\}$; $r=1$ when $\bs U=\emptyset$}
\lIf{$\Exp\in\{S_\mathcal U\mid\mathcal U\in\bc U\}$}{\Return{$\Exp$}}
\lElseIf{$\cal H(\cal{Z},\bc{U})=1$}{\Return{$1$}}
\lElseIf{$\cal{H}(\cal{Z},\bc{U})=\cal{Z}$}{\Return{$\sigma-\sigma(1-z/\sigma)$}}
\lElseIf{$\cal H\in\{\times,+\}$}{\Return{$\prod_{\cal{U}\in\bc{U}}S_{\cal{U}}$ or $\sum_{\cal{U}\in\bc{U}}S_{\cal{U}}$, depending on $\cal H$}}}
\ElseIf{$\cal H(\cal Z,\cal U)=\Set(\cal U)$ and $S_\cal U\rightarrow\infty$ as~$z\rightarrow\sigma$}{
\lIf{$S_\cal U$ is of type \eqref{eq:type2} or of type \eqref{eq:type1} with $k>1$}{\Return{$\Exp$}}
Let $a\in \mathbb C$, $P,Q$ in $\mathbb C\{Z^{1/r}\}$ be such that 
$S_\cal U=aL+LP+Q$ with $\val P>0$\\
Compute the expansion $E=\exp(LP)\exp(Q)$\tcp*[f]{composition with $\exp$ at~0}\\
\Return{$Z^{-a}E$}
}
\ElseIf(\tcp*[f]{$S_\mathcal U\rightarrow1$ as $z\rightarrow\sigma$}){$\cal{H}(\cal{Z},\bc{U})\in\{\Cyc(\cal{U}),\Seq(\cal{U})\}$ and $|\sigma|=\rho$}
{
Let $c>0$ and $P\in\mathbb C\{Z^{1/r}\}$ be such that $S_u=1-cZ^{1/r}P$ with $P(0)=1$\\
Compute $Q:=1/P$ by power series inversion\\
\lIf{$\cal H=\Seq$}{\Return{$Z^{-1/r}Q/c$}}
\lElse(\tcp*[f]{composition with $\ln$ at~1}){\Return $L/r-\ln c+\ln Q$}
}
\Else(\tcp*[f]{$\cal{H}(\cal{Z},\bc{U})\in\{\Set(\cal{U}),\Cyc(\cal{U}),\Seq(\cal{U})\}$ and $\lim S_\cal U$ finite})
{\lIf{$\cal H(\cal Z,\cal U)=\Set(\cal U)$}{$f:=s\mapsto\exp(s)$}
\lElseIf{$\cal H(\cal Z,\cal U)=\Seq(\cal U)$}{$f:=s\mapsto1/(1-s)$}
\lElse{$f:=s\mapsto\ln(1/(1-s))$}
{\Return{$f(S_{\cal{U}}+O((1-z/\sigma)^{p/r}))$}\tcp*[f]{compose series expansions}}}
\end{algorithm}

\begin{algorithm}[ht]
\SetAlgoRefName{Localbehavior}
\caption{Local expansions of the generating functions  (general case)\label{algo:localbehavior}}
\Input{$\bc{Y}=\bc{H}(\cal{Z},\bc{Y})$ a constructive system in normal form with $m$ equations;\\
\phantom{{\bf Input:}} $\sigma\in\mathbb C$ with $|\sigma|$ the radius of convergence of its generating function solution $\bs Y$;\\
\phantom{{\bf Input:}} $p>1$ precision used in truncated series expansions.}
\Output{$\{S_{\cal Y}\mid\cal Y\in\bc Y\} $ truncated local expansion of $\bs{Y}$ at $\sigma$,\\
\phantom{{\bf Output:}} or {\sc Exp} if $\bs Y$ has a superpolynomial behavior at~$\sigma$\vspace{10pt}}
    Decompose the system into its $k$ strictly connected components:
    ${}\quad\bc{Y}_{1:i_1}=\bc{H}_{1:i_1}(\cal{Z},\bc{Y}_{1:i_1}),
\dots,\bc{Y}_{i_{k-1}+1:i_k}=\bc{H}_{i_{k-1}+1:i_k}(\cal{Z},\bc{Y}_{1:i_k})$

    \For(\tcp*[f]{$i_0=0$}){$j=1$ to $k$}{
        \If(\tcp*[f]{at most one nonrecursive eqn}){$\partial\bs{H}_{i_{j-1}+1:i_j}/{\partial\bs{Y}_{i_{j-1}+1:i_j}}=0$}
        {$S_{i_j}:=\textsf{NR-Localbehavior}(\mathcal{Y}=\cal{H}_{i_j}(\mathcal{Z},\bc Y_{1:i_{j-1}}),\sigma, \bs S_{1:i_{j-1}}, p)$}
        \lElse{
        $\bs S_{i_{j-1}+1:i_{j}}:=\textsf{Irr-Localbehavior}(\bc{Y}=\bc{H}_{i_{j-1}+1:i_j}(\cal{Z},\bc Y_{1:i_{j-1}}, \bc{Y}), \sigma, \bs S_{1:i_{j-1}}, p)$
        }
    }
    \Return{$\bs S_{1:m}$}
\end{algorithm}

\subsection{Example}\label{sec:ex-singular-behavior}
For the last time, we consider the system of \cref{eq:colored_forest} given in \cref{ex:dependency_graph}. Again, we proceed along the strongly connected components.
Recall from \cref{ex:radius_sys1} that the radius of convergence is $\rho:=\rho_T\in[0.1703916,0.1703917]$ and it arises from the irreducible component defining $(T_r,T_b,T_g)$.
\paragraph{First and Second Components $\cal G,(\cal B, \cal R)$.}
Both components match the condition at ligne~\ref{algline:nonsing-nonlin} in Algorithm~\textsf{\ref{algo:localbehavior-irred}}. As $\rho_T$ is smaller than the radius of convergence of $G,B,R$, they have a Taylor expansion there. 

In this example, the equations defining the generating functions $G(z),B(z),R(z)$ are polynomials of small degrees, so that these expansions can also be obtained from 
closed-form expressions:
\[G(z)=\frac{1-\sqrt{1-4z}}{2},\quad
B(z)=\frac12-\frac12\sqrt{\frac{1-4z-z^3}{1-z^3}},\quad
R(z)=\frac{1+z^3}2-\frac{\sqrt{1-z^3}\sqrt{1-4z-z^3}}2.
\]

In general however, one would only have access to numerical approximations of the values of $G(\rho_T),B(\rho_T),R(\rho_T)$ with arbitrary precision obtained by Newton's iteration from their defining systems in  \cref{ex:radius_sys1}. From there, the expansions are obtained using either undetermined coefficients or Newton iteration over power series, giving
\begin{align*}
    G&=0.217851 - 0.301953u + 0.161573u^2 - 0.172913u^3 + O(u^4),\\
    B&=0.219356 - 0.309632u + 0.179975u^2 - 0.204836u^3 + O(u^4),\\
    R&=0.223218 - 0.319686u + 0.186075u^2 - 0.200419u^3 + O(u^4),
\end{align*}
where $u=1-z/\rho$. We display only the first 6~digits of each coefficients and expansions to arbitrary order are easily computed.


\paragraph{Third Component $(\cal T_r,\cal T_b,\cal T_g)$.}
This is the singular nonlinear case at line~\ref{algline:sing-nonlin} in Algorithm \textsf{\ref{algo:localbehavior-irred}}, from where the dominant singularity occurs. By \cref{prop:DLW}, the expansion is of the form
\[\begin{pmatrix}T_r(z)\\ T_b(z)\\ T_g(z)\end{pmatrix}=
\begin{pmatrix}T_r(\rho_T)\\ T_b(\rho_T)\\ T_g(\rho_T)\end{pmatrix}-\begin{pmatrix}C_r\\ C_b\\ C_g\end{pmatrix}\sqrt{1-z/\rho_T}+\dotsb
\]
and the values of $T_r(\rho_T),T_b(\rho_T),T_g(\rho_T)$ have been computed by Newton iteration in \cref{ex:radius_sys1}.

Again, since the equations are polynomial and of low degrees, the system can be solved symbolically in that case, leading to
\begin{multline*}
T_g=1-\frac B{T_b},\quad T_r=z+\frac{R}{1-T_b},\\
T_b=\frac{(1-z)(1+B)-G-R-\sqrt{((1-z)(1+B)-G-R)^2-4B(1-z-G)(1-z-R)}}{2(1-z-G)}.
\end{multline*}
The computation performed by the algorithm does not exploit the possible existence of such closed forms. Instead, the expansions are computed in two stages. First, $u$ is replaced by $v^2$ where $v$ will be used to perform expansions in $\sqrt{1-z/\rho_T}$. Replacing $B,G,R$ by their expansions (now in $v$) and $T_r$ by $T_r(\rho)+\delta_{T_r}v$ in the last two equations of the system \cref{eq:TrTbTg} produces a system with unknowns $T_b,T_g$ and variables $v,\delta_{T_r}$ whose Jacobian matrix is invertible at $(0,0)$ (see the proof of \cref{prop:DLW}). It can thus be solved in power series in $v$ with coefficients that are polynomials in $\delta_{T_r}$, e.g., by Newton iteration over power series, giving
\begin{align}\label{eq:TrTb}
T_b&=0.466295+ 1.276068\, \delta_{T_r}  v +\left(- 1.385784+ 6.594118\, \delta_{T_r}^{2}\right) v^{2}+\dotsb\\
T_g&=0.529578+ 1.287362\, \delta  v +\left(- 0.734024+ 3.129476\, \delta^{2}\right) v^{2}+\dotsb,
\end{align}
where again we display only the first 6 digits of the coefficients and expansions to arbitrary order are easily computed.

Next, these expansions are substituted into the first equation of \cref{eq:TrTbTg}. Expanding in $v$ gives
\[\epsilon_1+ \epsilon_2 \delta_{T_r}  v +\left( 1.855365- 7.558492\, \delta_{T_r}^{2}\right) v^{2}+\left( 12.543550\, \delta_{T_r} - 57.130799\, \delta_{T_r}^{3}\right) v^{3}+\dotsb.\]
The coefficients $\epsilon_1$ and $\epsilon_2$ are numerically small and are known to be exactly~0. They can thus be discarded, the resulting expansion can be divided by $v^2$ and can then be solved for $\delta_{T_r}$, e.g., by Newton iteration for power series, once the correct root of the leading term (the negative one) is selected. This gives
\[\delta_{T_r}=- 0.495447- 0.097917\, v - 0.003464\, v^{2}+ 0.021532\, v^{3}+\dotsb\]
whence
\[T_r=0.588633- 0.495447\,{u}^{1/2}- 0.097917\, u - 0.003464\, u^{{3}/{2}}+\dotsb\]
which can be injected into \cref{eq:TrTb}, giving
\begin{align}
T_b&=0.466295- 0.632224\,{u}^{1/2}+ 0.107910\, u + 0.232963 u^{{3}/{2}}+\dotsb,\\
T_g&=0.529578- 0.637819\,{u}^{1/2}- 0.091894\, u + 0.258036\, u^{{3}/{2}}+\dotsb.
\end{align}

\paragraph{Last Component $\cal F$.} The generating function of $\cal F$ is $F(z)=\exp(T_r(z))$. 
The expansion is obtain by composition of that of $\exp(z)$ with that of $T_r$, giving
\begin{equation}\label{eq:sing_exp_F}
F=1.801525- 0.892560\, {u}^{1/2}+ 0.044708\, u + 0.044641\, u^{{3}/{2}}+\dotsb.
\end{equation}



\section{Asymptotic Behavior}\label{sec:transfer}
Once the dominant singularities of the generating functions have been computed, together with the local behavior of the generating functions at these dominant singularities, the asymptotic behavior of their coefficients is determined by the method of \emph{singularity analysis}, described for instance in Chapter~6 of the book by Flajolet and Sedgewick~\cite{FlajoletSedgewick2009}. 
A version of these results that is sufficient for this article is the following.
\begin{proposition}(Direct corollary of \cite[Thm.VI.5]{FlajoletSedgewick2009})\label{prop:FlSeThm.VI.5}
If $f(z)$ is analytic in $|z|<\rho$, has finitely many singularities $\sigma_1,\dots,\sigma_p$ on~$|z|=\rho$, admits an analytic continuation in a disc $|z|<\rho'$ with $\rho'>\rho$ slit along the segments $(\rho,\rho')\exp(i\arg\sigma_j)$ and behaves locally as
\begin{multline*}
f(z)=\sum_{j=1}^{m_i}{\sum_{k=0}^\ell c_{i,j,k}(1-z/\sigma_i)^{e_{i,j}}\ln^k\!\left(\frac1{1-z/\sigma_i}\right)}\\
+O\!\left((1-z/\sigma_i)^{f_i}\ln^{b_i}\!\left(\frac1{1-z/\sigma_i}\right)\right),\quad z\rightarrow\sigma_i,\quad i=1,\dots,p,
\end{multline*}
then the asymptotic behavior of~$[z^n]f(z)$ is given by
\begin{equation}\label{eq:asympt-coeffs}
[z^n]f(z)=\sum_{i=1}^p\left(\sum_{j=1}^{m_i}\sum_{k=0}^\ell c_{i,j,k}\sigma_i^{-n}C_n({e_{i,j}},k)\right)+O(\rho^{-n}n^{-\min\Re{f_i}-1}(\ln n)^{\max b_i}),\quad n\rightarrow\infty
\end{equation}
where
\[C_n(\alpha,k)=\begin{cases}0,\quad\text{if $\alpha\in\mathbb N$ and $k=0$},\\
(-1)^k\frac{\partial^k}{\partial\alpha^k}\frac{\Gamma(n-\alpha)}{\Gamma(-\alpha)\Gamma(n+1)}
     \sim\begin{cases}
     \frac{n^{-\alpha-1}\ln^kn}{\Gamma(-\alpha)},\quad n\rightarrow\infty&\quad{\text{if $\alpha\in\mathbb C\setminus\mathbb N$}},\\
     (-1)^{\alpha+k}k\alpha!n^{-\alpha-1}\ln^{k-1}n,\quad n\rightarrow\infty&\quad{\text{if $\alpha\in\mathbb N$.}}
     \end{cases}
\end{cases}\]
\end{proposition}
\begin{example}Let $f_n$ be the number of forests with $n$ nodes from our running \cref{ex:colored_forest}. 
Its asymptotic behavior can now be deduced directly from the expansion of \cref{sec:singular-behavior}. Indeed, from the expansion of a generating function at its dominant singularity~$\rho$ of the form
\[c_0+c_1u^{1/2}+c_2u+c_3u^{3/2}+\dotsb,\]
with $u={1-z/\rho}$, \cref{prop:FlSeThm.VI.5} shows that the $n$th coefficient behaves asymptotically like
\[\frac{\rho^{-n}}{n!}\left(c_1\frac{\Gamma(n-1/2)}{\Gamma(-1/2)}+c_3\frac{\Gamma(n-3/2)}{\Gamma(-3/2)}+c_5\frac{\Gamma(n-5/2)}{\Gamma(-5/2)}+O(\Gamma(n-7/2))\right),\quad n\rightarrow\infty.\]
Expanding by Stirling's formula shows that this is asymptotically
\[\rho^{-n}n^{-3/2}\left(-\frac{c_{1}}{2 \sqrt{\pi}}-\frac{3 \left(c_{1}-4 c_{3}\right)}{16 \sqrt{\pi}\, n}-\frac{5 \left(5 c_{1}-72 c_{3}+96 c_{5}\right)}{256 \sqrt{\pi}\, n^{2}}+O\! \left(\frac{1}{n^{3}}\right)\right),\quad n\rightarrow\infty.\]
With the values of \cref{eq:sing_exp_F}, this leads to 
\begin{equation}\label{eq:final_asympt}
\frac{f_n}{n!}=\rho_T^{-n}n^{-3/2}\left(0.251787 + \frac{0.113309}n+ \frac{0.112583}{n^2}+\dotsb\right),\quad n\rightarrow\infty,
\end{equation}
with $\rho_T\in [0.1703916,0.1703917]$ root of the polynomial from \cref{pol-deg-21}. The coefficients, and $\rho_T$, can be computed to arbitrary precision efficiently. 

Already with $n$ as small as~10, one has the exact value $f_{10}/10!=5767537729/14175\approx 406881.$, while the first term of the asymptotic expansion above at~$n=10$ is approximately~$385966.$, giving a relative error of $5.4\%$. Using one more term gives~$403335.$ and a relative error of $0.9\%$. Using all three gives~$405060$ and a relative error of $0.45\%$. 
These errors decrease to~$4.5\,10^{-3},4.6\,10^{-5},1.6\,10^{-6}$ for~$n=100$.
\end{example}

\begin{example}\label{ex:ex27}\Cref{ex:imaginary-exponent} leads to the generating function
\[Y(z)=\exp\!\left(z\ln\!\left(\frac1{1-z^4}\right)\right)=1+z^5+\frac12z^9+\frac12z^{10}+
\frac13z^{13}+\frac12z^{14}+\dotsb=:\sum_{n\ge0}y_nz^n\]
that has 4 dominant singularities at~$\{\pm1,\pm i\}$. Its local behavior at these singularities is 
\begin{alignat*}{3}
f(z)&=\frac{1/4}{1-z}-\frac14\ln\!\left(\frac1{1-z}\right)
+\frac38+\frac{\ln 2}2
+O\!\left((1-z)\ln^2\!\left(\frac1{1-z}\right)\right),&\qquad&z\rightarrow 1,\\
f(z)&=4(1+z)+O\!\left((1+z)^2\ln\!\left(\frac1{1+z}\right)\right),&\qquad&z\rightarrow-1,\\
f(z)&=4^{-i}(1-z/i)^{-i}+O\!\left((1-z/i)^{-i+1}\ln\!\left(\frac1{1-z/i}\right)\right),&\qquad&z\rightarrow i
\end{alignat*}
and the conjugate of the last one as~$z\rightarrow-i$. These expansions contribute respectively to the asymptotic behavior of the coefficient of~$z^n$ in~$Y$ for
\[
\frac14-\frac1{4n}+O\!\left(\frac{\ln^2n}{n^2}\right),\quad
O\!\left(\frac{\ln n}{n^3}\right),\quad
4^{-i}i^{-n}\frac{\Gamma(n+i)}{\Gamma(n+1)\Gamma(i)}+O\!\left(\frac{\ln n}{n^2}\right),
\]
and the conjugate of the last one.
Adding these contributions, expanding them and beautifying the result with the value of $|\Gamma(i)|$~\cite[Eq.~5.4.3]{dlmf} gives the asymptotic behavior
\[y_n=\frac14+\left(2\sqrt{\frac{\sinh(\pi)}\pi}{\cos\left(\frac{\pi n}2-\ln(n)+\phi\right)-\frac1{4}}\right)\frac1n+O\!\left(\frac{\ln^2n}{n^2}\right),\quad n\rightarrow\infty,\]
with $\phi=2\ln 2+\arg\Gamma(i)\approx-0.486$.
Full asymptotic expansions can be derived the same way, preferably with the help of a computer algebra system.
\end{example}

That \cref{prop:FlSeThm.VI.5} applies in our context follows from the fact that if $\bs Y=\bs H(z,\bs Y)$ is a zero-free well-founded system, then the generating function solution~$\bs Y$ is analytic in a disk of positive radius~$\rho$ (by \cref{prop:well-founded-analytic}) and moreover, each of its irreducible components has solutions with finitely many dominant singularities by \cref{thm:BBY72i,thm:period-linear}. In summary, the following has been proved.
\begin{theorem}\label{thm:everything}
Let $\bs Y$ be a constructible generating function. For all coordinates of $\bs Y$ that have an algebraic-logarithmic dominant singularity, an asymptotic expansion in the form of \cref{eq:asympt-coeffs} can be computed given a \texttt{Constant Oracle} as in \cref{sec:constant-oracle}.
\end{theorem}
Again, note that when the constructive system contains neither \Cyc\ nor \Set, then the oracle is computable, while this computability is dependent on the truth of Schanuel's conjecture otherwise.

\section{Conclusion}\label{sec:conclusion}
The automation of analytic combinatorics was started many years ago \cite{FlajoletSalvyZimmermann1991}. An implementation was developed in Maple and many examples were treated automatically, including applications to the analysis of algorithms~\cite{FlajoletSalvyZimmermann1989}, which was an important motivation for the development of analytic combinatorics. The main limitation of that implementation was that it could only deal with a restricted set of systems of equations. Either the systems were completely nonrecursive, or the recursive equations were sufficiently simple to be solvable in closed form. The main contribution of the present article is to show how these limitations can be removed. 
We conclude by reviewing other parts of analytic combinatorics where more automation work is still needed.

\newcommand{\ind}{\mathfrak{S}}
\paragraph{Cardinality constraints}
For simplicity, we have focused on the species $\Seq,\Set,\Cyc$ in their versions without cardinality constraints. A large part of the theory naturally  extends to indexed versions $\Seq_{\ind},\Set_\ind,\Cyc_\ind$ with $\ind\subset\mathbb N$ a restriction on the cardinality of the sequence, set or cycle. The general case is recovered with $\ind=\mathbb N$ for $\Seq_\ind$ and $\Set_\ind$ and with $\ind=\mathbb N\setminus\{0\}$ for $\Cyc$. Of course, for effectivity reasons, the set~$\ind$ cannot be arbitrary. For the uses in~\cite{FlajoletSedgewick2009}, $\ind$~is either a finite set, or the union of a finite set and a finite number of arithmetic progressions. This is easily accommodated by the algorithms presented here. For instance, when $\mathcal H_\ind$ is one of $\Seq_\ind(\mathcal Y_j),\Set_\ind(\mathcal Y_j),\Cyc_\ind(\mathcal Y_j)$, Algorithm \textsf{\ref{algo:wellfoundedandleadingterm}} would set~$w_i$ to $(c_\ell v_{j,1}^\ell,\ell v_{j,2})$ with $\ell=\min \ind$ and $c_\ell$ one of~$1,1/\ell!,1/\ell$ depending on the case and provided~$v_{j,2}\neq0$; it would use ($\sum_{\ell\in \ind} c_\ell v_{j_1}^\ell,0)$ if $v_{j_2}=0$ and $\ind$ is finite; \textsf{Fail} otherwise. Similarly, Algorithm \textsf{\ref{algo:nr-radius}} returns~$r_U$ for $\Set_\ind$ and arbitrary~$\ind$ and for $\Cyc_\ind$ and $\Seq_\ind$ with finite~$\ind$. It uses its other branch in the other cases. Similar modifications can be made to Algorithms \textsf{\ref{algo:periods}} and \textsf{\ref{algo:dominant-nr}}. The adjustments to \textsf{\ref{algo:localbehavior-nonrecursive}} are equally straightforward, but slightly more technical.

\paragraph{Superpolynomial growth}
We have focused on algebraic-logarithmic singularities and have not discussed what should be done with generating functions that have a superpolynomial behavior. For generating functions whose growth is at least exponential at their dominant singularity or that are entire and non-polynomial, one can use the saddle-point method~\cite[Ch.~VIII]{FlajoletSedgewick2009}. This can be automated to some extent, thanks to Hayman's class of functions, to which membership is easy to check and for which a first order asymptotic expansion is easily computed~\cite{Hayman1956}. A full asymptotic expansion can be obtained for a more restricted class studied by Harris and Schoenfeld~\cite{HarrisSchoenfeld1968,OdlyzkoRichmond1985}. Wyman~\cite{Wyman1959} gives a very nice exposition that is more general, but more difficult to apply in practice. Already in the simple setting of Hayman's class, from the computer algebra point of view, a difficulty lies in the handling of suitable asymptotic scales, see for instance~\cite[Sec.~3]{SalvyShackell1999}.

There are also generating functions whose growth lies between the algebraic-logarithmic ones and the exponential ones. In view of \cref{thm:alg-log-real} and its proof, these are functions that behave at their dominant singularity (normalized here to be~1) like
\[\exp\!\left(\ln^k\frac1{1-z}\right)(1-z)^\alpha,\quad k\in\mathbb N\setminus\{0,1\},\quad \alpha\le 0,\]
or more complicated variants with a polynomial in $\ln$ in the exponential and extra logarithmic factors at the end. We are not aware of results on the asymptotic behavior of their coefficients, which would complete a classification of the possible asymptotic behaviors of the coefficients of the generating functions of constructible species.

\paragraph{Ordinary generating functions}
Ordinary generating functions are used for \emph{unlabeled} enumeration. The dictionary for the translation of a combinatorial specification into ordinary generating functions is very similar to that of \cref{tab:sum_esp_sg}, except for $\Cyc(\mathcal G)$ and $\Set(\mathcal G)$ that translate into
\[\sum_{k>0}\frac{\phi(k)}k\log\frac1{1-G(z^k)},\quad\exp\!\left(\sum_{k>0}\frac{G(z^k)}k\right),\]
where~$\phi$ is Euler's totient function~\cite[Thm.~I.1]{FlajoletSedgewick2009}. Flajolet and Sedgewick also have another construction \textsf{PowerSet} that is not a species, see the discussion in~\cite[\S10.2]{PivoteauSalvySoria2012}.

This translation shows that differences between ordinary and exponential generating functions only matter if $\Set$ or $\Cyc$ are used. Newton's iteration can be used to obtain numerical values of the generating functions~\cite{PivoteauSalvySoria2012}, but the general analysis of the location of singularities is more involved. In many cases, there are infinitely many singularities on the circle of convergence, preventing singularity analysis to apply. However, this is not the case when the radius of convergence is less than~1. This can sometimes be detected by comparison with the radius of convergence of a related species. A typical example is the class of P\'olya trees whose generating function satisfies
\[P(z)=z\exp\!\left(\frac{P(z)}1+\frac{P(z^2)}2+\dotsb\right).\]
Its radius of convergence is clearly smaller than~$e^{-1}$, the radius of convergence of the generating function of the Cayley trees, solution to~$T(z)=z\exp(T(z))$. Then one can compute the radius efficiently (see~\cite[\S15.1]{Odlyzko1995},\cite[pp.~475--477]{FlajoletSedgewick2009},\cite[\S9.4.2]{PivoteauSalvySoria2012}) and ultimately, a full asymptotic expansion of the coefficients. In this setting, we do not have computability results like \cref{th:radius-is-computable} but it would be interesting to investigate the automation of some of the techniques, along the lines of~\cite{Genitrini2016}.

\paragraph{Multivariate generating functions}
The combinatorial part of analytic combinatorics (the symbolic method) can easily deal with several sorts of atoms, which then translate into equations over multivariate generating functions.

As a typical example, in order to study the expected number of trees in the forests~$\mathcal F$ from \cref{ex:colored_forest}, one would replace the first equation by $\mathcal F=\Set(\mathcal U\times\mathcal T_r)$ for a new sort of atom~$\mathcal U$. This equation translates into the generating function equation $F=\exp(uT_r(z))$ whose coefficient of $u^kz^n$ is the number of forests with $n$ nodes using $k$ trees~$\mathcal T_r$, divided by $n!$. Differentiating with respect to~$u$ and setting~$u=1$ shows that the expected number of trees is given by $[z^n]T_rF/[z^n]F$. With the expansions computed in \cref{sec:ex-singular-behavior}, it is easy to compute the behavior of $T_rF$ as $z\rightarrow\rho_T$ and from there obtain the asymptotic behavior of the numerator, whence, after division by \cref{eq:final_asympt}, the asymptotic expansion of the expected number of trees in a forest $\mathcal F$ drawn uniformly at random among all forests $\mathcal F$ with $n$ nodes:
\[\mathbb E[\#\mathcal T_r\text{ in }\mathcal F]=1.588633+\frac{0.0241419}{n}+\frac{0.117434}{n^{2}}+\dotsb,\qquad n\rightarrow\infty.\]
Again, these coefficients are easily computed to arbitrary precision and more terms of the expansion are easily found.

Many other parameters can be studied in a similar way, see \cite[Ch.~III,IX]{FlajoletSedgewick2009}. Analytic combinatorics also allows to predict limiting distributions in many cases. Not much has been done in the automation of these aspects of the theory. A notable exception is the extension of Hayman's class to the bivariate setting~\cite{DrmotaGittenbergerKlausner2005}.

A direction that has received a lot of attention in recent years is the field of `Analytic Combinatorics in Several Variables'. There, the aim is to extract asymptotic information on the coefficient sequence of multivariate generating functions along a diagonal. Already the case of rational generating functions can become involved. In some cases, the approach can be automated~\cite{MelczerSalvy2021,HacklLuoMelczerSeloverWong2023,HacklLuoMelczerSchost2025}.\footnote{\url{https://github.com/ACSVMath/sage_acsv}} The recent new edition of the reference book on that topic~\cite{PemantleWilsonMelczer2024} will be a source of ideas and questions for computer algebraists with an interest in this area.

\paragraph{Implementation}
Maple worksheets illustrating some of the examples are available with this article.\anonymousvariant{}{\footnote{The Editflow system used to submit this article does not seem to allow anonymous submission of supplementary material.}}
Still, a complete implementation of the algorithms presented here should be attempted. 
For a long time, many of the necessary routines were available in the Algolib\footnote{\url{https://algo.inria.fr/libraries/}} Maple library of the Algorithms Project at Inria, but these have not been maintained and do not run well on modern versions of Maple. The Maple package NewtonGF\footnote{\url{https://perso.ens-lyon.fr/bruno.salvy/software/the-newtongf-package/}} provides the tools required for the combinatorial part and for generating function expansions~\cite{PivoteauSalvySoria2012}. It could serve as a basis to develop the analytic part. We are not aware of implementations of the asymptotic part of analytic combinatorics in other computer algebra systems. 


\clearpage
\appendix
\section{Proof of Theorem \ref{thm:gl-implicit-leading} : Characterization by Leading Terms}\label{appendix:prooftheorem1.7}
\glimplicitleading*
Leading terms can be computed inductively, thanks to the following.
\begin{lemma}\label{lemma:leading-leading}
For any species $\bc H,\bc G$ such that the composition $\bc H(\bc Z, \bc G)$ is well defined, the leading terms satisfy the equality
$\LT(\bc H(\bc Z, \bc G))=\LT(\bc H(\bc Z,\LT(\bc G)))$.
\end{lemma}
\begin{proof} 
This is a direct consequence of the definition of composition of species
\[\mathcal H(\mathcal G_1,\dots,\mathcal G_k)[U]=\sum_{\substack{\pi\text{ partition of }U\\ \pi_1+\dots+\pi_k=\pi}}\mathcal H[\pi_1,\dots,\pi_k]\times\prod_{p_1\in\pi_1}\mathcal G_1[p_1]\times\dots\times\prod_{p_k\in\pi_k}\mathcal G_k[p_k].\]
The smallest $\mathcal H(\mathcal G_1,\dots,\mathcal G_k)$-structures are therefore obtained from the smallest $\mathcal G_i$ structures. This applies to each of the coordinates~$\mathcal H$ of~$\bc H$.
\end{proof}

\begin{lemma}\label{lemma:leading-U}
Let $\bc H_{1:m}(\bc Z,\bc Y)$ be a vector of species such that the iteration \eqref{eq:ite-def-wf} is well defined for $n=0,\dots,m$ and let $\bc U^{[k]}$ be defined as in~\cref{thm:gl-implicit}, then $\bc U^{[k]}=\LT(\bc Y^{[k]})(\bs 0)$ for $k=0,\dots,m$.
\end{lemma}
\begin{proof} The proof is by induction. First, $\bc Y^{[0]}=\bs 0=\LT(\bs 0)=\bc U^{[0]}$. Next, 
if the property holds for~$k$, then
\begin{multline*}\LT(\bc Y^{[k+1]})(\bs 0)=
\LT(\bc H(\bc Z,\bc Y^{[k]}))(\bs 0)=
\LT(\bc H(\bc Z,\LT(\bc Y^{[k]})))(\bs 0)\\ =
\bc H(\bs 0,\LT(\bc Y^{[k]})(\bs 0))=
\bc H(\bs 0,\bc U^{[k]})=
\bc U^{[k+1]},\end{multline*}
where the second equality uses \cref{lemma:leading-leading} and the third one uses the fact that for any species~$\bc F$, $(\LT(\bc F))(\bs 0)=\bc F(\bs 0)$.
\end{proof}

\begin{lemma}\label{lemma:composition-leading-terms}
The composition $\bc H(\bc Z,\bc Y)$ is well defined if and only if the composition with the leading term, $\bc H(\bc Z,\LT(\bc Y))$, is well defined.
\end{lemma}
\begin{proof}
As $\LT(\bc Y)\subset\bc Y$, the composition $\bc H(\bc Z,\LT(\bc Y))$ is well defined when $\bc H(\bc Z,\bc Y)$ is.

Conversely, the composition with $\bc Y$ is well defined if and only if the composition with $\bc Y(0)$ is. As $\bc Y(0) = \LT(\bc Y(0))$, if $\bc H(\bc Z,\LT(\bc Y))$ is well defined, then so is $\bc H(\bc Z,\bc Y)$.
\end{proof}

\begin{lemma}\label{lemma:leading-terms}
Let $\bc Y^{[0]}$ and $\ell\in\mathbb N$ be such that the iteration $\bc Y^{[n+1]}=\bc H(\bc Z,\bc Y^{[n]})$ starting with~$\bc Y^{[0]}$ is well defined for $k=1,\dots,\ell$ and $\LT(\bc Y^{[\ell+1]})=\LT(\bc Y^{[\ell]})$. Then the iteration is well defined for all~$n\ge1$ and the value of~$\LT(\bc Y^{[n]})$ is constant for~$n\ge\ell$. If $\bc C$ is this value, then the minimal~$\ell$ such that $\LT(\bc Y^{[\ell]})=\bc C$ is at most the number of coordinates of~$\bc C$ different from those of~$\LT(\bc Y^{[0]})$.
\end{lemma}
\begin{proof}
As the iteration is well defined up to~$n=\ell$, by \cref{lemma:composition-leading-terms}, the composition~$\bc H(\bc Z,\bc  Y^{[n]})$ is well defined for all~$n$ if the leading terms of~$\bc Y^{[n]}$ is~$\LT(\bc Y^{[\ell]})$. By induction using \cref{lemma:leading-leading}, leading terms of  $\bc Y^{[n+1]}$ can only be built from leading therms of $\bc Y^{[n]}$, thus $\LT(\bc Y^{[n]})=\LT(\bc Y^{[\ell]})=\bc C$ for any $n\geq\ell$. So the iteration $\bc Y^{[n+1]}=\bc H(\bc Z,\bc Y^{[n]})$ is well defined for all~$n$.

Without loss of generality, let $\ell$ be minimal such that $\LT(\bc Y^{[\ell+1]})=\LT(\bc Y^{[\ell]})$ and let $J$ be the set of indices where $\LT(\bc Y^{[\ell]})\neq\LT(\bc Y^{[\ell-1]})$. By definition of~$\ell$, for $j\in J$, the coordinates $\mathcal C_i$ for $i\neq j$ do not depend on $\mathcal C_j$. 
This implies that the vector of coefficients~$(\mathcal C_i | i\not\in J)$ is the same as that of the system  $\bc Y=\bc H(\bc Z,\bc Y)|_{\cal Y_j=0, j\in J}$ in $|J|$ variables less. This new system has $\ell'=\ell-1$, as coefficients change at each iteration, and $|J|>0$ fewer variables. This construction can be repeated less times than the number of coordinates of~$\bc C$ different from those of~$\LT(\bc Y^{[0]})$, hence the result. 
\end{proof}

The next three properties follow directly from the definitions.
\begin{lemma}\label{lemma:leading-constants}
For any species $\cal F$ and $\cal G$, 
\begin{enumerate}
    \item $\cal F(0) = \LT(\cal F(0)) = \LT(\cal F)(0)$;
    \item if $\cal F(0)\neq 0$ then $\LT(\cal F) = \cal F(0)$;
    \item if $\cal F(0) = \cal F$ and $\cal F \subset \cal G$ then $\cal F \subset \LT(\cal G)$.
\end{enumerate}
\end{lemma}

\begin{lemma}\label{lemma:truncated-system}
If the iteration \eqref{eq:ite-def-wf} is well defined for $n=1,\dots,m$ and $\LT(\bc Y^{[m+1]})=\LT(\bc Y^{[m]})$, then let $\bc C=\LT(\bc Y^{[m]})$ and $\bc U=\bc C(\bs 0)$. Up to renumbering the coordinates, let $p$ be such that $\mathcal U_{1:p}=\bs 0$ and the coordinates $\mathcal U_i\neq0$ for $i=p+1,\dots,m$. Write $\bct Y=\bc Y_{1:p}$ and consider the system~$\bct Y=\bc H_{1:p}(\bc Z,\bct Y,\mathcal U_{p+1:m})$. Then, 
\begin{enumerate}
    \item the iteration $\bct Y^{[n+1]}=\bct H(\bc Z,\bct Y^{[n]})$ starting with~$\bct Y^{[0]}=\bs 0$ is well defined for $n=1,\dots,p$;
    \item it satisfies $\LT(\bct Y^{[p]})=\LT(\bct Y^{[p+1]})=\bc C_{1:p}$;
    \item the Jacobian matrix $\bpartial \bc H/\bpartial \bc Y(\bs 0,\bc U)$ is nilpotent if and only if 
    $\bpartial\bct H/\bpartial\bct Y(\bs 0,\bs 0)$ is nilpotent.
\end{enumerate}
\end{lemma}
\begin{proof}
Note that when $p=0$, the first two points are irrelevant and the third one just states that the matrix is  unconditionally nilpotent.

1. From \cref{lemma:leading-terms}, the composition $\bc H(\bc Z, \bc Y^{[k]})$ is well defined for all $k\geq 0$ and, from $(\bct Y^{[0]},\bc U_{p+1:m})\subset \bc Y^{[m]}$, it follows by induction that $(\bct Y^{[n]},\bc U_{p+1:m})\subset \bc Y^{[n+m]}$, thus $\bct H(\bc Z,\bct Y^{[n]})$ is well defined for all $n$.

2. Let $(\bch Y^{[n]})_n$ be the sequence defined by $\bch Y^{[n+1]}=\bc H(\bc Z,\bch Y^{[n]})$ with initial value $\bch Y^{[0]}=\bc U$. By induction, one has $\bc Y^{[n]}\subset \bch Y^{[n]}\subset \bc Y^{[n+m]}$ and from the convergence of $\LT(\bc Y^{[n]})$ to~$\bc C$, it follows that $\LT(\bch Y^{[n]})$ also converges to~$\bc C$. 

Next, we prove that $\LT(\bch Y_{p+1:m}^{[n]})=\bc U_{p+1:m}$. 
Inclusion is preserved by composition with~$\bc H$, thus $\bc U_{p+1:m} \subset \bch Y_{p+1:m}^{[n]}$ and all the coordinates of $\bch Y_{p+1:m}^{[n]}$ are non empty. And by definition of~$\bc U$ and the third point of \cref{lemma:leading-constants}, $\bc U_{p+1:m} \subset \LT(\bch Y_{p+1:m}^{[n]}$).   On the other hand,
\[
\LT(\bch Y_{p+1:m}^{[n]}) = \bch Y_{p+1:m}^{[n]}(0) \subset \bc Y_{p+1:m}^{[n+m]}(0) =\bc U_{p+1:m}.
\]
Both equalities come from  \cref{lemma:leading-constants} and the inclusion follows from $\bch Y^{[n]}\subset \bc Y^{[n+m]}$.

Then, by induction, $\LT(\bch Y^{[n]})=(\LT(\bct Y^{[n]}),\bc U_{p+1:m})$. Indeed, the right-hand side is the previous point. As for the left-hand side, it is true for $n=0$ by definition of~$\bc U$ and if the property holds for~$n$, then it also holds for~$n+1$, using \cref{lemma:leading-leading}.

It follows that $\LT(\bct Y^{[n]})$ converges to~$\bc C_{1:p}$. \Cref{lemma:leading-terms} then proves that the convergence requires at most~$p$ iterations.

3. By~\cref{lemma:leading-U}, $\bc U=\bc U^{[m]}$, with $\bc U^{[k]}$ defined as in~\cref{thm:gl-implicit} for $k=0,\dots,m$. 
Let $k\le m$ be the maximal index where $\bc B=\bc U^{[k]}-\bc U^{[k-1]}\neq\bs 0$ and let $\bc A=\bc U^{[k-1]}$. From 
$\bc H(0,\bc A+\bc B)=\bc H(0,\bc A)$,
it follows that $\partial\bc H/\partial\cal Y_i(\bs 0,\bc Y)=0$ for all~$i$ such that~$\cal B_i\neq0$. Without loss of generality, let these coordinates be the last $m-q$ ones ($1\le m-q\le m-p$ by definition of~$k$ and the hypothesis on~$\bc B$) and consider the new system formed with 
$\bc{K}:=\bc H_{1:q}(\bc Z,\bc Y_{1:q},\bc{U}_{q+1:m})$. Since $\bc H(\bs 0,\bc Y)$ does not depend on the last $m-q$ coordinates, the Jacobian matrix ${\bpartial \bc H}/{\bpartial\bc Y}(\bs 0,\bc U^{[m]})$ has a top-left $q\times q$ block that is the Jacobian matrix of~$\bc{K}$ at~$(\bs 0,\bc U^{[m]}_{1:q})$, while its next $m-q$ columns are filled with~0. Thus these two matrices are either both nilpotent or both nonnilpotent. Repeating this argument with the new system until~$\bct H$ is reached or no matrix is left (when $p=0$) proves the result.
\end{proof}

The proof of \cref{thm:gl-implicit-leading} can now be detailed.
\begin{proof}
Assume that the system is well founded. Then the iteration \eqref{eq:ite-def-wf} is well defined for all~$n$ and converges to a species~$\bc S$. The sequence $\LT(\bc Y^{[n]})$ converges in a finite number of steps: the sequence of valuations is a decreasing sequence of $m$-tuples of nonnegative integers or infinity; once the valuations have converged, the sequence of leading terms is increasing by inclusion and contained in the polynomial species~$\LT(\bc S)$. Thus there exists~$\ell$ such that $\LT(\bc Y^{[\ell]})=\LT(\bc Y^{[\ell+1]})=\LT(\bc S)$. By \cref{lemma:leading-terms} one can take $\ell\le m$. 

Conversely, assume that the leading terms of~$\bc Y^{[m]}$ and $\bc Y^{[m+1]}$ are equal. Then \cref{lemma:leading-terms} shows that the iteration is well defined for all~$k\ge m$ with $\LT(\bc Y^{[k]})=\LT(\bc Y^{[m]})=:\bc C$. This shows that if the system is well founded, then its solution $\bc S$ satisfies $\LT(\bc S)=\bc C$. All that remains is to prove it is well founded. From now on, without loss of generality, we assume that none of the coordinates of~$\bc C$ is~0. (Otherwise, the corresponding coordinates of~$\bc Y^{[k]}$ are~0 for all~$k$ and can be removed from the system.)

If $\bc C(\bs 0)=\bs 0$, then by Taylor's formula and induction, for any $q\geq 0 $,
\begin{multline*}
\bc Y^{[m+q]} = \bc H(\bc Z,\bc Y^{[m+q-1]})
\supset \left(\frac{\bpartial \bc H}{\bpartial\bc Y}(\bc Z,\bs{0})\right)\bc Y^{[m+q-1]}
\supset\dotsb \\
\supset \left(\frac{\bpartial \bc H}{\bpartial\bc Y}(\bc{Z},\bs{0})\right)^{\!\!q}\bc Y^{[m]}
\supset \left(\frac{\bpartial \bc H}{\bpartial\bc Y}(\bs{0},\bs{0})\right)^{\!\!q}\bc C.
\end{multline*}
If the matrix $\bpartial\bc H/\bpartial\bc Y(\bs 0,\bs 0)$ is not nilpotent, there exists at least one nonzero structure in $({\bpartial \bc H}/{\bpartial\bc Y}(\bs{0},\bs{0}))^q$, and its size is~$0$.
Thus, there are indices $i_1,i_2,\dots i_{q+1}$ and structures $\beta_{i_1} \in {\partial \cal H_{i_1}}/{\partial\cal Y_{i_2}}(\bs{0},\bs{0})$, $\beta_{i_2} \in {\partial \cal H_{i_2}}/{\partial\cal Y_{i_3}}(\bs{0},\bs{0})$, \dots, $\beta_{i_q} \in {\partial \cal H_{i_q}}/{\partial\cal Y_{i_{q+1}}}(\bs{0},\bs{0})$, such that all the $\beta_i$ are nonzero structures. There can be at most~$m$ such indices that are distinct, thus $i_k=i_\ell$ for some $k,\ell \in [1,q+1]$. 
As none of the coordinates of $\bc C$ is~0, one can choose a structure $\gamma$ in $\cal C_k$, and find infinitely many structures of size $|\gamma|$ of the form $(\beta_{i_k}\cdot \dotsc \cdot\beta_{i_{\ell-1}})^j\cdot \gamma$ with $j\geq 0$ belonging to the sequence $(\cal Y_k^{[n]})_{n\in\N}$. This is in contradiction with the convergence of the sequence $\LT(\bc Y^{[k]})$, showing that the matrix $\bpartial\bc H/\bpartial\bc Y(\bs 0,\bs 0)$ is nilpotent and therefore that the system is well founded when $\bc C(\bs 0)=\bs 0$. 

If $\bc C(\bs 0)\neq\bs 0$ then the hypotheses of \cref{lemma:truncated-system} hold with $p\ge0$. When $p\ge 1$, this lemma defines a new system with $p$ unknowns, $\bct{Y}=\bct H(\bc Z,\bct Y)$. By \cref{lemma:truncated-system}, the iteration defined by the new system is well defined, it satisfies $\LT(\bct Y^{[p]})=\LT(\bct Y^{[p+1]})$ and this is~$\bs 0$ at~$\bs 0$. Thus the previous part of the proof applies, showing that this system is well founded. Extending to $p\ge0$, \cref{lemma:truncated-system} then shows that~$\bc H$, from the original system has a nilpotent Jacobian matrix at~$(\bs 0,\bc C(\bs 0))$. By \cref{lemma:leading-U}, this is equal to~$(\bs 0,\bc U^{[m]}(\bs 0))$, with $\bc U^{[k]}$ the sequence defined in \cref{thm:gl-implicit}. This theorem then shows that the system is well founded, which concludes the proof.
\end{proof}

\section{Proof of Theorem \ref{thm:Lambda}: Dominant Eigenvalue}
\label{appendix:prooftheorem2.3}
\thLambda*

The positivity of $\rho$ follows from \cref{prop:well-founded-analytic}.
The rest of the proof relies on the ideas described by Bell, Burris and Yeats in the analysis of their well-conditioned systems~\cite{BellBurrisYeats2010} (see \cref{table:system_types}).

\subsection{Behavior at the Singularity}
A basic tool in several proofs for irreducible systems is the following.
\begin{lemma}\label{lemma:NoZeroInJac}
Let $\bs Y=\bs H(z,\bs Y)$ be an irreducible system. Then its solution satisfies a well-founded system~$\bs Y=\bs K(z,\bs Y)$ where all entries of the Jacobian matrix $\partial\bs K/\partial\bs Y$ are nonzero and positive.
\end{lemma}
\begin{proof}
Let $m$ be the dimension of the system.
Consider the sequence of systems $\bs Y=\bs H^{[i]}(z,\bs Y)$ with $\bs H^{[i]}$ defined by $\bs H^{[1]}=\bs H$ and $\bs H^{[i+1]}(z,\bs Y)=\bs H(z,\bs H^{[i]}(z,\bs Y))$. All these systems have the same domain of convergence and the same solution and so does the system~$\bs Y=\bs K(z,\bs Y)$, with
\begin{equation}\label{eq:defK}
\bs K(z,\bs Y):=\frac1m(\bs H^{[1]}(z,\bs Y)+\dots+\bs H^{[m]}(z,\bs Y)).
\end{equation}
All the entries of the Jacobian matrix of~$\bs K$ are nonzero and positive and depend on~$z$, since the effect of this sum is to add up all paths of length at most~$m$ in the dependency graph of the original system, which is irreducible.
\end{proof}

The following lemma implies that all coordinates of the solution have a similar behavior at the dominant singularity~$\rho$. It is similar to~\cite[Prop.~3 (iv)]{BellBurrisYeats2010} and parts of \cite{BellBurrisYeats2006}.

\begin{lemma}\label{lemma:finiteornot}
Let $\bs Y=\bs H(z,\bs Y)$ be an irreducible system and~$\rho$ be the radius of convergence of its generating function solution~$\bs Y$. Either all coordinates of $\bs Y(x)$ tend to infinity as $x\rightarrow\rho-$, or all of them have a finite limit. The former case occurs only for $\bs H$ that is linear in~$\bs Y$.
\end{lemma}
\begin{proof}
Let $m$ be the dimension of the system.
Observe first that an irreducible system is well-founded and zero-free, thus~$\partial\bs H/\partial z\neq \bs 0$ since~$\bpartial\bs H/\bpartial \bs Y(z,\bs Y)$ is not nilpotent, while its evaluation at~$(0,\bs U_m)$ with $\bs U_m=\bs Y(0)$ is.

\Cref{lemma:NoZeroInJac} implies that each of the coordinates of $\bs Y$ depends positively on the other ones and thus that if one of them tends to infinity, all the other ones do. 

In the nonlinear case, we follow the argument of~\cite[Prop.~66]{BellBurrisYeats2006}.
Since $\bs H$ is not linear in $\bs Y$, 
there is an index~$i$ such that the $i$th equation in \cref{eq:defK} contains a monomial of the form~$ c_{ij}z^kY_iY_1^{j_1}\dots Y_m^{j_m}$ with at least one $j_\ell\neq0$. Positivity of all the coefficients implies that $Y_i(z)$ is larger than that monomial. Dividing by~$Y_i$ on both sides and letting~$z$ tend to~$\rho$ thus implies that~$Y_{j_\ell}$ has a finite limit and therefore all of them do by the previous part of the proof.
\end{proof}
Note that the last statement of the lemma is not an equivalence. It is possible for a linear system to have generating function solutions with a finite limit at their singularity.

\begin{example}Consider the equation $y(z)=1+B(z)y(z)$ with $B(z)=(1-\sqrt{1-4z})/2$ the generating function of the Catalan numbers. Thus, $y$ is the generating function of sequences of binary trees. This equation satisfies the hypothesis of the lemma: it is well-founded since $B(0)=0$; the system containing only one equation, it is irreducible. This equation is linear in~$y(z)$, but at the singularity $\rho=1/4$, the generating function $1/(1-B(z))$ has finite value $1/(1-B(1/4))=2$.
\end{example}

\subsection{Ordering Solutions}
The following is a variant of~\cite[Lemma~13]{BellBurrisYeats2010}.
\begin{lemma}\label{lemma:dominance}
Let $\bs Y=\bs H(x,\bs Y)$ be an irreducible system and let~$\rho$ be the radius of convergence of its generating function solution~$\bs Y$. Let $\tau_i=\lim_{u\rightarrow\rho-}Y_i(u)$.  Let $(a,\bs B)\in\mathbb{R}_{\ge0}^{m+1}$ belong to the domain of convergence of $\bs H$ or its boundary and satisfy 
$\bs B=\bs H(a,\bs B)$.
\begin{enumerate}
    \item\label{lemma:part1} If for some index~$i$, $B_i\le \tau_i$, then $a\le\rho$;
    \item\label{lemma:part2} if $a\le\rho$, then $\bs Y(a)\le\bs B$;
    \item\label{lemma:part3} if $\bs B\le\bs\tau$ then $\bs Y(a)=\bs B$;
    \item\label{lemma:part4} if $\bs B=\bs\tau$ then $a=\rho$.
\end{enumerate}
where inequalities between vectors mean that they hold for all coordinates.
\end{lemma}
\begin{proof}
The idea of the proof is taken from the proof of Lemma~13 by Bell, Burris and Yeats~\cite{BellBurrisYeats2010}.

We first prove that $\bs B\ge \bs Y(0)$. Let $\bs H^{[m]}$ be the $m$th iterate of $\bs Y\mapsto\bs H(x,\bs Y)$ as in the proof of \cref{lemma:NoZeroInJac}.
By \cref{prop:well-founded-analytic},  $\bs Y(0)=\bs H^{[m]}(0,0)$. Up to relabeling the coordinates, assume that $B_j<Y_j(0)$ for $j=1,\dots,k$ and $Y_j(0)\le B_j$ otherwise. Since $\bpartial\bs H/\bpartial\bs Y$ is nilpotent at~$(0,\bs Y(0))$ (by definition of a zero-free well-founded system) and $\bs H^{[i]}(0,\bs Y(0))=\bs Y(0)$ for all $i>0$, it follows by induction that $\bpartial\bs H^{[m]}/\bpartial\bs Y(0,\bs Y(0))=\bs0$. Then by positivity of $Y_j(0)$ for $j\le k$, this implies that $\partial\bs H^{[m]}/\partial Y_j(0,\bs Y)$ is $0$ for those values of~$j$, i.e., at $x=0$, $\bs H^{[m]}$ is constant wrt $Y_1,\dots,Y_k$. Thus,
\begin{multline*}
B_1<Y_1(0)=H^{[m]}_1(0,Y_1(0),\dots,Y_m(0))\\
=H^{[m]}_1(0,B_1,\dots,B_k,Y_{k+1}(0),\dots,Y_m(0))\le H^{[m]}_1(0,\bs B)\le H^{[m]}_1(a,\bs B)=B_1,
\end{multline*}
a contradiction, which shows that $k=0$ and thus $\bs B\ge \bs Y(0)$.

Next, up to relabeling the coordinates, assume $Y_j(0)\le B_j\le\tau_j$ for $j=1,\dots,k$ and $B_j>\tau_j$ for $j>k$ ($k\le m$). If $k=0$, there is nothing to prove in Part~\ref{lemma:part1} of the lemma. Otherwise, by continuity, there exist $u_1,\dots,u_k$, all bounded by $\rho$, such that $B_j=Y_j(u_j)$ for $j\le k$. Without loss of generality, they can be ordered as $u_1\le \dots\le u_k\le \rho$. 
Since the $Y_j$'s are increasing, $Y_j(u_1)\le B_j$ for $j\le k$ and the same inequality holds for the other ones by definition of $k$. Let $\bs K$ be as in \cref{eq:defK}, which satisfies $\partial K_i/\partial z\neq0$ for all $i$. 
If $u_1<a$, then 
\[B_1=Y_1(u_1)=K_1(u_1,Y_1(u_1),\dots,Y_m(u_1))<K_1(a,B_1,\dots,B_m)=B_1,\]
a contradiction, which implies that $a\le u_1\le\dots\le \rho$ and proves \cref{lemma:part1} of the \namecref{lemma:dominance}.

Next, when $a\leq \rho$, for all $i$ and even if $k=0$,
\[Y_i(a)=H_i(a,Y_1(a),\dots,Y_m(a))\le H_i(a,Y_1(u_1),\dots,Y_k(u_k),B_{k+1},\dots,B_m)=B_i,\]
proving \cref{lemma:part2}.

The third part of the lemma is the case when~$k=m$. Then 
\[B_m=K_m(u_m,Y_1(u_m),\dots,Y_m(u_m))\ge K_m(a,B_1,\dots,B_m)=B_m\]
forces the inequality to be an equality and thus, by irreducibility, monotonicity and dependency in $z$, all the arguments must be equal: $B_i=Y_i(u_m)$ for all $i$ and $u_m=a$.
The fourth part of the lemma comes from the dependency in~$z$ in this inequality.
\end{proof}

Part~1 of Theorem~\ref{thm:Lambda} requires slightly more work.
\begin{corollary} Assume that $\bs Y=\bs A(z)+\bs J(z)\bs Y$ is an irreducible linear system. Let $(a,\bs B)\in\mathbb{R}_{\ge0}^{m+1}$ belong to the domain of convergence of $\bs H=\bs A+\bs J\bs Y$ or its boundary and satisfy $\bs B=\bs H(a,\bs B)$. Then, $a\le\rho$ and $\bs B=\bs Y(a)$.
\end{corollary}
\begin{proof} 
If $\bs Y$ tends to infinity as $x\rightarrow\rho$, then all its coordinates do by \cref{lemma:finiteornot} and the conclusion follows from \cref{lemma:part1} and \cref{lemma:part3} of the previous lemma.

If~$\bs\tau=\lim_{x\rightarrow\rho}\bs Y$ is finite, the point $(\rho,\bs\tau)$ must belong to the boundary of the domain of convergence of~$\bs H$, which implies that $a\le \rho$.
Assume by contradiction that $\det(\Id-\bs J(z))$ tends to~0 as $z\rightarrow\rho$, then from $\det(\Id-\bs J(z))\det((\Id-\bs J(z))^{-1})=1$ it follows that the second determinant tends to infinity and therefore at least one entry of the inverse matrix tends to infinity, which shows that 
$\bs Y(z)=(\Id-\bs J(z))^{-1}\bs A(z)$ has an infinite limit as $z\rightarrow\rho$, a contradiction. Thus in the case when $\bs Y$ has a finite limit as~$z\rightarrow\rho$, the matrix $\Id-\bs J(z)$ is invertible for $z\le\rho$. Subtracting $\bs B=\bs A(a)+\bs J(a)\bs B$ from $\bs Y(a)=\bs A(a)+\bs J(a)\bs Y(a)$ gives $(\Id-\bs J(a))(\bs Y(a)-\bs B)=0$ from which the invertibility of $\Id-\bs J$ implies $\bs B=\bs Y(a)$.
\end{proof}

\subsection{Dominant Eigenvalue}

We now turn to the rest of the theorem, starting with part~\ref{part2b} ($\Leftarrow$) for $a<\rho$.
\begin{lemma}\label{lemma:Lambdalessthan1}
Let $\bs Y=\bs H(z,\bs Y)$ be an irreducible  system and let $\rho$ be the radius of convergence of its generating function solution~$\bs Y$. If $0<a<\rho$, then $0\le\Lambda_{\bs H}(a)<1$.
\end{lemma}
This is \cite[Lemma~10 (b)]{BellBurrisYeats2010} for their {\em well-conditioned} systems. The proof adapts to our setting.
\begin{proof}
Since $(a,\bs Y(a))$ lies inside the domain of convergence of~$\bs H$, $\bs H$ and its partial derivatives are well defined and have nonnegative values at $(a,\bs Y(a))$. Differentiation of $\bs Y=\bs H(z,\bs Y)$ gives
\[\bs Y'(a)=\frac{\partial\bs H}{\partial z}(a,\bs Y(a))+\frac{\bpartial\bs H}{\bpartial\bs Y}(a,\bs Y(a))\cdot\bs  Y'(a).\]
By the Perron-Frobenius theorem, there exists a left eigenvector $v$ of $\frac{\bpartial\bs H}{\bpartial\bs Y}(a,\bs Y(a))$ with positive (and not merely nonnegative)
entries for the eigenvalue $\Lambda_{\bs H}(a)$. 
Multiplying on the left by this vector~$v$ leads to an equation over real values
\[\left(1-\Lambda_{\bs H}(a)\right)v\cdot\bs Y'(a)=v\cdot \frac{\partial\bs H}{\partial z}(a,\bs Y(a)).\]
Since $\partial\bs H/\partial z\neq \bs 0$, its value at a positive $(a,\bs Y(a))$ is nonzero, implying that the right-hand side is positive. This in turn implies that $\bs Y'(a)\neq0$ and finally $\Lambda_{\bs H}(a)<1$.
\end{proof}

\begin{corollary}\label{coro:lambda_at_most_1} Let $\rho$ be the radius of convergence of the generating function~$\bs Y$, solution of the irreducible system $\bs Y=\bs{H}(z,\bs Y)$. Let $\bs\tau=\lim_{u\rightarrow\rho-}\bs Y(u)$. Then $0\le\Lambda_{\bs H}(\rho,\bs\tau)\le 1$.
\end{corollary}
\begin{proof}
This follows from the previous lemma and the continuity of both $\bs Y$ and $\Lambda_{\bs H}$.
\end{proof}
\noindent Part~\ref{part2b} ($\Leftarrow$) of the theorem for $a=\rho$ then follows from the implicit function theorem (\cref{sec:ift}).
The other direction of Part~\ref{part2b} is an adaptation of a similar result of Bell, Burris, Yeats~\cite[Prop.~19]{BellBurrisYeats2010}.

\begin{lemma} Let $\bs Y=\bs H(x,\bs Y)$ be an irreducible system and let $\rho$ be the radius of convergence of its generating function solution~$\bs Y$.
Let $(a,\bs B)\in\mathbb{R}_{\ge0}^{m+1}$ belong to the domain of convergence of $\bs H$ or its boundary and satisfy $\bs B=\bs H(a,\bs B)$. If $\Lambda_{\bs H}(a,\bs B)\le 1$ then $a\le\rho$ and $\bs B=\bs Y(a)$.
\end{lemma}
\begin{proof}
First, if $a>\rho$ then by Lemma~\ref{lemma:dominance}(\ref{lemma:part1}), $\bs B>\bs\tau$. Since $(a,\bs B)$ belongs to the domain of convergence of~$\bs H$, this implies that so does $(\rho,\bs\tau)$. Then by monotonicity of $\Lambda_{\bs H}$ this would imply $\Lambda_{\bs H}(a,\bs B)>\Lambda_{\bs H}(\rho,\bs\tau)=1$, a contradiction.

Thus $a\le\rho$ and then, by Lemma~\ref{lemma:dominance}(\ref{lemma:part2}), $\bs B\ge\bs Y(a)$. 
If $\bs B\neq\bs Y(a)$, Bell, Burris and Yeats cleverly construct a nonnegative matrix~$\bs M$ with spectral radius $r(\bs M)$ such that $r(\bs M)\ge1$, while $\bs M<\bpartial\bs H/\bpartial\bs Y(a,\bs Y(a))$, which implies $1\le r(\bs M) < r(\bpartial\bs H/\bpartial\bs Y(a,\bs Y(a)))=1$, a contradiction. The construction is as follows: for each~$i$, let $\phi_i(t)=H_i(a,\bs Y(a)+t(\bs B-\bs Y(a)))$, so that $\phi_i(0)=Y_i(a)$ and $\phi_i(1)=B_i$. By the mean value theorem, $\phi_i(1)=\phi_i(0)+\phi_i'(t_i)$ for some~$t_i\in(0,1)$. This rewrites as
\[B_i=Y_i(a)+h_i\cdot(\bs B-\bs Y(a)),\]
where $h_i$ is the row vector $\bpartial H_i/\bpartial\bs Y(a,\bs Y(a)+t_i(\bs B-\bs Y(a)))$. The matrix $\bs M$ is chosen as the matrix with rows~$h_1,\dots,h_m$. By design, it satisfies~$\bs B-\bs Y(a)=\bs M\cdot(\bs B -\bs Y(a))$ and $\bs M<\bpartial\bs H/\bpartial\bs Y(a,\bs Y(a))$, whence the result.
\end{proof}
Finally, if $\Lambda_{\bs H}(a,\bs B)=1$, then this lemma implies that $\bs B=\bs Y(a)$ and $a\le \rho$, while \cref{lemma:Lambdalessthan1} shows that $a<\rho\Rightarrow\Lambda_{\bs H}(a,\bs B)<1$, so that $a=\rho$, proving Part~\ref{part2a} of \cref{thm:Lambda} and concluding its proof.

\section{Proof of Proposition \ref{thm:alerho}: Nonnegative Solutions}\label{appendix:proofprop3.11}
\alerho*
\begin{proof}
The proof is by induction on the decomposition of the system into irreducible components. If the system consists of a nonrecursive equation $Y=H(z)$, then by hypothesis the point $z=a$ is in the domain of convergence of~$H$, which is that of~$Y$. The only possible value for $\bs B$ is then $H(a)=Y(a)$.

If the system is irreducible and linear, then the conclusion comes from part~1 of \cref{thm:Lambda}.

If the system is irreducible and nonlinear, let $\bs\tau$ denote $\bs Y(\rho)$, which has finite coordinates, by \cref{lemma:finiteornot}.
If for some~$i$, $B_i\le \tau_i$, then $a\le\rho$ by \cref{lemma:dominance}(\ref{lemma:part1}). 
Otherwise, $\bs B>\bs \tau$. If $(\rho,\bs\tau)$ is not on the boundary of the domain of convergence of~$\bs H$, then the implicit function theorem implies that $\Lambda_{\bs H}(\rho,\bs\tau)\ge 1$, while part 2(a) of \cref{thm:Lambda} implies that $\Lambda_{\bs H}(\rho,\bs\tau)\le1$. Thus either $(\rho,\bs\tau)$ is on the boundary of the domain of convergence of~$\bs H$ or $\Lambda_{\bs H}(\rho,\bs\tau)=1$. The first case cannot occur with $a>\rho$: $(a,\bs B)>(\rho,\bs\tau)$ is beyond the boundary of the domain of convergence. Finally, it is not possible that $(a,\bs B)>(\rho,\bs \tau)$ and $\Lambda_{\bs H}(\rho,\bs\tau)=1$. This is obtained by contradiction, following the proof of \cite[Prop.~11]{BellBurrisYeats2010}: by Taylor expansion and positivity
\[\bs B-\bs\tau=\bs H(a,\bs B)-\bs H(\rho,\bs\tau)
\ge\frac{\bpartial\bs H}{\bpartial\bs Y}(\rho,\bs\tau)(\bs B-\bs\tau)
\]
and the inequality is strict for at least one coordinate, by nonlinearity.
By the Perron-Frobenius theory, there exists a left eigenvector~$\bs v$ with positive coordinates of $\bpartial\bs H/\bpartial\bs Y(\rho,\bs\tau)$ for the dominant eigenvalue $\Lambda_{\bs H}(\rho,\bs\tau)>0$. Multiplying by $\bs v$ on the left gives
\[\bs v(\bs B-\bs\tau)\ge\Lambda_{\bs H}(\rho,\bs\tau)\bs v(\bs B-\bs\tau).\]
Since the previous inequality is strict for at least one coordinate, it follows that $\Lambda_{\bs H}(\rho,\bs\tau)<1$,
a contradiction. Thus for irreducible systems $a\le\rho$ and \cref{lemma:dominance}(\ref{lemma:part2}) shows that $\bs Y(a)\le\bs B$.

Finally, if the system is not irreducible, it splits into two sub-systems
\[\bs Y_{1:i-1}=\bs H_{1:i-1}(z,\bs Y_{1:i-1}),\quad\bs Y_{i:m}=\bs H_{i:m}(z,\bs Y_{1:i-1}(z),\bs Y_{i:m})=:\hat{\bs H}_{i:m}(z,\bs Y_{i:m}),\]
where the second is an irreducible component. Since $(a,\bs B)$ belongs to the domain of convergence of~$\bs H$ or its boundary, $(a,\bs B_{1:i-1})$ belongs to the domain of convergence of $\bs H_{1:i-1}$ and thus by induction, $a$ is at most the radius of convergence of~$\bs Y_{1:i-1}$ and moreover $\bs Y_{1:i-1}(a)\le \bs B_{1:i-1}$. 

By~\cref{lemma:wf-irreducible-components}, $\bs Y_{i:m}=\hat{\bs H}_{i:m}(z,\bs Y_{i:m})$ is an irreducible system. 
Since $(a,\bs B)$ belongs to the domain of convergence of $\bs H$ or its boundary,  so does $(a,\bs Y_{1:i-1}(a),\bs B_{i:m})\le(a,\bs B)$ and therefore $(a,\bs B_{i:m})$ belongs to the domain of convergence of $\hat{\bs H}_{i:m}$ or its boundary. The conditions of the first part of the proof are met and therefore $a\le\rho$ and $\bs Y_{i:m}(a)\le \bs B_{i:m}$, which concludes the proof.
\end{proof}

\section{Proof of Proposition \ref{prop:DLW}: Drmota-Lalley-Woods Theorem}
\label{appendix:proofpropDLW}

\DLW*
\begin{proof}
The proof is an adaptation of Drmota's~\cite{Drmota1997}.

\paragraph{Singularity different from $\rho$}
When $\sigma\neq\rho$, by \cref{lemma:wf-irreducible-components,thm:BBY72i,thm:periods}, the singularity $\sigma$ is of the form $\rho\exp(2i\pi k/q)$, with $k\in\{0,\dots,q-1\}$ and $q$ a common period of the coordinates of the solution. 
Moreover, each coordinate of the generating function $\bs Y$ is of the form $Y_j(z)=z^{\val_0(Y_j)}V_j(z^q)$ with~$V_j(z)$ a function that is analytic in~$|z|<\rho^{1/q}$. As $z\rightarrow\sigma$, $z^q\rightarrow\rho$ so that the expansion of $Y_j(z)$ as $z\rightarrow\sigma$ is obtained by multiply its expansion as~$z\rightarrow\rho$ found in the next part of the proof by the expansion of $z^{\val_0(Y_j)}$.

\paragraph{Change of variables}
By \cref{lemma:finiteornot}, all $Y_i$ have a finite limit at~$\rho$, which is a real positive number since the system is irreducible and therefore zero-free. All the~$U_i$ such that $\partial \bs H/\partial U_i\neq0$ also have a finite limit, as they appear polynomially in $\bs H$. The value of~$\bs Y(\rho)$ is defined as the solution of~$\bs Y(\rho)=\bs H(\rho,\bs Y(\rho),\bs U(\rho))$. 
Let $z=\rho(1-t^r)$ and $\bs Y(z)=\bs Y(\rho)-\tilde{\bs Y}(t)$ and define $\tilde{\bs U}(t)$ similarly,
so that~$\tilde{\bs U}$ has coordinates that are \emph{convergent power series} in $\mathbb C[[t]]$.
The functions~$Y_i$ and~$U_i$ being increasing for~$0<z<\rho$, the leading coefficients in the expansions of~$\tilde{Y}_i$ and $\tilde{U}_i$ are positive.
The function~$\tilde{\bs Y}$ is solution to the equation
\begin{equation}\label{eq:Ytilde}
\tilde{\bs Y}(t)=\bs Y(\rho)-\bs H\!\left(\rho(1-t^r),\bs Y(\rho)-\tilde{\bs Y}(t),\bs U(\rho)-\tilde{\bs U}(t)\right)=:\tilde{\bs H}(t,\tilde{\bs Y}(t)),
\end{equation}
with 
\[
\tilde{\bs H}(t,\bs V)=\bs Y(\rho)-\bs H\!\left(\rho(1-t^r),\bs Y(\rho)-\bs V,\bs U(\rho)-\tilde{\bs U}(t)\right).
\]
By design, both sides of \cref{eq:Ytilde} are~0 at~$t=0$ and 
\[
\bpartial\tilde{\bs H}/\bpartial{\bs V}(t,\tilde{\bs Y}(t))=\bpartial{\bs H}/\bpartial{\bs Y}(z,\bs Y(z),\bs U(z))|_{z=\rho(1-t^r)}.
\] 

\paragraph{Regular case}
If $\Lambda_{\bs H}(\rho)<1$, the matrix $\Id-\bpartial\bs H/\bpartial\bs Y$ is invertible at~$(\rho,\bs Y(\rho),\bs U(\rho))$. Then, so is $\Id-\bpartial\tilde{\bs H}/\bpartial{\bs V}$ at~$(0,\bs 0)$ and therefore \cref{eq:Ytilde} satisfies the hypotheses of the implicit function theorem over convergent power series. It follows that $\tilde{\bs Y}(t)$ has convergent coordinates in $\mathbb C\{t\}$, so that those of~$\bs Y$ lie in $\mathbb C\{(1-z/\rho)^{1/r}\}$, which proves the theorem in that case. 

\paragraph{Halving of the exponent in the singular case}
If $\Lambda_{\bs H}(\rho)=1$, the matrix $\Id-\bpartial\bs H/\bpartial\bs Y$ is not invertible at~$(\rho,\bs Y(\rho),\bs U(\rho))$. Let $m$ be the dimension of the system. 
The simplest case is when~$m=1$. Then we have 
\[
\tilde H(0,0)=0,\quad\frac{\partial\tilde H}{\partial V}(0,0)=\frac{\partial H}{\partial Y}(\rho,Y(\rho),U(\rho))=1,\quad\frac{\partial^2\tilde H}{\partial V^2}(0,0)=-\frac{\partial^2H}{\partial Y^2}(\rho,Y(\rho),U(\rho))\neq0,\]
the last one by nonlinearity and positivity of the power series. Then the Weierstrass preparation theorem~\cite[Cor.~6.1.2]{Hormander1990} implies that the solution~$\tilde{Y}$ of the equation~$\tilde Y=\tilde H(t,\tilde Y)$ has a convergent Puiseux expansion in~$\mathbb C\{t^{1/2}\}$, giving the theorem in that case.

The case when $m>1$ reduces to the case when~$m=1$ by an argument due to Drmota~\cite[\S5]{Drmota1997},\cite[p.~64-65]{Drmota2009}, see also Mac~Millan~\cite{Mac-Millan1912} for an ancestor of this idea.
Write $\hat{\bs Y}$ for $(\tilde Y_2,\dots,\tilde Y_m)$ and similarly for~$\hat{\bs H}$. A consequence of the Perron-Frobenius theory is that the first minors of $\Id-\bpartial\bs H/\bpartial\bs Y$ are nonzero~\cite[XIII \S2.5]{Gantmacher1959}; in particular, $\Id-\bpartial\hat{\bs H}/\bpartial\hat{\bs Y}$ is invertible at $(0,\bs 0)$. 
Then the implicit function theorem can be applied to solve the system $\hat{\bs Y}=\hat{\bs H}(t,Y,\hat{\bs Y})$ for unique $\hat{\bs Y}$ that are convergent power series in~$Y$ and~$t$. Moreover, differentiation gives 
\[
\frac{\partial\hat{\bs Y}}{\partial Y}=\left(\Id-\frac{\bpartial\hat{\bs H}}{\bpartial\hat{\bs Y}}\right)^{-1}\frac{\partial\hat{\bs H}}{\partial Y}.\]
Injecting these solutions~$\hat{\bs Y}$ into the first equation~$\tilde Y_1(t)=\tilde H_1(t,\tilde Y_1(t),\hat{\bs Y})$ gives a single equation in the unknown~$\tilde Y_1$,
\[\tilde Y_1=G(t,\tilde Y_1),\]
with
$G(t,Y)=\tilde H_1(t,Y,\hat{\bs Y}(t,Y))$,
a convergent power series in $(t,Y)$ with nonnegative coefficients and $G(0,0)=0$. The situation is reduced to the case when $m=1$, and it is now sufficient to show that $\partial G/\partial Y(0,0)=1$ and $\partial^2 G/\partial Y^2(0,0)\neq0$.
Differentiation gives
\[\frac{\partial G}{\partial Y}=
\frac{\partial \tilde H_1}{\partial Y_1}(t,Y,\hat{\bs Y}(t,Y))
+\left(\frac{\bpartial \tilde H_1}{\bpartial\hat{\bs Y}}(t,Y,\hat{\bs Y}(t,Y))\right)^T\frac{\partial\hat{\bs Y}}{\partial Y}.
\]
Since the matrix $\Id-\bpartial\tilde{\bs H}/\bpartial \bs Y(0,\bs0)$ is singular, it follows that $1-\partial G/\partial Y(0,0)=0$:
the Schur complement of the invertible block $\Id-\bpartial\hat{\bs H}/\bpartial\hat{\bs Y}(0,\bs 0)$ in this singular matrix
 is 
 \begin{equation}\label{eq:Schur}
0=\left(1-
\frac{\partial \tilde H_1}{\partial Y_1}
-\left(\frac{\bpartial\tilde{H}_1}{\bpartial\hat{\bs Y}}\right)^{T}
\left(\Id-\frac{\bpartial\hat{\bs H}}{\bpartial\hat{\bs Y}}\right)^{-1}
\frac{\partial\hat{\bs H}}{\partial{Y}}\right)(0,\bs 0)=1-\frac{\partial G}{\partial Y}(0,\bs 0).
\end{equation}
Positivity of all the power series involved implies that~$\partial^2 G/\partial Y^2\neq0$. The conclusion for $m=1$ is that the solutions~$\tilde{Y}_1$ of this equation have convergent Puiseux expansions in~$\mathbb C\{t^{1/2}\}$. Therefore, so do the corresponding~$\tilde Y_2,\dots$ by composition. Thus all coordinates~$Y_i$ belong to~$\mathbb C\{Z^{r/2}\}$.

\paragraph{Appearance of the square-root.}
The Taylor expansion of~$\tilde{\bs H}(t,\bs V)$ with respect to~$\bs V$, at~$(0,\bs 0)$, evaluated at~$\bs V=\tilde{\bs Y}(t)$, gives
\begin{equation}\label{eq:expansion-tilde-Y}
\tilde{\bs Y}(t)=
\tilde{\bs H}(t,\bs 0)
+\frac{\bpartial\tilde{\bs H}}{\bpartial\bs V}(0,\bs 0)\cdot \tilde{\bs Y}(t)
+\frac12\frac{\bpartial^2\tilde{\bs H}}{\bpartial\bs V^2}(0,\bs 0)(\tilde{\bs Y}(t),\tilde{\bs Y}(t))+\text{s.o.t.}.
\end{equation}
By \cref{lemma:NoZeroInJac}, we can assume that all entries of $\bs\partial\tilde{\bs H}/\bs\partial\bs V(0,\bs 0)$ are positive. Thus the leading coefficients of each of the coordinates of the first two terms in the right-hand side are nonnegative and the valuation of each coordinate of the right-hand side comes from these two terms. Therefore the valuation of each coordinate~$\tilde Y_i$ is at most that of~$\tilde{H}_i(t,0)$. 
Using \cref{lemma:NoZeroInJac} if necessary, we can assume that all the entries of~${\bpartial\tilde{\bs H}}/{\bpartial\bs V}(0,\bs 0) = \bpartial{\bs H}/\bpartial{\bs Y}(\rho,\bs Y(\rho),\bs U(\rho))$ are nonzero, which implies that, on the right-hand side, the minimal valuation among the~$\tilde Y_i$ appears  for each coordinate, and therefore also on the left, showing that all these valuations are identical.
Let~$\alpha$ be their common value. 

The next step shows that~$\alpha=1/2$. Let~$\bs v$ be a \emph{left} eigenvector of~$\bpartial\tilde{\bs H}/\bpartial\bs V(0,\bs 0)$ for the eigenvalue~1.  It can be taken with positive coordinates by the Perron-Frobenius theorem, so that no cancellation occurs when
multiplying \cref{eq:expansion-tilde-Y}  on the left by~$\bs v$, giving
\begin{equation}\label{eq:left-by-v}
\bs v\cdot\frac{\partial\tilde{\bs H}}{\partial t}(0,\bs 0)t=- \frac12\bs v\cdot\frac{\bpartial^2\tilde{\bs H}}{\bpartial\bs V^2}(0,\bs 0)(\tilde{\bs Y}(t),\tilde{\bs Y}(t))+\text{s.o.t.}.
\end{equation}
By positivity, $\bpartial^2\tilde{\bs H}/\bpartial\bs V^2(0,\bs 0)\neq0$, so that the right-hand side has valuation~$2\alpha$. The linear term of the left-hand side with respect to~$t$ comes from
\begin{equation}\label{eq:dHt/dt00}
\bs W:=\frac{\partial\tilde{\bs H}}{\partial t}(0,\bs 0)=\rho r [t^{r-1}]_{t=0}\frac{\partial\bs H}{\partial z}(\rho,\bs Y(\rho),\bs U(\rho))
+\frac{\bs\partial\bs H}{\bs\partial\bs U}(\rho,\bs Y(\rho),\bs U(\rho))\tilde{\bs U}'(0).
\end{equation}
There, $\tilde{\bs U}'(0)=-\bs T$ (with $\bs T$ from the hypothesis) has at least one nonzero coordinate, that is positive, and the corresponding~$\partial \bs H/\partial U_i$ is nonzero and positive too, so that $\bs W$ is nonzero.
Thus the valuation $\alpha$ of all coordinates of~$\tilde{\bs Y}$ is $1/2$. 

\paragraph{Coefficients~$C_i$.}
Let finally $\bs C$ be such that $\tilde{\bs Y}(t)\sim -\bs Ct^{1/2}$. Extracting the terms with exponent~1/2 in \cref{eq:expansion-tilde-Y} gives
\[-\bs C=-\frac{\bpartial\tilde{\bs H}}{\bpartial\bs V}(0,\bs 0)\cdot \bs C,\]
i.e., $\bs C$ is a \emph{right} eigenvector of~$\bpartial\tilde{\bs H}/\bpartial{\bs V}(0,\bs 0)= \bpartial{\bs H}/\bpartial{\bs Y}(\rho,\bs Y(\rho),\bs U(\rho))$ for the eigenvalue~1. By the Perron-Frobenius theorem, there is a one-dimensional choice of such vectors. Let~$\bs w$ be a right eigenvector with positive coordinates, and let~$\lambda>0$ be such that~$\bs C=\lambda \bs w$. (The sign of~$\lambda$ is dictated by the fact that $\bs Y(z)$ is increasing for~$z<\rho$.) The constant factor~$\lambda$ is obtained by considering the terms with exponent~1 in \cref{eq:left-by-v}, giving the quadratic equation
\[
\bs v\cdot\bs W
=- \lambda^2\frac12\bs v\cdot\frac{\bpartial^2\tilde{\bs H}}{\bpartial\bs V^2}(0,\bs 0)({\bs w},{\bs w})
=\lambda^2\frac12\bs v\cdot\frac{\bpartial^2{\bs H}}{\bpartial\bs Y^2}(\rho,\bs Y(\rho),\bs U(\rho))({\bs w},{\bs w})
,
\]
which concludes the proof of the proposition, in view of \cref{eq:dHt/dt00}.
\end{proof}

\anonymousvariant{
\section*{Acknowledgment} This work has been supported in part by the 
French project NuSCAP (ANR-20-CE48-0014).
}{}

\bibliographystyle{plain}
\bibliography{biblio}

\end{document}